\documentclass[11pt]{article}

\usepackage[margin=1in]{geometry}
\usepackage{amsmath,amssymb,amsthm}
\usepackage{graphicx}
\usepackage{booktabs,longtable}
\usepackage{array}
\usepackage{placeins}
\usepackage{natbib}
\usepackage{xr-hyper}
\usepackage{hyperref}

\hypersetup{
  colorlinks=true,
  linkcolor=blue,
  citecolor=blue,
  urlcolor=blue,
  pdftitle={Elimination Geometry II: Integrability, Representation, and Structural Limits},
  pdfauthor={Mian Huang and Xueqin Wang}
}

\newtheorem{theorem}{Theorem}[section]
\newtheorem{proposition}[theorem]{Proposition}
\newtheorem{corollary}[theorem]{Corollary}

\newtheorem{definition}[theorem]{Definition}

\newtheorem{remark}[theorem]{Remark}
\newtheorem{example}[theorem]{Example}
\numberwithin{equation}{section}

\newcommand{\E}{\mathbb E}
\newcommand{\R}{\mathbb R}
\newcommand{\cA}{\mathcal A}
\newcommand{\cC}{\mathcal C}
\newcommand{\cD}{\mathcal D}

\newcommand{\cH}{\mathcal H}
\newcommand{\cJ}{\mathcal J}
\newcommand{\cM}{\mathcal M}

\newcommand{\cR}{\mathcal R}

\newcommand{\argmin}{\operatorname*{arg\,min}}

\newcommand{\KL}{\operatorname{KL}}
\newcommand{\TC}{\operatorname{TC}}
\newcommand{\diag}{\operatorname{diag}}
\newcommand{\tr}{\operatorname{tr}}
\newcolumntype{L}[1]{>{\raggedright\arraybackslash}p{#1}}

\title{\Large Elimination Geometry II: Integrability, Representation, and Structural~Limits}

\author{Mian Huang and Xueqin Wang}

\date{\today}

\begin{document}

\maketitle

\begin{abstract}
Elimination geometry studies what remains after an auxiliary variable is
optimized out: an exact, nonnegative defect that prices every nonoptimal
auxiliary choice in the units of the original objective.
Part I introduces a constructive elimination--coupling calculus that takes a
declared certified lift as its interface.  This Part II audits that
interface, taking a certified elimination
\(\cJ(\Theta)=\inf_A\cH_0(\Theta,A)\) and its nonnegative defect
\(\cD_\Theta(A)=\cH_0(\Theta,A)-\cJ(\Theta)\) as primitive, in three
theorem families.

First, integrability.  Locally reported gaps form a repeated
partial-elimination tower exactly when they are nonnegative, touching, and
have zero circulation on every elementary elimination square; approximate
curvature bounds elimination-order ambiguity.  In finite spaces, Hodge
projection and a witness-constrained quadratic program construct a nearest
certified repair and split its cost exactly into curvature and
certification taxes, plus a globally supported period tax on a general
complex that is invisible to every local face test.

Second, representation and price.  A gate handles invertible coordinates
and surjective refinements, the latter being valid exactly when coherence
descends by fiberwise minimization; consequently no nonconstant strictly
convex separable raw-weight penalty survives duplication of
output-equivalent components.  Unrestricted lift design is vacuous: dummy
lifts make coherence free and target-calling lifts manufacture any local
positive-semidefinite metric.  Conditional on an admissible lift, the
fidelity--coherence frontier
\(\mathfrak F_\Theta(r)=\inf\{\cD_\Theta(A):\cR(A)\le r\}\)
separates irreducible coherence cost from external-design inefficiency.  An
augmented modal Gram block is a complete fingerprint of the unconstrained
quadratic lift--coupling germ under origin-preserving linear equivalence,
placing target-shift value, gradient, and transmitted outer curvature in one
positive-semidefinite exchange account.  Its contraction filter yields the
aligned Gaussian risk ledger and, for smoothers commuting with the penalty,
an exact reading of external inefficiency,
because every positive-definite self-adjoint contraction in the
unconstrained quadratic class is the variational solution of an induced
penalty.  Algebraic validity is not
semantic validity.  Bernstein spectral functions map every finite weighted
graph Laplacian to a possibly nonlocal Laplacian.  Conversely, a tunable
weighted three-node path gives an exact failure interval whenever a spectral
function grows locally faster than linearly.  Hence
\((I+\lambda L^\nu)^{-1}\) is Markov for every weighted graph and every
\(\lambda>0\) exactly for \(0<\nu\le1\), and the only graph-universal
polynomial spectral penalties are nonnegative multiples of \(L\).  This
bounds linear mechanisms, not probability-valued estimators: Part I's
nonlinear CDF repair escapes it, retains an exact
graph-effective-dimension risk bound, and attains the matching minimax CRPS
rate for every real graph-source smoothness on growing paths.

Third, converses and impossibility.  Target-\(k\)-jet preservation after
freezing a coupled lift forces subtraction of the defect \(k\)-jet, uniquely
modulo a \(k\)-flat term.  The canonical model splices the target's
low-order jet to the lift's higher-order tail; in a conjugate lift an exact
transmission map identifies which auxiliary perturbations reach X curvature,
and an affine natural-parameter map transmits none.  A perturbation theorem for
approximate defect jets yields quantitative score error, quasi-descent, and
stationary limits under summable certification error.  Fixed-marginal
coupling cannot improve additive coordinate targets and obeys a
dependence--benefit frontier.  In
entropy-regularized control, reversing the certified KL direction destroys
every dimension-free performance modulus.  In a noncommutative Gibbs fiber,
measured KL is only a lower certificate; one fixed measurement preserves an
entire faithful family's Umegaki defects exactly if and only if that family
commutes.  A classical Polyakov fiber checks the same ledger on a
degenerate continuum elimination.  All definitions used here are restated,
so Part II can be read
independently.
\end{abstract}


\section{Introduction}
\label{sec:intro}

Part I of this series proposes a constructive calculus for algorithms built
from a declared certified auxiliary lift.  It starts from that lift as an
explicit interface and asks how pointwise elimination, external
globalization, exactification, variational coupling, and fixed-marginal
coupled-law design can be constructed and certified.  The present Part II
asks the reverse structural questions: when the declared local certificates
actually integrate to one objective, when a coupling claim is invariant to
representation, which converses are forced by target preservation, and which
extensions are impossible.

The audit is self-contained.  Elimination hides a geometry: whenever an
objective is obtained by profiling an auxiliary variable,
\begin{equation}
 \cJ(\Theta)=\inf_{A\in\cA}\cH_0(\Theta,A),
 \label{eq:introelimination}
\end{equation}
every nonoptimal \(A\) remains above the same target point with height
\begin{equation}
 \cD_\Theta(A)
 =\cH_0(\Theta,A)-\cJ(\Theta)\ge0.
 \label{eq:introdefect}
\end{equation}
Thus elimination does not destroy the auxiliary coordinate: it quotients
out a fiber whose vertical coordinate is already denominated in exact units
of objective loss.  Coupling several such fibers is therefore not primarily
a smoothing problem.  It is a constrained movement through a defect
geometry.

For orientation, the five mechanism labels used throughout the series can be
defined without appealing to Part I.  P is exact pointwise elimination.  G
replaces the pointwise auxiliary oracle by an externally coupled field and
evaluates the original lift there.  X corrects such a frozen coupled lift to
preserve a declared target jet.  V inserts an auxiliary coherence cost before
elimination and thereby changes the profiled target.  C chooses dependence
inside a joint-law fiber after the coordinate marginals have been fixed.
Thus P/G/X/V locate where deterministic coupling enters, while C is an
orthogonal fixed-marginal operation; these labels describe architectures, not
guarantees.

Four structural questions must then be answered.  When do locally reported
gaps belong to one common eliminated objective?  Which properties survive a
change of auxiliary coordinates?  Which scalar and matrix invariants govern
the price and leverage of coupling?  What exactness remains possible under
approximate correction, projection, or measurement?  These questions concern
the admissibility and limits of the calculus, not the behavior of one
particular algorithm.

Fix a coherence functional \(\cR:\cA\to[0,+\infty]\).  Conditional on a
declared lift, the scalar design object is the lower boundary of attainable
pairs
\[
 \bigl(\cR(A),\cD_\Theta(A)\bigr).
\]
Before that boundary can be compared, however, two gates must be passed:
the local gaps must integrate to a certified scalar objective, and the lift
must be compared only modulo coupling-covariant reparameterization.  The
paper is therefore organized around three theorem families.

\begin{enumerate}
\item \textbf{Certified integrability and repair.}
Section~\ref{sec:flatness} proves that local gaps form one repeated
partial-elimination tower exactly when they are nonnegative, touching, and
flat on every elementary square.  Approximate curvature bounds
elimination-order ambiguity.  In finite spaces, Hodge projection and
certification constraints give a constructive two-tax repair.

\item \textbf{Representation validity, frontier geometry, and leverage
filters.}
Section~\ref{sec:liftgate} identifies when lift--coupling packages represent
the same coupled object, including an iff transport rule for surjective
refinements and a duplication no-go.  Section~\ref{sec:frontier} studies the
resulting fidelity--coherence frontier and external-design inefficiency.
Section~\ref{sec:fingerprint} classifies the local quadratic germ by a
complete augmented modal fingerprint, derives its
value--gradient--curvature exchange block, and shows that an external
positive-definite self-adjoint contraction is the variational solution of
an induced penalty in the unconstrained quadratic class, so that design
inefficiency compares induced and declared filters at matched roughness.
Section~\ref{sec:vsemantics}
then subjects such filters to a semantic V passport.  It identifies the
three-node Dirichlet cone, proves a graph-universal positive passport for
Bernstein spectral functions, and supplies a scalable weighted-path reverse
witness.  The resulting spectral-power iff classification and polynomial
no-go separate probability validity, finite-hop locality, and approximation
order.  Their collision with the companion's real-smoothness
effective-dimension and minimax CRPS theorem separates a genuine linear
representation obstruction from the statistical cost of nonlinear repair.

\item \textbf{Coupling consequences and impossibility tests.}
Sections~\ref{sec:exactification}--\ref{sec:examples} give the converse and
approximate robustness of defect-jet correction, its jet-splicing normal
form, and a conjugate-lift curvature-transmission no-go, followed by an
aligned Gaussian three-currency filter account, the fixed-marginal
dependence--benefit boundary, and three hard tests: two cases in which a
plausible certificate
has the wrong currency or inequality direction, and a classical continuum
fiber on which the exact ledger survives a degenerate elimination.  These results yield a
scoped no-free first-order coupling principle, not a universal scalar
conservation law: native dependence, overlap, and computational costs become
comparable with target loss only after an explicit exchange condition is
proved.
\end{enumerate}

\emph{Elimination Geometry I:
A Certified Calculus of Auxiliary-Field Coupling}
\citep{FrameworkPartI} develops the constructive P/G/X/V/C calculus,
its model instances, and its numerical diagnostics.  It treats a declared
certified lift as an explicit interface assumption and states the conditions
needed by each mechanism where they are used.  This Part II does not supply
missing premises needed to read Part I; it audits that interface by proving
when the certificates integrate, which representation changes are
legitimate, what target preservation forces, and where proposed extensions
fail.  Part I also constructs a nonlinear CDF projection repair for a
higher-order signed graph proposal.
Section~\ref{sec:vsemantics} identifies the obstruction it repairs: the raw
higher-order linear resolvent cannot have a graph-universal probability
passport, while an active nonlinear projection no longer belongs to the
exact linear fingerprint class.  Part I additionally proves that the
projected empirical-CDF estimator retains the raw proposal's
graph-effective-dimension risk bound and attains the matching minimax CRPS
rate for every real source smoothness \(s>0\) on growing paths; that
statistical theorem is not an unstated premise of the structural barrier
proved here.  A two-mode quadratic laboratory is used throughout Part II to show how
the integrability gate, frontier, leverage spectrum, exactification, and
fixed-marginal layer fit together.  Section~\ref{sec:related} separates
classical ingredients, constructions owned by Part I, and the structural
converses or impossibility statements proved here.

\begin{table}[t]
\centering
\small
\caption{Cross-paper proof-dependency ledger.  The structural results in
the middle column are proved without the constructive or statistical
conclusions in the left column.}
\label{tab:crosspaperdependencies}
\begin{tabular}{L{0.29\textwidth}L{0.29\textwidth}L{0.32\textwidth}}
\toprule
Construction owned by Part~I & Result owned by Part~II & Direction of use \\
\midrule
CRPS effective-dimension bound, Assouad lower bound, and nonlinear CDF
repair
& Weighted-path obstruction and spectral-power passport
& Part~II supplies only the linear probability-validity barrier; Part~I's
rate and repair theorem are not premises of that barrier. \\
Ambient induced-penalty identity for the Richardson filter
& Contraction body and graph Dirichlet passport
& Substitution identifies the same raw linear operator; an active CDF
projection lies outside the modal theorem in both papers. \\
Accepted X descent and inexact-X implementation
& Converse uniqueness of defect-jet subtraction
& Part~I is the constructive example; the converse proof is
self-contained and is not used to certify descent. \\
Spectral, Gibbs, and LogDet block updates
& Representation, leverage, and measurement limits
& The limits audit the scope of the updates after the fact and do not
serve as hidden premises. \\
\bottomrule
\end{tabular}
\end{table}

\section{Certified elimination and the defect coordinate}
\label{sec:elimination}

Let \(\Theta\in\mathsf T\) denote an outer parameter field and
\(A\in\cA\) an auxiliary field, possibly indexed by observations, tasks,
spatial sites, temperatures, replicas, or combinations thereof.  We work in a
minimization orientation; likelihood problems use negative log-likelihood.

\begin{definition}[Certified minimum elimination]
\label{def:certified}
A lifted criterion \(\cH_0:\mathsf T\times\cA\to
\R\cup\{+\infty\}\) is a \emph{certified minimum elimination} of
\(\cJ:\mathsf T\to\R\cup\{+\infty\}\) if
\begin{equation}
 \cJ(\Theta)=\inf_{A\in\cA}\cH_0(\Theta,A)
 \label{eq:certifiedelim}
\end{equation}
and the nonnegative gap
\begin{equation}
 \cD_\Theta(A)
 :=\cH_0(\Theta,A)-\cJ(\Theta)
 \label{eq:defectcoordinate}
\end{equation}
is available as an exact identity or as a computable certificate.  Define
the oracle set by
\begin{equation*}
 \cA^\star(\Theta)=\argmin_{A\in\cA}\cH_0(\Theta,A).
\end{equation*}
It is assumed nonempty whenever an oracle is invoked.
\end{definition}

The word ``certified'' excludes an arbitrary residual that merely vanishes
at a selected optimizer.  The scalar in \eqref{eq:defectcoordinate} must be
the actual insertion cost in the lifted objective.  This distinction is what
allows different coupling mechanisms to share one ledger.

\begin{example}[Conjugate-complete Bregman lifts]
\label{ex:bregmanlift}
Suppose
\[
 \cH_0(\Theta,A)
 =\sum_{t,i}\omega_{it}h_{it}(A_{it},\theta_t),
 \qquad \omega_{it}\ge0,
\]
and each inner problem is conjugate complete with oracle
\(a_{it}^\star(\theta_t)\).  Then
\begin{equation}
 \cD_\Theta(A)
 =\sum_{t,i}\omega_{it}
 D_{\Phi_{it}}\!\left(A_{it}\middle\|a_{it}^\star(\theta_t)\right).
 \label{eq:bregmandefect}
\end{equation}
Negative entropy produces relative-entropy defects, Burg entropy produces
LogDet defects on matrix cones, and log-partition duality produces
information-geometric defects for exponential families.  Bregman centroids
and their directed asymmetry are classical \citep{BanerjeeEtAl2005}; the role
of \eqref{eq:bregmandefect} here is to identify which orientation is forced
by the declared elimination.
\end{example}

\begin{remark}[Representation dependence]
\label{rem:representation}
The target \(\cJ\) may admit several certified lifts, and their defect
coordinates need not agree.  The theory below is conditional on the declared
lift.  A change of auxiliary representation is therefore a genuine design
decision, whereas a bijective reparameterization of one fixed lift merely
rewrites the same scalar defect.
\end{remark}

\subsection{A running quadratic laboratory}
\label{subsec:quadraticlab}

One model will be carried through the three theorem families.  Let
\(\Theta\in\R^p\), \(A\in\R^q\), and consider the exactly quadratic lift
\begin{equation}
 \cH_0(\Theta,A)
 =
 \frac12\Theta^\top H_L\Theta
 +\Theta^\top C_LA
 +\frac12A^\top B_LA,
 \qquad B_L\succ0,
 \label{eq:quadraticlab}
\end{equation}
with affine terms suppressed.  Partial elimination gives
\begin{align}
 A^\star(\Theta)&=-B_L^{-1}C_L^\top\Theta,
 \label{eq:quadraticlaboracle}\\
 \cJ(\Theta)
 &=\frac12\Theta^\top
 \left(H_L-C_LB_L^{-1}C_L^\top\right)\Theta,
 \label{eq:quadraticlabtarget}\\
 \cD_\Theta(A)
 &=\frac12\{A-A^\star(\Theta)\}^\top
 B_L\{A-A^\star(\Theta)\}.
 \label{eq:quadraticlabdefect}
\end{align}
Let the coherence cost be
\begin{equation}
 \cR(A)=\frac12A^\top K_LA,\qquad K_L\succeq0.
 \label{eq:quadraticlabroughness}
\end{equation}
When \(A=(A_1,A_2)\), its canonical conditional defects form a flat
two-block elimination square.  Section~\ref{sec:flatness} shows how a
locally touching perturbation of one edge creates detectable curvature.
Section~\ref{sec:liftgate} audits its admissible reparameterizations;
Sections~\ref{sec:frontier}--\ref{sec:fingerprint} derive its frontier,
modal leverage, and intrinsic variational filter:
\[
 A_\lambda(\Theta)=-(B_L+\lambda K_L)^{-1}C_L^\top\Theta.
\]
Section~\ref{sec:exactification}
freezes such a field and audits exact or approximate jet correction.  If
\(B_L\) is block diagonal and the two auxiliary blocks are randomized,
Section~\ref{sec:C} shows that fixed coordinate laws preserve expected
additive defect while cross-covariance may change graph roughness.  Thus the
laboratory is not an additional application; it is a common calculation
showing how the abstract gates compose.

Let \(\cR:\cA\to[0,+\infty]\) be a coherence cost.  Its zero set
\begin{equation}
 \cC:=\{A:\cR(A)=0\}
 \label{eq:coherenceset}
\end{equation}
may encode full consensus, graphwise agreement, a low-dimensional shared
family, or another declared notion of perfect coherence.  The geometry of
\(\cR\) is selected by the analyst; the geometry of \(\cD_\Theta\) is
generated by the lift.

\section{Flatness of the elimination lattice}
\label{sec:flatness}

Definition~\ref{def:certified} begins with a global lifted criterion.  A
multi-stage method is often specified in the opposite direction: it supplies
one local surrogate gap for each block and tacitly assumes that these gaps
belong to a common objective.  Local majorization and touching do not by
themselves justify that assumption.  The missing condition is integrability
over elimination order.

Fix \(\Theta\) and suppress it below.  Let
\(\cA=\prod_{i=1}^K\cA_i\), with every \(\cA_i\) nonempty, and let
\([K]=\{1,\ldots,K\}\).  If a top criterion
\(H:\cA\to\R\) is given, its partial values are
\begin{equation}
 V_S(a_S):=\inf_{a_{S^c}}H(a),
 \qquad S\subseteq[K],
 \label{eq:partialvalues}
\end{equation}
whenever the displayed infima are finite.  Thus
\(V_\varnothing=J:=\inf_aH(a)\) and \(V_{[K]}=H\).  The directed edge
\(S\to S\cup\{j\}\) carries the conditional defect
\begin{equation}
 d_{j\mid S}(a_S,a_j)
 :=V_{S\cup\{j\}}(a_S,a_j)-V_S(a_S).
 \label{eq:edgedefect}
\end{equation}
Commutation of infima gives
\begin{equation}
 d_{j\mid S}\ge0,
 \qquad
 \inf_{a_j}d_{j\mid S}(a_S,a_j)=0.
 \label{eq:edgenontouch}
\end{equation}
The Boolean lattice of retained coordinate sets is therefore equipped with
a nonnegative, touching edge field.

\begin{theorem}[Flat certified-elimination criterion]
\label{thm:flatness}
For every \(S\subseteq[K]\) and \(j\notin S\), suppose a real-valued local
gap
\[
 e_{j\mid S}:\cA_{S\cup\{j\}}\to\R
\]
is given.  The following are equivalent.

\begin{enumerate}
\item There are potentials \(V_S:\cA_S\to\R\), unique after fixing the
constant \(V_\varnothing\), such that
\begin{equation}
 e_{j\mid S}=V_{S\cup\{j\}}-V_S,
 \qquad
 V_S(a_S)=\inf_{a_j}V_{S\cup\{j\}}(a_S,a_j)
 \label{eq:globaltower}
\end{equation}
on every edge.  Consequently, with \(H=V_{[K]}\),
\[
 V_S(a_S)=\inf_{a_{S^c}}H(a)
\]
for all \(S\).

\item The edge field satisfies all three local certificates:
\begin{align}
 e_{j\mid S}(a_S,a_j)&\ge0,
 \label{eq:localnonnegative}\\
 \inf_{a_j}\left|e_{j\mid S}(a_S,a_j)\right|&=0,
 \label{eq:localtouching}
\end{align}
and, for every distinct \(i,j\notin S\),
\begin{align}
 &e_{i\mid S}(a_S,a_i)
 +e_{j\mid S\cup\{i\}}(a_S,a_i,a_j)
 \nonumber\\
 &\quad =
 e_{j\mid S}(a_S,a_j)
 +e_{i\mid S\cup\{j\}}(a_S,a_i,a_j).
 \label{eq:squareflatness}
\end{align}
\end{enumerate}
\end{theorem}

\par\medskip

\begin{example}[A curved report of the quadratic laboratory]
\label{ex:quadraticlabcurvature}
Fix \(\Theta\) in \eqref{eq:quadraticlab}, split
\(A=(A_1,A_2)\), and let \(d_{j\mid S}\) be the canonical conditional
defects generated by partial minimization of \(\cH_0(\Theta,\cdot)\).
They are nonnegative, touching, and flat by
Theorem~\ref{thm:flatness}.  Now leave three edges unchanged and report
\[
 e_{2\mid\{1\}}=(1+\eta)d_{2\mid\{1\}},
 \qquad \eta>0,
\]
on the fourth.  Every reported edge remains nonnegative and touching at
the same conditional oracle, but its signed square circulation becomes
\[
 \eta\,d_{2\mid\{1\}},
\]
for the orientation defined below.  It is nonzero away from the conditional
oracle.  Thus even a harmless-looking rescaling of one locally valid
quadratic gap destroys the existence of one common scalar objective.
The example will reappear below in a different guise: a simultaneous
coordinate transport of the entire lift is valid, whereas changing only
one reported edge is not.
\end{example}

\begin{example}[Locally valid edges without a global objective]
\label{ex:curvededges}
Let \(K=2\), \(\cA_1=\cA_2=\{0,1\}\), and define
\[
 e_{1\mid\varnothing}(a_1)=a_1,\qquad
 e_{2\mid\varnothing}(a_2)=a_2,
\]
\[
 e_{2\mid\{1\}}(a_1,a_2)=\mathbf 1\{a_1\ne a_2\},
 \qquad
 e_{1\mid\{2\}}(a_1,a_2)=0.
\]
Every edge is nonnegative and can touch zero in its newly retained
coordinate.  At \((a_1,a_2)=(1,0)\), however, the two path sums are
\[
 e_{1\mid\varnothing}(1)+e_{2\mid\{1\}}(1,0)=2,
 \qquad
 e_{2\mid\varnothing}(0)+e_{1\mid\{2\}}(1,0)=0.
\]
Thus the square curvature is \(2\), and no scalar \(H(a_1,a_2)\) can have
all four functions as its conditional elimination defects.
\end{example}

Theorem~\ref{thm:flatness} separates three failures that are otherwise easy
to conflate.  A negative edge is a majorization failure.  A positive
edgewise infimum is a touching failure.  A nonzero square circulation is
more global: all individual edges may be valid surrogates, yet their
collection cannot be the increments of one certified scalar objective.

\begin{definition}[Defect curvature]
\label{def:defectcurvature}
For an arbitrary edge field \(e\), define its elementary curvature by
\begin{align}
 \Omega_{ij\mid S}
 :={}&e_{i\mid S}+e_{j\mid S\cup\{i\}}
 \nonumber\\
 &-e_{j\mid S}-e_{i\mid S\cup\{j\}},
 \label{eq:defectcurvature}
\end{align}
with every term evaluated at the compatible restriction of the same
\(a_{S\cup\{i,j\}}\).
\end{definition}

\begin{proposition}[Quantitative elimination holonomy]
\label{prop:holonomy}
For a permutation \(\pi\) of \([K]\), let
\[
 L_\pi(a)
 :=
 \sum_{k=1}^K
 e_{\pi_k\mid S_{k-1}^\pi}(a_{S_k^\pi}).
\]
If
\(
\sup_{S,i,j,a}|\Omega_{ij\mid S}(a)|\le\varepsilon
\),
then for any permutations \(\pi,\sigma\),
\begin{equation}
 |L_\pi(a)-L_\sigma(a)|
 \le d_{\rm KT}(\pi,\sigma)\varepsilon
 \le \frac{K(K-1)}2\varepsilon,
 \label{eq:holonomybound}
\end{equation}
where \(d_{\rm KT}\) is adjacent-transposition distance.
\end{proposition}

\par\medskip

\subsection{Curvature repair and the nearest certified objective}
\label{subsec:curvaturerepair}

The preceding result audits order dependence but does not yet repair it.
For finite auxiliary spaces, the same geometry gives both the closest
scalar ledger and the closest genuine elimination tower.

For each \(i\), let \(\operatorname{Star}(\cA_i)\) be the graph with one
central ``unassigned'' vertex and leaves \(a_i\in\cA_i\).  The product
cubical complex
\[
 \mathsf X=\prod_{i=1}^K\operatorname{Star}(\cA_i)
\]
has partial assignments \((S,a_S)\) as vertices, the retained-coordinate
moves as oriented edges, and the compatible elimination squares in
\eqref{eq:squareflatness} as two-cells.  It is contractible.  Let
\[
 \mathbf d_0:C^0(\mathsf X)\to C^1(\mathsf X),
 \qquad
 \mathbf d_1:C^1(\mathsf X)\to C^2(\mathsf X)
\]
be the vertex-to-edge incidence and square-circulation operators.  Then
\begin{equation}
 \mathbf d_1\mathbf d_0=0,
 \qquad
 \ker\mathbf d_1=\operatorname{im}\mathbf d_0.
 \label{eq:eliminationcomplexexact}
\end{equation}
Thus an edge field is the increment field of a scalar ledger exactly when
its square curvature vanishes.

Let \(W\succ0\) be diagonal edge weights,
\(\|u\|_W^2=u^\top Wu\), and define
\begin{equation}
 \Omega=\mathbf d_1e,\qquad
 L_2=\mathbf d_1W^{-1}\mathbf d_1^\top.
 \label{eq:curvaturelaplacian}
\end{equation}
A touching-witness map chooses, for every \((S,a_S,j)\), one state
\(\omega_{j\mid S}(a_S)\in\cA_j\).  Define the witness-certified cone
\begin{align}
 \mathcal C_\omega
 :=\bigl\{f=\mathbf d_0V:\;&
 f_{j\mid S}(a_S,a_j)\ge0\ \text{for every edge},
 \nonumber\\[-2mm]
 &f_{j\mid S}\{a_S,\omega_{j\mid S}(a_S)\}=0
 \ \text{for every edge family}\bigr\}.
 \label{eq:certifiedcone}
\end{align}
It is nonempty, closed, and convex.

\begin{theorem}[Curvature floor and nearest certified repair]
\label{thm:curvaturerepair}
Suppose every \(\cA_i\) is finite and nonempty.  The closest flat ledger
\[
 e^\flat
 :=\argmin_{f\in\operatorname{im}\mathbf d_0}\|e-f\|_W^2
 =\mathbf d_0V^\flat
\]
is unique as an edge field.  In the Moore--Penrose constant gauge,
\begin{equation}
 V^\flat
 =
 (\mathbf d_0^\top W\mathbf d_0)^\dagger
 \mathbf d_0^\top We,
 \label{eq:flatrepairpotential}
\end{equation}
and
\begin{align}
 e-e^\flat
 &=
 W^{-1}\mathbf d_1^\top L_2^\dagger\Omega,
 \label{eq:flatrepairresidual}\\
 \|e-e^\flat\|_W^2
 &=
 \Omega^\top L_2^\dagger\Omega.
 \label{eq:curvaturefloor}
\end{align}
Consequently every scalar objective, through any exact edge field
\(f=\mathbf d_0V\), obeys the impossibility bound
\begin{equation}
 \|e-f\|_W^2
 \ge
 \Omega^\top L_2^\dagger\Omega.
 \label{eq:globalobjectivefloor}
\end{equation}

For every declared witness map \(\omega\), the closest certified field
\begin{equation}
 e_\omega^{\rm cert}
 :=
 \argmin_{f\in\mathcal C_\omega}\|e-f\|_W^2
 \label{eq:certifiedrepair}
\end{equation}
also exists and is unique.  Its potentials form a genuine partial-minimum
tower and, with \(H^{\rm cert}=V_{[K]}^{\rm cert}\), satisfy
\begin{equation}
 V_S^{\rm cert}(a_S)
 =
 \min_{a_{S^c}}H^{\rm cert}(a).
 \label{eq:repairedtower}
\end{equation}
The minimum repair cost has the exact two-tax decomposition
\begin{equation}
 \boxed{
 \|e-e_\omega^{\rm cert}\|_W^2
 =
 \underbrace{\Omega^\top L_2^\dagger\Omega}_{\textnormal{curvature tax}}
 +
 \underbrace{\|e^\flat-e_\omega^{\rm cert}\|_W^2}
 _{\textnormal{certification tax}}.
 }
 \label{eq:twotaxrepair}
\end{equation}
\end{theorem}

\par\medskip

The flat repair requires one sparse weighted graph-Laplacian solve and does
not enumerate the \(K!\) elimination orders.  With declared witnesses,
\eqref{eq:certifiedrepair} is a convex quadratic program in the vertex
potentials with linear equality and inequality constraints.  Without
witnesses, all certified fields form a finite union of the cones
\(\mathcal C_\omega\); the global repair becomes a witness-enumeration or
mixed-integer problem.  The tensor-product state count
\(\prod_i(1+|\cA_i|)\) is a genuine computational boundary.

\begin{corollary}[Global certified repair without declared witnesses]
\label{cor:globalcertifiedrepair}
Under the finite-space assumptions of
Theorem~\ref{thm:curvaturerepair}, let
\[
 \mathcal C_{\rm cert}:=\bigcup_\omega\mathcal C_\omega
\]
where the union ranges over all touching-witness maps.  Then
\(\mathcal C_{\rm cert}\) is exactly the set of certified flat edge
fields, and the global problem
\[
 e^{\rm cert}\in
 \argmin_{f\in\mathcal C_{\rm cert}}\|e-f\|_W^2
\]
has at least one solution.  It need not be unique.  Every global solution
is the conditional-defect field of a genuine partial-minimum tower and
obeys
\begin{equation}
 \|e-e^{\rm cert}\|_W^2
 =
 \Omega^\top L_2^\dagger\Omega
 +\|e^\flat-e^{\rm cert}\|_W^2.
 \label{eq:globaltwotaxrepair}
\end{equation}
\end{corollary}

\par\medskip

The product-of-stars complex is contractible, which is why square flatness
alone certifies a potential there.  On a general complex a second, globally
supported obstruction appears.

\begin{theorem}[Periods and the three-tax repair on a general complex]
\label{thm:threetax}
Let \(C^0\xrightarrow{\ \mathbf d_0\ }C^1
\xrightarrow{\ \mathbf d_1\ }C^2\) be finite-dimensional cochain data with
\(\mathbf d_1\mathbf d_0=0\), let \(W\succ0\) metrize \(C^1\), and define
the weighted harmonic space
\begin{equation}
 \mathcal H^1_W
 :=\ker\mathbf d_1\cap\ker(\mathbf d_0^\top W).
 \label{eq:harmonicspace}
\end{equation}
Then \(C^1=\operatorname{im}\mathbf d_0
\oplus W^{-1}\operatorname{im}\mathbf d_1^\top\oplus\mathcal H^1_W\),
\(W\)-orthogonally.  For an edge field \(e\), let
\(\Omega=\mathbf d_1e\), \(L_2=\mathbf d_1W^{-1}\mathbf d_1^\top\), and let
\(h\) be the \(W\)-orthogonal projection of \(e\) onto
\(\mathcal H^1_W\).  Then
\begin{equation}
 \operatorname{dist}_W^2(e,\operatorname{im}\mathbf d_0)
 =\Omega^\top L_2^\dagger\Omega+\|h\|_W^2,
 \label{eq:perioddistance}
\end{equation}
Let \(\mathcal C_{\rm cert}\subseteq\operatorname{im}\mathbf d_0\) be
any explicitly declared nonempty closed convex set of exact increments
satisfying the intended certification constraints, and define
\begin{equation}
 e^{\rm cert}
 :=P^W_{\mathcal C_{\rm cert}}e^\flat
 =\argmin_{f\in\mathcal C_{\rm cert}}\|e^\flat-f\|_W^2,
 \label{eq:generalcertifiedprojection}
\end{equation}
where \(e^\flat=P^W_{\operatorname{im}\mathbf d_0}e\).  Then
\(e^{\rm cert}\) is also the unique nearest point in
\(\mathcal C_{\rm cert}\) to \(e\), and the minimum certified repair
cost has the exact three-tax decomposition
\begin{equation}
 \boxed{
 \|e-e^{\rm cert}\|_W^2
 =
 \underbrace{\Omega^\top L_2^\dagger\Omega}_{\textnormal{curvature tax}}
 +\underbrace{\|h\|_W^2}_{\textnormal{period tax}}
 +\underbrace{\|e^\flat-e^{\rm cert}\|_W^2}
  _{\textnormal{certification tax}},
 }
 \label{eq:threetaxrepair}
\end{equation}
The curvature tax is visible on filled
elementary cycles; the period tax is invisible to every local face test.
On the product-of-stars complex, \(\mathcal H^1_W=\{0\}\) and
choosing \(\mathcal C_{\rm cert}=\mathcal C_\omega\) makes
\eqref{eq:threetaxrepair} reduce to \eqref{eq:twotaxrepair}.  On a general
elimination complex, nonnegativity and touching are not inherited from the
product-of-stars construction: any touching-witness maps and their induced
constraints must be declared separately as part of
\(\mathcal C_{\rm cert}\).
\end{theorem}

\par\medskip

\begin{example}[A period with zero curvature]
\label{ex:circleperiod}
Let the edge set form one oriented cycle of length \(m\) with no filled
two-cell, so \(C^2=\{0\}\) and \(\Omega=\mathbf d_1e=0\) vacuously for
every \(e\).  For any potential, the signed sum of \(\mathbf d_0V\) around
the cycle telescopes to zero, whereas the constant circulation
\(e\equiv c\ne0\) has cycle sum \(mc\ne0\).  Hence no scalar ledger
reproduces \(e\), although every local test passes.  Local square flatness
certifies a global objective only after every period vanishes; on a
general elimination complex, \(\mathcal H^1_W\) carries exactly this
obstruction, in target-loss units.
\end{example}

\begin{corollary}[A local falsification bound]
\label{cor:squarefalsification}
For any elementary square \(q\) and any scalar ledger \(f\),
\begin{equation}
 \max_{\alpha\in\partial q}|e_\alpha-f_\alpha|
 \ge \frac{|\Omega_q|}{4}.
 \label{eq:squarefalsification}
\end{equation}
\end{corollary}

\par\medskip

\begin{corollary}[Path-independent multi-layer defect ledger]
\label{cor:multilayerledger}
For the canonical edge field \eqref{eq:edgedefect}, every permutation
\(\pi\) satisfies
\begin{equation}
 H(a)-J
 =
 \sum_{k=1}^K
 d_{\pi_k\mid S_{k-1}^\pi}(a_{S_k^\pi}).
 \label{eq:multilayerchain}
\end{equation}
Consequently, if \(R_m:\cA\to[0,+\infty]\) are declared coherence costs,
\(\lambda_m\ge0\), and
\[
 J_{\boldsymbol\lambda}
 :=
 \inf_a\left\{H(a)+\sum_{m=1}^M\lambda_mR_m(a)\right\},
\]
then
\begin{equation}
 J_{\boldsymbol\lambda}-J
 =
 \inf_a\left\{
 \sum_{k=1}^K
 d_{\pi_k\mid S_{k-1}^\pi}(a_{S_k^\pi})
 +\sum_{m=1}^M\lambda_mR_m(a)
 \right\},
 \label{eq:multilayervshift}
\end{equation}
and the right side is independent of \(\pi\).
\end{corollary}

\par\medskip

\begin{remark}[Order-flat value, order-curved work]
\label{rem:ordercurvedwork}
Corollary~\ref{cor:multilayerledger} makes the accumulated value exactly
order-independent.  Intermediate computational size is not.  For a binary
model whose factors form a star with center \(c\) and \(n\) leaves,
eliminating the leaves first keeps every intermediate factor unary in
\(c\); eliminating the center first creates one factor jointly involving
all \(n\) leaves, with \(2^n\) entries.  The final value is identical, and
the work differs exponentially.  Zero value holonomy therefore does not
imply zero computational holonomy: induced width is an independent entry
of \(C_{\mathrm{comp}}\) in the lift package and is invisible to the
defect ledger \citep{Dechter1999}.
\end{remark}

The zero-curvature statement is a discrete exactness condition, and cycle
criteria for scalar potentials are classical.  The additional content here
is the elimination certificate: nonnegativity and edge touching turn the
integrated potential into a repeated partial-infimum tower.  Approximate
curvature is then an auditable ambiguity in target units, not just a
topological analogy.

\begin{remark}[Canonical defects versus reported divergences]
\label{rem:canonicalinteraction}
If every stage is a partial elimination of one \(H\), its canonical
conditional defects add along every path with no interaction remainder.  A
Bregman three-point term may nevertheless appear when an analyst decomposes
\(D_\Phi(b,a^\star)\) through an arbitrary \(\widetilde a\):
\[
 D_\Phi(b,a^\star)
 =
 D_\Phi(b,\widetilde a)+D_\Phi(\widetilde a,a^\star)
 +\langle b-\widetilde a,
 \nabla\Phi(\widetilde a)-\nabla\Phi(a^\star)\rangle.
\]
The remainder signals that the two reported legs are not, in general, the
canonical edge defects of one flat elimination tower.  Opposite divergence
orientations, a nonoptimal intermediate state, or a projection outside the
tower can all produce this symptom.
\end{remark}

\section{The representation gate}
\label{sec:liftgate}

A frontier cannot be compared across auxiliary representations until the
coupling is transported with the representation.  The relevant object is
therefore not a lift alone but a \emph{lift--coupling package}
\[
 L=(\cA,\cH_0,\cR,C_{\mathrm{comp}}),
\]
where \(C_{\mathrm{comp}}\) records the computational primitives needed to
evaluate the lift, its oracle, and its coupled block.  The cost component is
not used in the algebra below, but it prevents an exact identity from being
mistaken for a free algorithm.

\subsection{Equivalent packages and unrestricted-lift collapse}
\label{subsec:equivalentpackages}

\begin{definition}[Coupling-equivalent lift packages]
\label{def:couplingequivalent}
Two packages \(L_1=(\cA_1,\cH_1,\cR_1,C_1)\) and
\(L_2=(\cA_2,\cH_2,\cR_2,C_2)\) for the same target are
\emph{coupling-equivalent} if there is a \(\Theta\)-independent bijection
\(T:\cA_1\to\cA_2\) such that
\begin{equation}
 \cH_2\{\Theta,T(A)\}=\cH_1(\Theta,A),
 \qquad
 \cR_2\{T(A)\}=\cR_1(A)
 \label{eq:couplingequivalence}
\end{equation}
for every feasible \((\Theta,A)\).  When computational comparisons are
made, \(T\) and \(T^{-1}\) must additionally have declared negligible or
accounted cost.
\end{definition}

\begin{proposition}[Representation invariance and unrestricted collapse]
\label{prop:liftgate}
The following statements hold.
\begin{enumerate}
\item[(i)] Coupling-equivalent packages have identical defects under \(T\),
identical fidelity--coherence frontiers, and identical penalized excess
values.  Every externally inserted field and every defect-jet correction is
transported without changing its scalar ledger.
\item[(ii)] Optimization over all exact lifts of a fixed target is vacuous.
Indeed, for any smooth map \(g:\mathsf T\to\R^q\),
\begin{equation}
\cH_g(\Theta,A)
 =\cJ(\Theta)+\frac12\|A-g(\Theta)\|^2
 \label{eq:artificiallift}
\end{equation}
is a certified elimination with oracle \(A^\star(\Theta)=g(\Theta)\).
For any coherence cost with nonempty zero set, choosing
\(g(\Theta)\equiv A_0\in\cC\) creates a perfectly coherent dummy oracle
without changing the outer problem.  At an anchor \(\bar\Theta\), freezing
\(\bar A=g(\bar\Theta)\) gives
\begin{equation}
 \nabla^2_\Theta
 \cD^{\,g}_\Theta(\bar A)\big|_{\Theta=\bar\Theta}
 =Dg(\bar\Theta)^\top Dg(\bar\Theta).
 \label{eq:manufacturedmetric}
\end{equation}
Consequently, in local outer coordinates of dimension \(d\), every
positive-semidefinite metric \(M\in\mathbb R^{d\times d}\) can be
manufactured by a target-calling lift of the form
\eqref{eq:artificiallift}, provided the auxiliary dimension satisfies
\(q\ge\operatorname{rank}M\).
\end{enumerate}
\end{proposition}

\par\medskip

\subsection{Surjective refinements and quotient-faithful coherence}
\label{subsec:refinementgate}

Definition~\ref{def:couplingequivalent} treats invertible changes of
auxiliary coordinates.  Redundant expansions, component splitting, latent
augmentations, and gauge descriptions are typically noninjective.  Their
coherence cost must be transported by minimizing inside each redundant
fiber.

\begin{definition}[Surjective refinement and quotient coherence]
\label{def:surjectiveref}
Let \(\widetilde{\cA}\) and \(\cA\) be compact Hausdorff spaces and let
\(q:\widetilde{\cA}\to\cA\) be a continuous surjection.  A lifted package
\((\widetilde{\cA},\widetilde{\cH}_0,\widetilde{\cR})\) is a
\emph{\(q\)-refinement} of \((\cA,\cH_0)\) if
\begin{equation}
 \widetilde{\cH}_0(\Theta,\widetilde A)
 =
 \cH_0\{\Theta,q(\widetilde A)\}.
 \label{eq:qrefinementlift}
\end{equation}
For a lower-semicontinuous
\(\widetilde{\cR}:\widetilde{\cA}\to[0,+\infty]\), define its
\emph{fiberwise infimal descent} by
\begin{equation}
 (q_!\widetilde{\cR})(A)
 :=
 \min_{\widetilde A:\,q(\widetilde A)=A}
 \widetilde{\cR}(\widetilde A).
 \label{eq:infimaldescent}
\end{equation}
The minimum is attained because every nonempty \(q\)-fiber is compact and
\(\widetilde{\cR}\) is lower semicontinuous.
Relative to a declared coherence cost \(\cR\) on \(\cA\), the refinement
is \emph{quotient faithful} if
\begin{equation}
 q_!\widetilde{\cR}=\cR.
 \label{eq:quotientfaithful}
\end{equation}
\end{definition}

\begin{theorem}[The non-bijective representation gate]
\label{thm:refinementgate}
Under Definition~\ref{def:surjectiveref}, let
\(\cD_\Theta=\cH_0(\Theta,\cdot)-\cJ(\Theta)\) and let
\(\widetilde{\cD}_\Theta\) be the defect of the refined lift.  Locally
write
\[
 \widetilde{\mathfrak F}_\Theta(r)
 :=
 \inf_{\widetilde{\cR}(\widetilde A)\le r}
 \widetilde{\cD}_\Theta(\widetilde A),
 \qquad
 \widetilde{\mathfrak V}_\Theta(\lambda)
 :=
 \inf_{\widetilde A}
 \{\widetilde{\cD}_\Theta(\widetilde A)
 +\lambda\widetilde{\cR}(\widetilde A)\},
 \quad \lambda>0.
\]
These are the refined constrained and penalized excess values; the
corresponding reduced values and their geometry are formalized in
Section~\ref{sec:frontier}.  Then:
\begin{enumerate}
\item[(i)] the refined defect descends exactly,
\begin{equation}
 \widetilde{\cD}_\Theta=\cD_\Theta\circ q,
 \label{eq:qdefectdescent}
\end{equation}
\(q_!\widetilde{\cR}\) is lower semicontinuous, and, for every \(r\ge0\),
\begin{equation}
 q\{\widetilde A:\widetilde{\cR}(\widetilde A)\le r\}
 =
 \{A:(q_!\widetilde{\cR})(A)\le r\};
 \label{eq:qsublevelidentity}
\end{equation}

\item[(ii)] the refined frontier and penalized excess value satisfy
\begin{align}
 \widetilde{\mathfrak F}_\Theta(r)
 &=
 \inf\{\cD_\Theta(A):(q_!\widetilde{\cR})(A)\le r\},
 \label{eq:qfrontieridentity}\\
 \widetilde{\mathfrak V}_\Theta(\lambda)
 &=
 \inf_A\{
 \cD_\Theta(A)+\lambda(q_!\widetilde{\cR})(A)\},
 \qquad \lambda>0.
 \label{eq:qpenalizedidentity}
\end{align}
At \(\lambda=0\), both penalized excess values are defined directly as
\(\inf\cD_\Theta=0\), rather than through a \(0\cdot+\infty\)
convention.  Consequently a
quotient-faithful refinement has exactly the same
fidelity--coherence frontier and V path as the reduced package;

\item[(iii)] quotient faithfulness is necessary and sufficient for
universal frontier invariance, where ``universal'' is essential: the
claim is uniform over all continuous certified defects, not a converse
from one fixed defect.  More precisely, for a lower-semicontinuous
\(\cR:\cA\to[0,+\infty]\), the following are equivalent:
\begin{enumerate}
\item \(q_!\widetilde{\cR}=\cR\);
\item the two sides of \eqref{eq:qsublevelidentity}, with
\(q_!\widetilde{\cR}\) replaced by \(\cR\), agree for every \(r\);
\item for every continuous \(d:\cA\to[0,+\infty)\) with
\(\min_{\cA}d=0\), and every \(r\),
\begin{equation}
 \inf_{\widetilde{\cR}(\widetilde A)\le r}
 d\{q(\widetilde A)\}
 =
 \inf_{\cR(A)\le r}d(A).
 \label{eq:universalfrontierfaithfulness}
\end{equation}
\end{enumerate}
Here and in this equivalence, the infimum over an empty feasible set is
\(+\infty\).

\item[(iv)] fiberwise infimal descent composes.  If
\(\widehat{\cA}\) is compact Hausdorff,
\(p:\widehat{\cA}\to\widetilde{\cA}\) is another continuous surjection, and
\(\widehat{\cR}\) is lower semicontinuous, then
\begin{equation}
 (q\circ p)_!\widehat{\cR}
 =
 q_!(p_!\widehat{\cR}).
 \label{eq:qdescentcomposition}
\end{equation}
Thus a chain of quotient-faithful refinements remains quotient faithful.
\end{enumerate}
\end{theorem}

\par\medskip

\subsection{Duplication no-go and transported units}
\label{subsec:duplicationgate}

\begin{corollary}[Duplication no-go for raw-weight penalties]
\label{cor:duplicationnogo}
Let
\(\Delta_m=\{w\in[0,1]^m:\sum_iw_i=1\}\), and suppose that duplicating a
component produces a \(q_j\)-refinement of the lifted criterion.  Let
\(q_j:\Delta_{m+1}\to\Delta_m\) merge two copies of component \(j\) by
replacing their weights by their sum.  Consider separable penalties
\begin{equation}
 \cR_m(w)=\sum_{i=1}^m\phi(w_i),
 \label{eq:separableweightpenalty}
\end{equation}
for every \(m\ge2\), where the same dimension-independent generator
\(\phi:[0,1]\to[0,+\infty)\) is continuous and convex with
\(\phi(0)=0\).  Then
\begin{equation}
 (q_{j!}\cR_{m+1})(w)
 =
 \sum_{i\ne j}\phi(w_i)+2\phi(w_j/2).
 \label{eq:duplicationdescent}
\end{equation}
If the family is quotient faithful under every such duplication, then
\begin{equation}
 \phi(x)=cx,\qquad 0\le x\le1,
 \label{eq:duplicationlinear}
\end{equation}
for some \(c\ge0\).  Consequently \(\cR_m\equiv c\) on every simplex:
no nonconstant strictly convex separable raw-weight penalty is universally
duplication safe within this separable class.  For example, with
\(\phi(x)=x^2\), splitting one
component of total weight \(x\) into \(k\) identical copies changes its
infimal penalty from \(x^2\) to \(x^2/k\), although the synthesized output
and its insertion defect are unchanged.
\end{corollary}

\par\medskip

\begin{remark}[Scope and classical operation]
\label{rem:qrefinementnovelty}
Infimal projection and infimal postcomposition are classical operations in
variational analysis
\citep{RockafellarWets1998,CombettesEtAl2019Interpolation}.  The role here
is a necessary-and-sufficient representation gate for certified coupling
frontiers, together with its composition rule and duplication no-go.
Duplicate-model effects in weighted model averaging have been
recognized previously \citep{GarthwaiteMubwandarikwa2010}; the statement
here isolates the uniform quotient-faithfulness obstruction for the
dimension-independent separable raw-weight class.  This is still a
profiled minimum over a redundant auxiliary fiber,
not the functional-integral pushforward excluded in
Section~\ref{sec:boundaries}.  Nor does the theorem declare two
statistically or computationally different output maps equivalent merely
because one surjects onto the other: factorization
\eqref{eq:qrefinementlift} and the computational account remain required.
\end{remark}

\begin{remark}[No-free rescaling]
\label{rem:nofreerescaling}
Let \(B=cA\), \(c>0\), and
\(\cH_c(\Theta,B)=\cH_0(\Theta,B/c)\).  Transporting the coherence cost as
\(\cR_c(B)=\cR(B/c)\) makes this a coupling-equivalent package and changes
nothing.  If instead the same coordinate formula is reused and
\(\cR(cA)=c^2\cR(A)\), the apparent frontier and path obey
\begin{equation}
 \mathfrak F^{\,c}_\Theta(r)
 =\mathfrak F_\Theta(r/c^2),
 \qquad
 \mathfrak V^{\,c}_\Theta(\lambda)
 =\mathfrak V_\Theta(c^2\lambda).
 \label{eq:rescalingpathology}
\end{equation}
The apparent improvement is only a change of units and a relabeling of
\(\lambda\).  A numerical coherence budget has no representation-free
meaning unless its semantics are transported.
\end{remark}

\subsection{Admissible representation classes}
\label{subsec:admissibleL}

\begin{remark}[Admissible L is a representation prefix]
\label{rem:admissibleL}
Proposition~\ref{prop:liftgate} rules out a universal ``best lift'' over all
exact identities.  A meaningful design class
\(\mathfrak L_{\mathrm{adm}}\) must declare, at minimum, a common target,
permitted primitives, tractable block oracles, transported coherence
semantics, and computational cost; data-adaptive selection also needs a
capacity or sample-splitting protocol.  After quotienting by
Definition~\ref{def:couplingequivalent}, the object to compare is an
attainable fidelity--coherence--leverage--cost surface, not defect alone.
Symbolically,
\[
 L\triangleright\{P,G,X,V\}\times C,
\]
so L selects the admissible representation on which the subsequent
coupling mechanism acts.  It is a prefix to the calculus, not a sixth flat
coupling letter.
\end{remark}

\section{Frontier geometry and supported V paths}
\label{sec:frontier}

Once the representation gate has been passed, the defect and coherence
cost define an intrinsic scalar tradeoff for the declared package.  Fix
\(\Theta\) throughout this section and abbreviate
\(\cD=\cD_\Theta\).  The pointwise oracle minimizes \(\cD\) without regard
to coherence.  Hard consensus minimizes \(\cD\) over \(\cC\).  Between
them lies a constrained value function.

\subsection{The scalar frontier and its supported dual}
\label{subsec:scalarfrontier}

\begin{definition}[Fidelity--coherence frontier]
\label{def:frontier}
For \(r\ge0\), define
\begin{equation}
 \mathfrak F_\Theta(r)
 :=\inf\{\cD_\Theta(A):A\in\cA,\ \cR(A)\le r\}.
 \label{eq:frontier}
\end{equation}
Define also the penalized excess value
\begin{equation}
 \mathfrak V_\Theta(\lambda)
 :=\inf_{A\in\cA}
 \{\cD_\Theta(A)+\lambda\cR(A)\},
 \qquad \lambda>0.
 \label{eq:penalizedexcess}
\end{equation}
Define \(\mathfrak V_\Theta(0):=\inf_A\cD_\Theta(A)=0\) directly,
without adopting a convention for \(0\cdot+\infty\).
When
\(
 \cJ_\lambda(\Theta)
 =\inf_A\{\cH_0(\Theta,A)+\lambda\cR(A)\}
\), the exact identity
\begin{equation}
 \cJ_\lambda(\Theta)-\cJ(\Theta)
 =\mathfrak V_\Theta(\lambda)
 \label{eq:Vtargetshift}
\end{equation}
holds for \(\lambda>0\); at \(\lambda=0\), both sides are defined by the
uncoupled problem.
\end{definition}

Thus V is not merely ``regularized auxiliary estimation.''  It is the
infimal regularization of the defect coordinate itself.

\begin{theorem}[Frontier duality and endpoints]
\label{thm:frontier}
Suppose the oracle set is nonempty.  Then:
\begin{enumerate}
\item[(i)] \(\mathfrak F_\Theta\) is nonincreasing and nonnegative.  If
\[
 r_\star:=\inf\{\cR(A):A\in\cA^\star(\Theta)\},
\]
then \(\mathfrak F_\Theta(r)=0\) for every \(r>r_\star\), and also at
\(r=r_\star\) when the infimum defining \(r_\star\) is attained.
\item[(ii)] For every \(\lambda>0\),
\begin{equation}
 \mathfrak V_\Theta(\lambda)
 =\inf_{r\ge0}
 \{\mathfrak F_\Theta(r)+\lambda r\}.
 \label{eq:frontiertransform}
\end{equation}
Consequently \(\mathfrak V_\Theta\) is nondecreasing and concave on
\((0,\infty)\).  The transform extends continuously to its declared value
\(\mathfrak V_\Theta(0)=0\) if
\begin{equation}
 \inf\{\cD_\Theta(A):\cR(A)<\infty\}=0.
 \label{eq:finitepenaltyoracleapprox}
\end{equation}
\item[(iii)] Weak Lagrange duality gives
\begin{equation}
 \mathfrak F_\Theta(r)
 \ge
 \sup_{\lambda\ge0}
 \{\mathfrak V_\Theta(\lambda)-\lambda r\}.
 \label{eq:frontierweakdual}
\end{equation}
If \(\cA\) is convex, \(\cD_\Theta\) and \(\cR\) are proper lower
semicontinuous convex functions, and the constraint in
\eqref{eq:frontier} satisfies a standard strong-duality qualification, then
equality holds.  On any region where this recovery holds,
\(\mathfrak F_\Theta\) is convex as well as nonincreasing.
\item[(iv)] If \(\lambda>0\) and \(A_\lambda\) attains
\eqref{eq:penalizedexcess}, then it
lies on the frontier at its realized roughness:
\begin{equation}
 \cD_\Theta(A_\lambda)
 =\mathfrak F_\Theta\{\cR(A_\lambda)\}.
 \label{eq:supportedfrontier}
\end{equation}
\item[(v)] Suppose \(\cA\) is compact, \(\cD_\Theta\) and \(\cR\) are
lower semicontinuous, \(\cC\ne\varnothing\), and
\(\mathfrak F_\Theta(0)<\infty\).  Then
\begin{equation}
 \lim_{\lambda\to\infty}\mathfrak V_\Theta(\lambda)
 =\mathfrak F_\Theta(0)
 =\min_{A\in\cC}\cD_\Theta(A).
 \label{eq:frontierendpoint}
\end{equation}
Every cluster point of penalized minimizers is a hard-coherence minimizer.
\end{enumerate}
\end{theorem}

\par\medskip

\subsection{External design inefficiency}
\label{subsec:externalinefficiency}

The frontier provides a comparison between an externally designed G field
and the intrinsic V optimum at the same coherence budget.

\begin{definition}[External-globalization inefficiency]
\label{def:ginefficiency}
For any externally produced field \(A^{\rm G}\), define
\begin{equation}
 r_{\rm G}=\cR(A^{\rm G}),
 \qquad
 \operatorname{Ineff}_\Theta(A^{\rm G})
 =\cD_\Theta(A^{\rm G})-\mathfrak F_\Theta(r_{\rm G})\ge0.
 \label{eq:ginefficiency}
\end{equation}
\end{definition}

\begin{corollary}[Irreducible frontier cost plus design inefficiency]
\label{cor:gdecomp}
Every G field obeys the exact decomposition
\begin{equation}
 \cD_\Theta(A^{\rm G})
 =\mathfrak F_\Theta(r_{\rm G})
 +\operatorname{Ineff}_\Theta(A^{\rm G}).
 \label{eq:gfrontierdecomp}
\end{equation}
The first term is the smallest defect compatible with at least as much
coherence; the second is the price of the chosen external operator.  A
penalized V minimizer has zero inefficiency at its realized roughness.
\end{corollary}

This comparison is stricter than comparing two smoothers by bandwidth or
graph energy alone.  It asks whether the G field wastes defect relative to
the best field with the same declared coherence budget.

Section~\ref{subsec:leveragefilter} turns the comparison into an identity in
the unconstrained exact quadratic class: a positive-definite self-adjoint
contraction is itself the V solution of an induced penalty, and, when the
smoother commutes with the declared penalty,
\eqref{eq:ginefficiency} then compares its squared shrinkage with that of the
declared filter at matched roughness, weighted by oracle energy.

\begin{example}[A classical frontier: the common-corrector price]
\label{ex:homogenization}
Let \(A(y)\) be a symmetric, uniformly elliptic, \(Y\)-periodic
coefficient field, and for a macroscopic load \(\xi\) and a mean-zero
periodic corrector \(\phi\) set
\(H(\xi,\phi)=\tfrac12\int_Y(\xi+\nabla\phi)^\top A(y)(\xi+\nabla\phi)\,dy\).
Profiling the corrector is a certified elimination:
\(W_{\rm hom}(\xi)=\inf_\phi H(\xi,\phi)
=\tfrac12\xi^\top A_{\rm hom}\xi\), and the insertion defect of a trial
corrector is exactly the elliptic energy
\(\tfrac12\int_Y\nabla(\phi-\phi_\xi^\star)^\top
A\nabla(\phi-\phi_\xi^\star)\), with uniform ellipticity the exchange rate
to \(H^1\) corrector error \citep{JikovKozlovOleinik1994}.  For loads
\(\xi_t\) with weights \(w_t\) and mean \(\bar\xi\), forcing one shared
corrector across the family costs exactly
\begin{equation}
 \mathfrak F(0)
 =\frac12\sum_tw_t(\xi_t-\bar\xi)^\top
 \bigl(\langle A\rangle-A_{\rm hom}\bigr)(\xi_t-\bar\xi)\ \ge0,
 \qquad
 \langle A\rangle=\int_YA,
 \label{eq:correctorfrontier}
\end{equation}
because the pointwise bias--variance split gives
\(\inf_\phi\sum_tw_tH(\xi_t,\phi)
=W_{\rm hom}(\bar\xi)+\tfrac12\sum_tw_t
\delta\xi_t^\top\langle A\rangle\delta\xi_t\), while
\(\sum_tw_tW_{\rm hom}(\xi_t)\) expands with \(A_{\rm hom}\) in place of
\(\langle A\rangle\), and \(\langle A\rangle\succeq A_{\rm hom}\) is the
Voigt bound.  The hard-coherence price of corrector sharing is the Voigt
gap along the load dispersion.  The components are classical; the
frontier reading is the framework's.
\end{example}

\subsection{Work and susceptibility}
\label{subsec:frontierwork}

\begin{theorem}[Work and susceptibility on the supported frontier]
\label{thm:work}
Suppose on an open interval \(I\subset(0,\infty)\) that the minimizer
\(A_\lambda\) of \eqref{eq:penalizedexcess} is unique and interior,
\(\cR\) is finite-valued near the path,
\(\cD_\Theta\) and \(\cR\) are twice continuously differentiable, and
\begin{equation}
 M_\lambda
 :=\nabla^2\cD_\Theta(A_\lambda)
 +\lambda\nabla^2\cR(A_\lambda)\succ0.
 \label{eq:frontiermetric}
\end{equation}
Then
\begin{align}
 \mathfrak V_\Theta'(\lambda)
 &=\cR(A_\lambda),
 \label{eq:frontierslope}\\
 -\mathfrak V_\Theta''(\lambda)
 &=\nabla\cR(A_\lambda)^\top
 M_\lambda^{-1}\nabla\cR(A_\lambda)\ge0.
 \label{eq:frontiersusceptibility}
\end{align}
Hence
\begin{equation}
 \mathfrak V_\Theta(\lambda_1)-\mathfrak V_\Theta(\lambda_0)
 =\int_{\lambda_0}^{\lambda_1}\cR(A_\lambda)\,d\lambda.
 \label{eq:frontierwork}
\end{equation}
If the differentiable path extends over \((0,\infty)\), has the endpoint
limits in Theorem~\ref{thm:frontier}(v), satisfies the finite-roughness
approximation \eqref{eq:finitepenaltyoracleapprox} --- so that
\(\mathfrak V_\Theta(\lambda)\downarrow0\) as \(\lambda\downarrow0\) ---
and the improper integral is well
defined, then the hard-coherence defect is
\(\int_0^\infty\cR(A_\lambda)d\lambda\).
\end{theorem}

\par\medskip

\section{Invariant leverage and the complete quadratic fingerprint}
\label{sec:fingerprint}

The scalar frontier records how much objective value is paid for coherence,
but not which auxiliary modes transmit that payment back to the outer
geometry.  In the quadratic germ, that missing information is carried by an
invariant Schur channel and its modal resolution.

\subsection{Curvature allocation and the local coupling spectrum}
\label{subsec:curvatureallocation}

An exact lift allocates curvature between the retained outer problem and
the eliminated fiber.  The following matrix measures that allocation.
The Schur-complement identity in part~(i) is the profile-curvature formula
of Part~I's proposition ``Profile curvature is quotient geometry''
\citep{FrameworkPartI}; parts~(ii)--(iii) give the structural
invariance and frozen-defect consequences used here.

\begin{theorem}[Auxiliary leverage and curvature allocation]
\label{thm:liftleverage}
Work in finite-dimensional local coordinates.  Suppose \(\cH_0\) is
twice continuously differentiable near
\((\bar\Theta,\bar A)\), where
\(\bar A=A^\star(\bar\Theta)\), the implicit oracle branch remains the
unique profiled minimizer in a neighborhood of \(\bar\Theta\), and suppose
\[
 B_L:=\nabla^2_{AA}\cH_0(\bar\Theta,\bar A)\succ0,
 \qquad
 C_L:=\nabla^2_{\Theta A}\cH_0(\bar\Theta,\bar A).
\]
Define the \emph{auxiliary-leverage matrix}
\begin{equation}
 \cM_L:=C_LB_L^{-1}C_L^\top\succeq0.
 \label{eq:auxiliaryleverage}
\end{equation}
Then:
\begin{enumerate}
\item[(i)] the profiled and lifted curvatures satisfy the exact allocation
identity
\begin{equation}
 \nabla^2\cJ(\bar\Theta)
 =\nabla^2_{\Theta\Theta}\cH_0(\bar\Theta,\bar A)-\cM_L;
 \label{eq:curvatureallocation}
\end{equation}
\item[(ii)] \(\cM_L\) is invariant under every smooth
\(\Theta\)-independent auxiliary reparameterization in
Definition~\ref{def:couplingequivalent};
\item[(iii)] if the oracle is frozen at \(\bar A\), then its defect has
\begin{equation}
 \nabla^2_\Theta\cD_\Theta(\bar A)
 \big|_{\Theta=\bar\Theta}
 =\cM_L.
 \label{eq:frozendefectleverage}
\end{equation}
Consequently the canonical first-order exactification with no added
stabilizer has curvature
\begin{equation}
 \nabla^2\!\left[
 \cH_0(\Theta,\bar A)-d(\bar\Theta)
 -\langle\nabla d(\bar\Theta),\Theta-\bar\Theta\rangle
 \right]_{\Theta=\bar\Theta}
 =\nabla^2\cJ(\bar\Theta)+\cM_L
 =\nabla^2_{\Theta\Theta}\cH_0(\bar\Theta,\bar A).
 \label{eq:xretainsleverage}
\end{equation}
\end{enumerate}
\end{theorem}

\par\medskip

\begin{proposition}[The local coupling spectrum]
\label{prop:couplingspectrum}
At fixed \(\Theta\), suppose auxiliary minimization is unconstrained and the
local defect and coherence cost are exactly quadratic,
\begin{equation}
 \cD_\Theta(A)
 =\frac12(A-A^\star)^\top B(A-A^\star),
 \qquad
 \cR(A)=\frac12A^\top K A,
 \label{eq:quadraticpackage}
\end{equation}
where \(B\succ0\) and \(K\succeq0\).  Let
\[
 G=B^{-1/2}KB^{-1/2},
 \qquad z^\star=B^{1/2}A^\star.
\]
For each distinct eigenvalue \(\gamma\) of \(G\), let \(P_\gamma\) be its
orthogonal spectral projector and set
\(w_\gamma=\|P_\gamma z^\star\|^2\).  Then
\(\{(\gamma,w_\gamma)\}\) is invariant under every invertible linear
coupling-equivalence.  Moreover,
\begin{align}
 \mathfrak V_\Theta(\lambda)
 &=
 \frac12\sum_\gamma
 \frac{\lambda\gamma}{1+\lambda\gamma}w_\gamma,
 \label{eq:spectrumvalue}\\
 \cD_\Theta(A_\lambda)
 &=
 \frac12\sum_\gamma
 \left(\frac{\lambda\gamma}{1+\lambda\gamma}\right)^2w_\gamma,
 \label{eq:spectrumdefect}\\
 \cR(A_\lambda)
 &=
 \frac12\sum_\gamma
 \frac{\gamma}{(1+\lambda\gamma)^2}w_\gamma.
 \label{eq:spectrumroughness}
\end{align}
Thus an honest coordinate preconditioner cannot improve the local frontier:
a genuine local L change must alter the coupling spectrum, the admissible
semantics, or the computational cost.
\end{proposition}

\par\medskip

\subsection{The complete augmented modal fingerprint}
\label{subsec:completefingerprint}

\begin{theorem}[Complete linear fingerprint of an unconstrained quadratic
lift--coupling germ]
\label{thm:completequadraticfingerprint}
Work at the anchor of Theorem~\ref{thm:liftleverage}.  Let
\(h=\Theta-\bar\Theta\in\mathbb R^p\), let \(A\in\mathbb R^q\), and
assume that auxiliary minimization is unconstrained and that
\begin{align*}
 \cH_0(\bar\Theta+h,A)
 &=
 F(h)+h^\top C_LA
 +\frac12A^\top B_LA+a_L^\top A,\\
 \cR(A)&=\frac12A^\top K_LA,
 \qquad B_L\succ0,\quad K_L\succeq0.
\end{align*}
Here \(F\) is an arbitrary outer-only term, held fixed when packages are
compared.  For \(\lambda\ge0\), define
\[
 \cJ_\lambda(\bar\Theta+h)
 :=
 \inf_{A\in\mathbb R^q}
 \{\cH_0(\bar\Theta+h,A)+\lambda\cR(A)\},
 \qquad
 \mathfrak V_{\bar\Theta+h}(\lambda)
 :=
 \cJ_\lambda(\bar\Theta+h)-\cJ(\bar\Theta+h).
\]
Choose any invertible whitening \(W\) with \(W^\top W=B_L\), and set
\begin{equation}
 G=\widetilde K_L:=W^{-\top}K_LW^{-1},
 \qquad
 T=T_L:=C_LW^{-1},
 \qquad
 z_0=WA^\star(\bar\Theta).
 \label{eq:modalwhitening}
\end{equation}
Then the whitened oracle germ is
\begin{equation}
 z^\star(\bar\Theta+h)=z_0-T^\top h.
 \label{eq:whitenedquadraticgerm}
\end{equation}
For every distinct eigenvalue \(\gamma\) of \(G\), let \(P_\gamma\) be its
spectral projector, let \(m_\gamma=\operatorname{rank}P_\gamma\), and
define
\begin{align}
 \cM_{L,\gamma}&:=TP_\gamma T^\top,
 &
 b_\gamma&:=TP_\gamma z_0,
 &
 w_\gamma&:=z_0^\top P_\gamma z_0,
 \label{eq:augmentedmodalblocks}\\
 \mathcal I_\gamma
 &:=
 \begin{pmatrix}
  \cM_{L,\gamma}&b_\gamma\\
  b_\gamma^\top&w_\gamma
 \end{pmatrix}
 =
 \begin{pmatrix}T\\z_0^\top\end{pmatrix}
 P_\gamma
 \begin{pmatrix}T^\top&z_0\end{pmatrix}
 \succeq0.
 \label{eq:augmentedfingerprint}
\end{align}
Then:
\begin{enumerate}
\item[(i)] the modal leverages resolve the total leverage,
\begin{equation}
 \cM_L=\sum_\gamma\cM_{L,\gamma},
 \label{eq:modalleveragesum}
\end{equation}
and the marked collection
\begin{equation}
 \mathfrak I(G,T,z_0)
 :=
 \{(\gamma,m_\gamma,\mathcal I_\gamma)\}_\gamma
 \label{eq:completefingerprint}
\end{equation}
is independent of the whitening and invariant under every invertible
origin-preserving linear auxiliary reparameterization transporting the
defect, coherence cost, and oracle germ simultaneously;

\item[(ii)] the collection is complete within this origin-preserving
linear class, modulo the common outer-only term \(F(h)\).  Two whitened
germs \((G_1,T_1,z_{01})\) and \((G_2,T_2,z_{02})\), in the same outer
coordinates, have the same collection
\eqref{eq:completefingerprint} if and only if there is an orthogonal map
\(O\) such that
\begin{equation}
 G_2=OG_1O^\top,\qquad
 T_2=T_1O^\top,\qquad
 z_{02}=Oz_{01}.
 \label{eq:orthogonalquadracticequivalence}
\end{equation}
Before whitening, this is equivalent to the existence of an invertible
origin-preserving linear auxiliary reparameterization;

\item[(iii)] for every \(\lambda\ge0\), with
\[
 \sigma(\gamma,\lambda)=\frac{\lambda\gamma}{1+\lambda\gamma},
 \qquad
 x_h=\begin{pmatrix}-h\\1\end{pmatrix},
\]
the entire local penalized target-shift germ is
\begin{equation}
 \mathfrak V_{\bar\Theta+h}(\lambda)
 =
 \frac12\sum_\gamma
 \sigma(\gamma,\lambda)\,
 x_h^\top\mathcal I_\gamma x_h.
 \label{eq:completeVgerm}
\end{equation}
In particular,
\begin{align}
 \nabla_\Theta\mathfrak V_{\bar\Theta}(\lambda)
 &=-\sum_\gamma \sigma(\gamma,\lambda)b_\gamma,
 \label{eq:Vgradientfingerprint}\\
 \nabla_\Theta^2\mathfrak V_{\bar\Theta}(\lambda)
 &=\sum_\gamma \sigma(\gamma,\lambda)\cM_{L,\gamma}
 =\nabla^2\cJ_\lambda-\nabla^2\cJ;
 \label{eq:Vhessianfingerprint}
\end{align}

\item[(iv)] value, first-order response, and curvature exchange only through
the same positive-semidefinite account:
\begin{equation}
 \begin{pmatrix}
  \nabla_\Theta^2\mathfrak V_{\bar\Theta}(\lambda)
  &-\nabla_\Theta\mathfrak V_{\bar\Theta}(\lambda)\\
  -\nabla_\Theta\mathfrak V_{\bar\Theta}(\lambda)^\top
  &2\mathfrak V_{\bar\Theta}(\lambda)
 \end{pmatrix}
 =
 \sum_\gamma \sigma(\gamma,\lambda)\mathcal I_\gamma
 \succeq0.
 \label{eq:fingerprintexchangeblock}
\end{equation}
Equivalently, for every outer direction \(u\),
\begin{equation}
 \left|
 u^\top\nabla_\Theta\mathfrak V_{\bar\Theta}(\lambda)
 \right|^2
 \le
 2\mathfrak V_{\bar\Theta}(\lambda)\,
 u^\top
 \{\nabla^2\cJ_\lambda-\nabla^2\cJ\}u.
 \label{eq:globalfingerprintexchange}
\end{equation}
Mode by mode this reduces to
\(b_\gamma b_\gamma^\top\preceq
w_\gamma\cM_{L,\gamma}\).
\end{enumerate}
\end{theorem}

\par\medskip

\begin{example}[Why energy and leverage are not a complete fingerprint]
\label{ex:missingmodalorientation}
Let the outer and auxiliary spaces both be one-dimensional, fix
\(G=[\gamma]\) with \(\gamma>0\), and compare
\[
 (T,z_0)=(1,1)
 \qquad\textnormal{and}\qquad
 (T',z_0')=(1,-1).
\]
Both germs have the same
\((\gamma,w_\gamma,\cM_{L,\gamma})=(\gamma,1,1)\), but their cross blocks
are \(b_\gamma=1\) and \(b_\gamma'=-1\).  No one-dimensional orthogonal
change maps both \(T\) and \(z_0\) simultaneously, and
\[
 \mathfrak V_{\bar\Theta+h}(\lambda)
 =
 \frac{\sigma(\gamma,\lambda)}2(1-h)^2,
 \qquad
 \mathfrak V'_{\bar\Theta+h}(\lambda)
 =
 \frac{\sigma(\gamma,\lambda)}2(1+h)^2.
\]
Thus pointwise value and curvature do not determine the first-order outer
response.  The augmented block \(\mathcal I_\gamma\) is the missing
orientation datum.
\end{example}

\begin{remark}[Scope and classical linear algebra]
\label{rem:fingerprintnovelty}
Spectral decomposition, Gram-factor equivalence, and orthogonal
classification are classical linear algebra
\citep{HornJohnson2013MatrixAnalysis}.  The framework statement is the
complete fingerprint, under origin-preserving linear auxiliary
equivalence and modulo a common outer-only term, of the certified
unconstrained quadratic lift--coupling germ.  Its three blocks are modal
leverage, oracle energy, and their first-order exchange channel.  The
classification does not cover constrained auxiliary domains,
affine translations of the auxiliary origin, nonlinear
reparameterizations, or computational cost, and it is not a new spectral
theorem.
\end{remark}

\subsection{Leverage filtering and contraction-body design}
\label{subsec:leveragefilter}

\begin{proposition}[The leverage filter law]
\label{prop:leveragefilter}
Assume the lifted criterion and coherence cost are exactly quadratic,
\begin{align*}
 \cH_0(\Theta,A)
 &=
 \frac12\Theta^\top H_L\Theta+\Theta^\top C_LA
 +\frac12A^\top B_LA+\textnormal{affine terms},\\
 \cR(A)&=\frac12A^\top K_LA,
 \qquad B_L\succ0,\quad K_L\succeq0,
\end{align*}
and define
\(
\cJ_\lambda(\Theta)
=\inf_A\{\cH_0(\Theta,A)+\lambda\cR(A)\}
\).
With the whitening, projectors, and modal leverages of
Theorem~\ref{thm:completequadraticfingerprint}, for every \(\lambda\ge0\),
\begin{align}
 \nabla^2\cJ_\lambda
 &=
 H_L-\sum_\gamma
 \frac{1}{1+\lambda\gamma}\cM_{L,\gamma}
 \label{eq:leveragefilterlevel}\\
 &=
 \nabla^2\cJ+
 \sum_\gamma
 \frac{\lambda\gamma}{1+\lambda\gamma}\cM_{L,\gamma}.
 \label{eq:leveragefilter}
\end{align}
Thus the same scalar filter
\begin{equation}
 \sigma(\gamma,\lambda)
 :=\frac{\lambda\gamma}{1+\lambda\gamma}
 \label{eq:leveragefilterscalar}
\end{equation}
that prices the modal value shift in \eqref{eq:spectrumvalue} reallocates
outer profile curvature.  The curvature gain is positive semidefinite and
nondecreasing in \(\lambda\) in the Loewner order.  At the hard-coherence
endpoint,
\begin{equation}
 \nabla^2\cJ_\infty=H_L-\cM_{L,0},
 \label{eq:leveragefilterendpoint}
\end{equation}
where the zero-mode leverage is unchanged along the entire path.
If the columns of \(U_0\) form a basis of \(\ker K_L\), then
\begin{equation}
 \cM_{L,0}
 =
 C_LU_0(U_0^\top B_LU_0)^{-1}U_0^\top C_L^\top,
 \label{eq:zeromodeleverage}
\end{equation}
the leverage term in the direct Schur complement of the
consensus-constrained lift.
\end{proposition}

\par\medskip

\begin{corollary}[The contraction-body parameterization]
\label{cor:contractionbody}
For \(\lambda>0\), define the whitened filter operator
\begin{equation}
 F_\lambda
 :=\lambda\widetilde K_L(I+\lambda\widetilde K_L)^{-1},
 \qquad 0\preceq F_\lambda\prec I.
 \label{eq:filteroperator}
\end{equation}
In the exact quadratic class,
\begin{equation}
 \mathfrak V_\Theta(\lambda)
 =\frac12(z^\star)^\top F_\lambda z^\star,
 \qquad
 \nabla^2\cJ_\lambda-\nabla^2\cJ
 =T_LF_\lambda T_L^\top,
 \label{eq:filtertwocurrencies}
\end{equation}
where \(z^\star=WA^\star(\Theta)\).  Conversely, every symmetric
\(F\) with \(0\preceq F\prec I\) is realized, for fixed \(\lambda>0\),
by the positive-semidefinite whitened penalty
\begin{equation}
 \widetilde K_L
 =\lambda^{-1}F(I-F)^{-1}.
 \label{eq:filterinverse}
\end{equation}
Requiring a declared consensus subspace \(\cC\) to remain unpenalized is
equivalent to \(\cC\subseteq\ker F\); requiring it to be the exact zero set
of the induced penalty is equivalent to \(\ker F=\cC\).
\end{corollary}

\par\medskip

The contraction body also settles the matched-budget comparison of
Section~\ref{subsec:externalinefficiency}, and in doing so shows what the
G/V distinction really is inside a precise subclass.  Within the
unconstrained exact quadratic class, every positive-definite self-adjoint
contraction is a V solution for \emph{some} penalty.  Outside that subclass
--- for nonlinear, nonsymmetric, singular, or constraint-active maps --- no
such conclusion follows.  Inside it, what distinguishes an external map
from the declared V path is the penalty against which the resulting point
is scored.

\begin{proposition}[Positive-definite self-adjoint contractions have induced penalties,
and inefficiency measures declared-target mismatch]
\label{prop:inducedpenaltyineff}
Work in the defect-whitened coordinate of the unconstrained exact quadratic class of
Proposition~\ref{prop:couplingspectrum}, so that
\(\cD_\Theta(z)=\tfrac12\|z-z^\star\|^2\) and
\(\cR(z)=\tfrac12z^\top Gz\) with \(G\succeq0\).  Let a linear external
operator produce \(\widetilde z=(I-F)z^\star\) for a symmetric
\(F\) with \(0\preceq F\prec I\); equivalently, its smoother
\(I-F\) is a positive-definite self-adjoint contraction.  Then:
\begin{enumerate}
\item[(i)] \emph{(Induced penalty.)} For every \(\lambda>0\),
\(\widetilde z\) is the exact penalized minimizer at coupling level
\(\lambda\) for the induced coherence cost
\(\tfrac12z^\top G_Fz\) with
\begin{equation}
 G_F=\lambda^{-1}F(I-F)^{-1}\succeq0,
 \label{eq:inducedpenaltyfromfilter}
\end{equation}
where \(\ker G_F=\ker F\).  Thus this is an admissible coupling for a
declared coherence set \(\cC\) exactly when \(\ker F=\cC\); otherwise it is
an algebraic variational representation with a different zero set.  Its
inefficiency \emph{against that} cost is zero.  Only the product
\(\lambda G_F\) is identified; \(\lambda\) and \(G_F\) separately are
not.
\item[(ii)] \emph{(Inefficiency against the declared cost.)} Suppose in
addition that \(F\) commutes with \(G\).  Let \(\{e_j\}\) be an orthonormal
basis diagonalizing both, with \(Ge_j=\gamma_je_j\),
\(Fe_j=f_je_j\), and \(w_j=\langle z^\star,e_j\rangle^2\).  Then there is a
\(\lambda_\sharp\in[0,\infty]\) matching the realized roughness,
\begin{equation}
 \sum_j\gamma_j(1-f_j)^2w_j
 =\sum_j\frac{\gamma_j}{(1+\lambda_\sharp\gamma_j)^2}w_j,
 \label{eq:matchedbudget}
\end{equation}
and, with \(\sigma(\gamma,\lambda)=\lambda\gamma/(1+\lambda\gamma)\) and
the endpoint convention
\(\sigma(\gamma,\infty):=\mathbf 1\{\gamma>0\}\),
\begin{equation}
 \operatorname{Ineff}_\Theta(\widetilde z)
 =\frac12\sum_j
 \Bigl\{f_j^2-\sigma(\gamma_j,\lambda_\sharp)^2\Bigr\}w_j
 \ \ge0.
 \label{eq:inefficiencyspectral}
\end{equation}
Equality holds if and only if \(f_j=\sigma(\gamma_j,\lambda_\sharp)\) for
every \(j\) with \(w_j>0\).
\end{enumerate}
Thus, within this contraction subclass, the external smoother \(I-F\) lies
on the one-parameter V path generated by the declared cost precisely when
\(F=F_\lambda\) for some \(F_\lambda\) in
Corollary~\ref{cor:contractionbody}.  Elsewhere in the subclass,
the inefficiency is exactly the oracle-energy-weighted \emph{signed
difference of squared shrinkages} in
\eqref{eq:inefficiencyspectral}, after roughness has been matched.  The full
sum is nonnegative, but its individual modal summands need not be.
Agreement is required only on the support of the oracle energy, not between
the two penalties.
\end{proposition}

\par\medskip

\begin{remark}[What this settles about the G/V distinction]
\label{rem:GVdistinction}
A recurring objection to any mechanism taxonomy is that the labels rename
one computation.  Proposition~\ref{prop:inducedpenaltyineff} concedes the
objection at the level of realized fields and then locates what is left:
within the unconstrained quadratic class they produce the same maps only on
the positive-definite self-adjoint contraction body.  On that body, the
remaining distinction is which coherence cost the resulting point is scored
against.  That is not vacuous, because the
declared cost is what fixes the target
\(\cJ_\lambda\), the frontier, and the reported budget.
Part I's raw signed graph smoother is a worked ambient instance: with
\(F=I-S_{\tau,m}\) and \(\lambda=\tau\),
\eqref{eq:inducedpenaltyfromfilter} is exactly Part I's induced penalty
\(G_F=\tau^{-1}(I-S_{\tau,m})S_{\tau,m}^{-1}
=\tau^{-1}(S_{\tau,m}^{-1}-I)=K_\tau\), so before CDF
repair the smoother has this induced-penalty representation.  Because
\(K_\tau\) varies with \(\tau\), the identification is pointwise in
\(\tau\), not a continuation path of one fixed declared penalty.  If
repair is active, the
feasible G output and the constrained CDF-valued V solution are
projections in different metrics and need not coincide; the ambient
identity diagnoses the linear obstruction
but does not erase the domain distinction.
\end{remark}

\subsection{Frontier insufficiency and leverage on symmetry quotients}
\label{subsec:quotientleverage}

\begin{example}[A pointwise frontier does not determine leverage]
\label{ex:frontierinsufficient}
At a fixed anchor, take two whitened auxiliary modes with equal oracle
energies \(w_1=w_2>0\), distinct penalty eigenvalues
\(\gamma_1\ne\gamma_2\), and scalar outer leverages
\(m_1\ne m_2\).  Swapping \(\gamma_1,\gamma_2\) between the two modes
leaves the multiset \(\{(\gamma,w_\gamma)\}\) unchanged, so
\eqref{eq:spectrumvalue}--\eqref{eq:spectrumroughness} give the same
pointwise fidelity--coherence frontier at that anchor.  The curvature gains
differ by
\[
 \{\sigma(\gamma_1,\lambda)-\sigma(\gamma_2,\lambda)\}(m_1-m_2)\ne0.
\]
Thus the pointwise frontier \(\mathfrak F_{\bar\Theta}\) is not sufficient
for local coupling design.  This does not apply to the full germ
\(\Theta\mapsto\mathfrak F_\Theta\), whose outer derivatives can encode the
missing leverage.
\end{example}

\begin{remark}[Coupling energy is not outer leverage]
\label{rem:modalrelevance}
The pair \((\gamma,w_\gamma)\) describes how much a V path moves a realized
auxiliary field, whereas \(\cM_{L,\gamma}\) describes whether that
eigenspace is locally visible to the outer parameter.  A mode can have large
oracle energy and be strongly penalized while carrying negligible outer
leverage.  Smoothing such a mode may dramatically change the auxiliary plot
without producing a comparable second-order change in the profiled problem.
\end{remark}

\begin{proposition}[Quotient leverage for orbit-degenerate oracles]
\label{prop:quotientleverage}
Let a compact group \(G\) act smoothly on the auxiliary space, let
\(\cH_0\) be \(C^2\) and \(G\)-invariant in \(A\), and suppose that near
the anchor the oracle set is the orbit
\(\mathcal O_\Theta=G\cdot a^\star(\Theta)\) of a free Morse--Bott
critical manifold: the auxiliary Hessian vanishes on the orbit tangent and
is positive definite transverse to the orbit.  At an oracle representative
\(A\), first form the intrinsic quotient tangent space
\begin{equation}
 \mathsf Q_A
 :=T_A\mathcal A/T_A(G\!\cdot\!A).
 \label{eq:orbitquotienttangent}
\end{equation}
The auxiliary Hessian and mixed Hessian descend, respectively, to a
positive-definite bilinear form \(B_{\rm quot}\) on \(\mathsf Q_A\) and a
linear cross map \(C_{\rm quot}\) from \(\mathsf Q_A\) to the outer
tangent dual.  Define intrinsically
\begin{equation}
 \cM_{\rm quot}
 :=C_{\rm quot}B_{\rm quot}^{-1}C_{\rm quot}^\top.
 \label{eq:quotientleverage}
\end{equation}
For any local slice with normal tangent \(N\), the quotient map
\(N\to\mathsf Q_A\) is an isomorphism; if the pulled-back blocks are
written \(B_\perp=\nabla^2_{NN}\cH_0\) and
\(C_\perp=\nabla^2_{\Theta N}\cH_0\), then
\(\cM_{\rm quot}=C_\perp B_\perp^{-1}C_\perp^\top\).
Then:
\begin{enumerate}
\item[(i)] differentiating the invariance
\(\cH_0(\Theta,g\cdot A)=\cH_0(\Theta,A)\) in \(t\) along
\(g=\exp(t\xi)\) and then in \(\Theta\) shows that
\(\nabla^2_{\Theta A}\cH_0\) annihilates every orbit-tangent direction;
hence the quotient blocks are well defined, and the displayed slice
formula is independent of the slice and of the orbit representative;
\item[(ii)] the curvature-allocation and frozen-defect identities of
Theorem~\ref{thm:liftleverage} hold with \(\cM_L\) replaced by
\(\cM_{\rm quot}\);
\item[(iii)] a \(G\)-invariant coherence cost descends to the orbit
space, and the filter law of Proposition~\ref{prop:leveragefilter}
applies to the normal modes;
\item[(iv)] if \(\cR(g\cdot A)\ne\cR(A)\) for some \(g\), the reported
coherence price differs between representatives of the same zero-defect
oracle: it is not defined on the eliminated object at all.
\end{enumerate}
\end{proposition}

\par\medskip

\begin{example}[Discriminant barrier]
\label{ex:discriminant}
For \(H_\theta(a)=(a^2-\theta)^2\) with \(\theta>0\), the two oracle
branches \(a^\star_\pm=\pm\sqrt\theta\) have
\(B_\perp=8\theta\),
\(C_\perp=\mp4\sqrt\theta\), and
\(da^\star_\pm/d\theta=\pm1/(2\sqrt\theta)\).  Thus
the normal curvature \(B_\perp\) vanishes, while \(B_\perp^{-1}\) and
the oracle sensitivity diverge as \(\theta\downarrow0\); the quotient
leverage remains
\(C_\perp^2/B_\perp=2\).  Leverage therefore does not recover the loss of
normal curvature or branch separation.  Distance from the discriminant is
an independent exchange currency, and quotient statements hold stratum by
stratum.
\end{example}

\begin{proposition}[Symmetry selection rule for leverage]
\label{prop:selectionrule}
Let a compact group \(G\) act linearly on the outer and auxiliary tangent
spaces at a symmetric anchor, with \(B_L\), \(K_L\), and \(C_L\)
equivariant.  Decompose the complexified spaces into isotypic components,
\(E_A=\bigoplus_\alpha V_\alpha\otimes\mathbb C^{m_\alpha}\) and
\(E_\Theta=\bigoplus_\alpha V_\alpha\otimes\mathbb C^{n_\alpha}\).  If an
auxiliary symmetry type is absent from the outer space,
\(n_\alpha=0\), then
\begin{equation}
 \cM_{L,\alpha}=0.
 \label{eq:selectionrule}
\end{equation}
Such modes can carry oracle energy and consume coherence budget, but they
transmit no outer curvature through the Schur channel; scalar per-mode
formulas additionally require multiplicity-free components.
\end{proposition}

\par\medskip

\section{Semantic admissibility of variational coupling}
\label{sec:vsemantics}

Positive-semidefinite quadratic penalties generate algebraically valid
filters, but algebraic validity does not by itself certify the semantics of
the represented auxiliary state.  This section separates contraction,
proximal, and probability-preserving claims, then gives a sharp witness for
the final distinction.

\subsection{The operator passport}
\label{subsec:vpassport}

\begin{proposition}[The V passport]
\label{prop:vpassport}
On a Euclidean auxiliary space, read a candidate V mechanism as the
one-parameter family \(\{T_\lambda\}_{\lambda>0}\) it deploys, and classify
it by the strongest of the following properties, understood throughout at
the family level: the first two classes require the member property at
every \(\lambda\), the third requires one shared \(\cR\) across all
\(\lambda\), and the fourth requires the Markov property at every
\(\lambda\).  The
four classes of families are nested, and each inclusion is strict:
\[
\begin{gathered}
 \{T_\lambda\ \textnormal{nonexpansive}\}
 \supsetneq
 \{T_\lambda\ \textnormal{firmly nonexpansive}\}\\
 \supsetneq
 \{T_\lambda=\operatorname{prox}_{\lambda\cR}:\cR\
 \textnormal{proper lsc convex}\}
 \supsetneq
 \{\textnormal{probability-valid quadratic V}\}.
\end{gathered}
\]
Here an everywhere-defined firmly nonexpansive map is exactly the
resolvent \((I+\mathsf A)^{-1}\) of a maximal monotone \(\mathsf A\), and
the last class consists of proximal resolvent families that are Markov ---
positivity- and mass-preserving on probability vectors --- for every
\(\lambda>0\).  The first two strict inclusions are witnessed by explicit
constant families.  For the first, take \(T_\lambda\equiv T=-I\): every
member is nonexpansive, yet
\(\langle Td,d\rangle=-\|d\|^2<\|Td\|^2\), so it is not firmly
nonexpansive.  The remaining two separations are as follows:
\begin{enumerate}
\item[(i)] In \(\R^2\), let \(R_\theta\) be a rotation with
\(0<\theta<\pi/2\), so that \(\cos\theta\in(0,1)\).  Then
the constant family \(T_\lambda\equiv T=\alpha R_\theta\), with
\(0<\alpha\le\cos\theta\), satisfies
\(\|Td\|^2=\alpha^2\|d\|^2\le\alpha\cos\theta\,\|d\|^2
=\langle Td,d\rangle\), hence is firmly nonexpansive; its Jacobian is
nonsymmetric, so it is the proximal map of no scalar convex penalty.
\item[(ii)] With
\(u=\tfrac1{\sqrt2}(1,-1,0)^\top\),
\(v=\tfrac1{\sqrt6}(1,1,-2)^\top\), and
\(K=uu^\top+4vv^\top\), one has \(K\succeq0\) and
\(\ker K=\operatorname{span}\{\mathbf 1\}\): the algebraically ideal
consensus penalty.  Yet
\begin{equation}
 \bigl\{(I+\lambda K)^{-1}\bigr\}_{12}
 =\frac{\lambda(8\lambda-1)}{6(1+\lambda)(1+4\lambda)}
 <0,
 \qquad 0<\lambda<\tfrac18,
 \label{eq:dirichletcounterexample}
\end{equation}
so the resolvent maps a simplex vertex to a vector with a negative
coordinate.  On a finite state space with symmetric \(K\), the resolvent
family is
Markov for every \(\lambda>0\) exactly when the quadratic form of \(K\) is
a conservative Dirichlet form --- equivalently, \(K\) is a weighted graph
Laplacian with nonpositive off-diagonal entries
\citep{FukushimaOshimaTakeda2011}.  The Laplacian default passes this
gate; a positive-semidefinite penalty in general need not.
\end{enumerate}
\end{proposition}

\par\medskip

\begin{remark}[Semantic gates for claimed V mechanisms]
\label{rem:vpassportuse}
Corollary~\ref{cor:contractionbody} parameterizes the algebraic design
space; Proposition~\ref{prop:vpassport} stratifies its semantics.  A
reported smoothing operator earns a V interpretation only at the passport
level matching its claim: stability needs contraction; a scalar objective
needs cyclic monotonicity; action on probability-valued auxiliary fields
--- responsibility fields included --- needs the Dirichlet gate.  The
Gaussian rules of Theorem~\ref{thm:stochasticleverage} act on the whitened
tangent space; deploying them on a simplex-valued field inherits the same
gate.
\end{remark}

\subsection{Probability validity versus Hilbert-scale order}
\label{subsec:threepath}

The passport asks a yes/no question of one operator.  On a fixed graph the
admissible penalties form a cone, and the question becomes quantitative: how
far can a penalty push its own modal spectrum before it exits that cone?  We
answer this first for all penalties sharing the graph's eigenbasis, then
specialize to spectral powers, and finally record what survives on a general
graph.

Throughout this subsection let
\begin{equation}
 L_{\rm path}
 =
 \begin{pmatrix}
  1&-1&0\\
  -1&2&-1\\
  0&-1&1
 \end{pmatrix},
 \qquad
 \operatorname{spec}(L_{\rm path})=\{0,1,3\},
 \label{eq:threepathlaplacian}
\end{equation}
with orthonormal eigenvectors
\(\mathbf 1/\sqrt3\),
\(\tfrac1{\sqrt2}(1,0,-1)^\top\), and
\(\tfrac1{\sqrt6}(1,-2,1)^\top\), and let
\(P_0,P_1,P_3\) be the corresponding spectral projectors; explicitly
\(P_0=\mathbf 1\mathbf 1^\top/3\), while the first rows of \(P_1\) and
\(P_3\) are \((\tfrac12,0,-\tfrac12)\) and
\((\tfrac16,-\tfrac13,\tfrac16)\).

\begin{proposition}[The Dirichlet cone on three nodes and the exact order
ceiling]
\label{prop:dirichletcone}
Let \(K=k_1P_1+k_3P_3\) be any symmetric
positive-semidefinite penalty commuting with \(L_{\rm path}\) and with
\(\ker K=\operatorname{span}\{\mathbf 1\}\), so that \(k_1,k_3>0\).  Define
its \emph{modal order}
\begin{equation}
 \nu_K:=\frac{\log(k_3/k_1)}{\log 3},
 \label{eq:modalorder}
\end{equation}
so that \(K=k_1L_{\rm path}^{\,\nu_K}\); consequently \(\nu_K=\nu\)
exactly on the ray \(\{cL_{\rm path}^{\,\nu}:c>0\}\), and in particular
when \(K=L_{\rm path}^{\,\nu}\).  Then:
\begin{enumerate}
\item[(i)] \emph{(The cone.)} \(K\) is a conservative Dirichlet form ---
equivalently a weighted graph Laplacian --- if and only if
\begin{equation}
 k_3\le 3k_1,
 \qquad\text{equivalently}\qquad
 \nu_K\le1.
 \label{eq:threenodecone}
\end{equation}
Its off-diagonal entries are
\(K_{12}=K_{23}=-k_3/3\) and
\(K_{13}=-k_1/2+k_3/6\).
\item[(ii)] \emph{(Exact failure interval.)} For every \(\lambda>0\),
\begin{equation}
 \bigl\{(I+\lambda K)^{-1}\bigr\}_{13}
 =
 \frac{\lambda\left(3k_1-k_3+2\lambda k_1k_3\right)}
 {6(1+\lambda k_1)(1+\lambda k_3)},
 \label{eq:conethresholdentry}
\end{equation}
all other entries are positive, and row sums equal one.  Hence the
resolvent family is Markov exactly for
\begin{equation}
 \lambda\ \ge\
 \lambda_\star(K)
 :=
 \frac{(k_3-3k_1)_+}{2k_1k_3}
 =
 \frac{1}{2k_1}\left(1-3^{\,1-\nu_K}\right)_+ .
 \label{eq:conethreshold}
\end{equation}
\item[(iii)] \emph{(Order ceiling.)} Consequently, among all penalties
commuting with \(L_{\rm path}\), probability validity for every coupling
level is exactly modal order at most one.  Any strictly higher modal order
buys its extra order by leaving the cone, and pays for it on the explicit
interval \(0<\lambda<\lambda_\star(K)\).  The scale factor \(1/k_1\) in
\eqref{eq:conethreshold} is a change of units in the sense of
Remark~\ref{rem:nofreerescaling}; only the shape
\(1-3^{1-\nu_K}\) is intrinsic.
\end{enumerate}
\end{proposition}

\par\medskip

\begin{figure}[t]
\centering
\includegraphics[width=\textwidth]{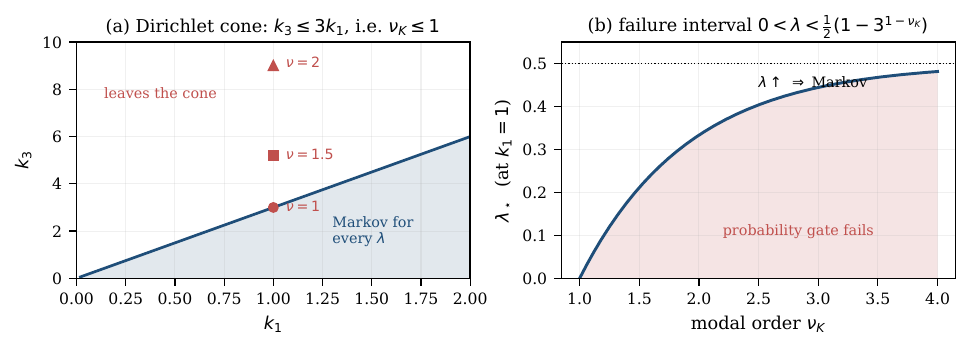}
\caption{The semantic gate on three nodes.
(a) Penalties commuting with \(L_{\rm path}\) are parameterized by their
two nonzero modal values \((k_1,k_3)\).  The shaded region
\(k_3\le3k_1\) of Proposition~\ref{prop:dirichletcone}(i) is exactly the
set of conservative Dirichlet forms, whose resolvent families are Markov at
every coupling level; markers show the spectral powers
\(K=L_{\rm path}^{\,\nu}\), which sit on the boundary at \(\nu=1\) and
leave the cone for \(\nu>1\).
(b) For a penalty outside the cone the failure is not global: the resolvent
is Markov above the threshold \eqref{eq:conethreshold} and fails below it,
so extra approximation order is paid for at weak coupling.}
\label{fig:dirichletcone}
\end{figure}

Specializing to \(K=L_{\rm path}^{\,\nu}\) makes the trade-off concrete and
adds the Hilbert-scale side of the statement, which
Proposition~\ref{prop:dirichletcone} does not see.
Figure~\ref{fig:dirichletcone} displays the cone and the failure
threshold.

\begin{proposition}[A three-node probability--approximation-order frontier]
\label{prop:threepathqualification}
With \(L_{\rm path}\) as in \eqref{eq:threepathlaplacian} and, for a real
spectral order \(\nu\ge1\), define
\begin{equation}
 \cR_\nu(a)=\frac12a^\top L_{\rm path}^{\,\nu}a,
 \qquad
 P_{\lambda,\nu}
 =(I+\lambda L_{\rm path}^{\,\nu})^{-1}.
 \label{eq:polyLaplacianresolvent}
\end{equation}
Here \(L_{\rm path}^{\,\nu}\), and below \(L_{\rm path}^{\,s}\), are
defined by spectral calculus.
Write \(\rho_\nu=3^\nu\).
For the quadratic defect
\(\cD(a)=\tfrac12\|a-a^\star\|_2^2\),
\(P_{\lambda,\nu}a^\star\) is the supported V solution.  Within this
spectral-power family:
\begin{enumerate}
\item[(i)] \(\nu=1\) is probability valid for every \(\lambda>0\):
\(P_{\lambda,1}\) is entrywise nonnegative and
\begin{equation}
 P_{\lambda,1}\mathbf 1=\mathbf 1.
 \label{eq:qoneMarkov}
\end{equation}

\item[(ii)] every \(\nu>1\) has an exact failure interval.  More precisely,
\begin{equation}
 \{L_{\rm path}^{\,\nu}\}_{13}
 =
 \frac{3^\nu-3}{6}>0,
 \label{eq:polyLaplacianpositiveoffdiag}
\end{equation}
\begin{equation}
 \{P_{\lambda,\nu}\}_{13}
 =
 \frac{\lambda(3-\rho_\nu+2\rho_\nu\lambda)}
 {6(1+\lambda)(1+\rho_\nu\lambda)},
 \label{eq:polyLaplaciannegativekernel}
\end{equation}
whereas
\[
 \{P_{\lambda,\nu}\}_{12}
 =
 \{P_{\lambda,\nu}\}_{23}
 =
 \frac{\rho_\nu\lambda}{3(1+\rho_\nu\lambda)}>0
\]
and all diagonal entries are positive.  Thus, among \(\lambda>0\),
\(P_{\lambda,\nu}\) is entrywise nonnegative if and only if
\begin{equation}
 \lambda\ge
 \frac{\rho_\nu-3}{2\rho_\nu}
 =
 \frac12(1-3^{1-\nu}).
 \label{eq:polyLaplacianthreshold}
\end{equation}
It fails the probability gate exactly below this threshold.
Throughout the failure interval it preserves total mass but not positivity.
Take the first two task laws to be \(\delta_x\) and the third to be
\(\delta_y\), with \(x\ne y\); the first output assigns mass
\(\{P_{\lambda,\nu}\}_{13}<0\) to \(\{y\}\).

\item[(iii)] in the declared \(L_{\rm path}\)-Hilbert scale, the
approximation residual has operational order \(\nu\): for every
\(0\le s\le\nu\),
\begin{equation}
 \|(I-P_{\lambda,\nu})a\|_2
 \le
 \lambda^{s/\nu}\|L_{\rm path}^{\,s}a\|_2.
 \label{eq:polyLaplacianqualification}
\end{equation}
No bound of this form, with a constant uniform as
\(\lambda\downarrow0\), extends to \(s>\nu\) on a nonzero
positive-eigenvalue mode.

\item[(iv)] separately, in the unconstrained or tangent Gaussian
observation model
\(\widehat a=a^\star+\epsilon\), with
\(\E\epsilon=0\) and
\(\E(\epsilon\epsilon^\top)=\sigma^2I\), the exact risk ledger is
\begin{align}
 \E\|P_{\lambda,\nu}\widehat a-a^\star\|_2^2
 &=
 \|(I-P_{\lambda,\nu})a^\star\|_2^2
 +\sigma^2\tr(P_{\lambda,\nu}^2)
 \label{eq:polyLaplacianrisk}\\
 &\le
 \lambda^{2s/\nu}\|L_{\rm path}^{\,s}a^\star\|_2^2
 +\sigma^2\tr(P_{\lambda,\nu}^2),
 \qquad 0\le s\le\nu.
 \nonumber
\end{align}
This risk ledger is not itself a probability-valued experiment; applying
the linear rule to probability-valued fields still requires item~(ii).
\end{enumerate}
Consequently this spectral-power resolvent family displays an exact frontier
between representation validity and approximation order.  Its first-order
member passes the Markov gate but has order one, whereas every higher-order
member has a larger uniform source-smoothness range and leaves the
probability cone exactly for
\(0<\lambda<\tfrac12(1-3^{1-\nu})\).
\end{proposition}

\par\medskip

The three-node computation is exact but small.  The part of it that
transfers to an arbitrary graph is the sign test, and it transfers with a
one-line proof.

\begin{corollary}[Positive off-diagonal entries defeat the probability gate
on any graph]
\label{cor:generalgraphgate}
Let \(K\succeq0\) be symmetric with \(K\mathbf 1=0\) on \(T\) nodes.  If
\(K_{st}>0\) for some \(s\ne t\), then the resolvent family
\(\lambda\mapsto(I+\lambda K)^{-1}\) preserves row sums for every
\(\lambda>0\) but has a strictly negative \((s,t)\) entry for all
sufficiently small \(\lambda>0\); consequently \(K\) fails the Dirichlet
gate of Proposition~\ref{prop:vpassport}.  For \(K=\mathbf L^m\) the
hypothesis holds whenever
\begin{enumerate}
\item[(a)] \(m=2\) and the connected graph is not complete, since then some
pair \(s,t\) is at distance exactly two and
\(\{\mathbf L^2\}_{st}=\sum_{k\notin\{s,t\}}w_{sk}w_{kt}>0\); or
\item[(b)] the graph is the unweighted path \(P_T\), \(T\ge3\), and
\(m\ge2\), since \(\{\mathbf L^m\}_{s,s+2}>0\) for
\(1\le s\le T-2\).  Moreover, entries sufficiently far from the boundary
obey
\(\{\mathbf L^m\}_{s,s+k}=(-1)^k\binom{2m}{m+k}\) for \(|k|\le m\) and
\(m+2\le s\le T-m-1\).
\end{enumerate}
The unweighted complete graph is a genuine exception: for \(K_T\) one has
\(\mathbf L^2=T\mathbf L\), which is again a Laplacian.
\end{corollary}

\par\medskip

\subsection{Graph-universal spectral passports and exact weighted witnesses}
\label{subsec:universalpassport}

The preceding sign test applies to any fixed penalty matrix.  For spectral
penalties one can ask a stronger, graph-independent question.  Call a
function \(\psi:[0,\infty)\to[0,\infty)\), with \(\psi(0)=0\),
\emph{graph-universally probability valid} if, for every finite connected
undirected weighted-graph Laplacian \(L\) and every \(\lambda>0\),
\[
 (I+\lambda\psi(L))^{-1}
\]
is a Markov matrix.  The next proposition gives a large positive class.
Its subordination ingredient is classical; its role here is to provide the
positive half of a validity--locality--order frontier.

\begin{proposition}[Bernstein spectral passport]
\label{prop:bernsteinpassport}
Let \(\psi\) be a Bernstein function satisfying
\(\psi(0)=0\) and \(\psi(x)>0\) for \(x>0\).  Then, for every finite
connected weighted-graph Laplacian \(L\), the matrix \(\psi(L)\) is itself
the Laplacian of a possibly dense weighted graph,
\[
 \psi(L)\succeq0,\qquad
 \psi(L)\mathbf 1=0,\qquad
 \{\psi(L)\}_{st}\le0\quad(s\ne t),
 \qquad
 \ker\psi(L)=\operatorname{span}\{\mathbf1\}.
\]
Consequently \(\psi\) is graph-universally probability valid.
More explicitly, if
\begin{equation}
 \psi(x)
 =b x+\int_{(0,\infty)}(1-e^{-ux})\,\varrho(du),
 \qquad
 b\ge0,\qquad
 \int(1\wedge u)\,\varrho(du)<\infty,
 \label{eq:bernsteinrepresentation}
\end{equation}
then the induced long-range edge weights are
\begin{equation}
 w_{st}^{\psi}
 =
 b w_{st}
 +\int_{(0,\infty)}\{e^{-uL}\}_{st}\,\varrho(du),
 \qquad s\ne t.
 \label{eq:subordinatedweights}
\end{equation}
\end{proposition}

\par\medskip

The reverse direction need not attempt to characterize every admissible
function.  A tunable three-node path already supplies an exact obstruction
whenever \(\psi(rx)>r\psi(x)\) at some point \(x>0\) and some spectral
ratio \(r\ge3\) --- the ratio range realizable by weighted three-node
paths.

\begin{proposition}[Exact weighted-path spectral obstruction]
\label{prop:weightedpathobstruction}
Fix \(x>0\) and \(r\ge3\), and set
\begin{equation}
 a,b
 =
 \frac{x}{4}
 \left[
  1+r\ \mathbin{\pm}\
  \sqrt{\frac{(3r-1)(r-3)}{3}}
 \right].
 \label{eq:weightedpathweights}
\end{equation}
For
\[
 L_{a,b}
 =
 \begin{pmatrix}
  a&-a&0\\
  -a&a+b&-b\\
  0&-b&b
 \end{pmatrix},
\]
one has \(a,b>0\) and
\(\operatorname{spec}(L_{a,b})=\{0,x,rx\}\).
Let \(\psi(0)=0\), put
\(\mu_1=\psi(x)>0\) and \(\mu_2=\psi(rx)>0\), and define
\(K_\psi=\psi(L_{a,b})\).  Then
\begin{align}
 \{K_\psi\}_{13}
 &=
 \frac{ab}{r x^2(r-1)}
 \left(\mu_2-r\mu_1\right),
 \label{eq:weightedpathpenaltyentry}\\
 \bigl\{(I+\lambda K_\psi)^{-1}\bigr\}_{13}
 &=
 \frac{ab\lambda
 \{r\mu_1-\mu_2+\lambda(r-1)\mu_1\mu_2\}}
 {r x^2(r-1)(1+\lambda\mu_1)(1+\lambda\mu_2)}.
 \label{eq:weightedpathresolvententry}
\end{align}
In particular, if \(\psi(rx)>r\psi(x)\), then \(K_\psi\) has a positive
off-diagonal entry and its resolvent has a negative \((1,3)\) entry
exactly on
\begin{equation}
 0<\lambda<
 \lambda_\star(x,r,\psi)
 :=
 \frac{\mu_2-r\mu_1}
 {(r-1)\mu_1\mu_2}.
 \label{eq:weightedpaththreshold}
\end{equation}
The unweighted three-node path is the endpoint \(x=1,r=3,a=b=1\), so
Proposition~\ref{prop:dirichletcone} is recovered by taking
\((\mu_1,\mu_2)=(k_1,k_3)\).
\end{proposition}

\par\medskip

\begin{corollary}[Universal spectral-power classification]
\label{cor:universalpowerclassification}
For a real \(\nu>0\), the following are equivalent:
\begin{enumerate}
\item[(i)] for every finite connected weighted graph \(L\) and every
\(\lambda>0\), the matrix
\((I+\lambda L^\nu)^{-1}\) is Markov;
\item[(ii)] \(0<\nu\le1\).
\end{enumerate}
When \(\nu>1\), the weighted path in
Proposition~\ref{prop:weightedpathobstruction} is an exact witness, and
its \((1,3)\) entry is negative precisely for
\begin{equation}
 0<\lambda<
 \frac{r^\nu-r}{(r-1)r^\nu x^\nu}.
 \label{eq:universalpowerthreshold}
\end{equation}
\end{corollary}

\par\medskip

\begin{corollary}[Graph-universal polynomial no-go]
\label{cor:polynomialnogo}
Let \(\psi\) be a real polynomial with \(\psi(0)=0\).  If
\(\psi(L)\) is a weighted-graph Laplacian for every finite connected
weighted-graph Laplacian \(L\), then
\[
 \psi(x)=cx
 \qquad\text{for some }c\ge0.
\]
Conversely every such function has the stated property.
\end{corollary}

\par\medskip

\begin{remark}[Scope and ownership of the spectral frontier]
\label{rem:threepathscope}
Qualification and saturation of spectral regularization filters are
classical \citep{EnglHankeNeubauer1996,ChungEspanolNguyen2014}.
Graph poly-Laplacian regularization can attain higher-smoothness regression
rates under its own sampling and geometric assumptions
\citep{GarciaTrillosMurrayThorpe2023}; higher-order graph generators are
also known to lose immediate Markov positivity
\citep{GregorioMugnolo2020}.  Bernstein subordination and classical
preservation results for Stieltjes matrices are established
\citep{SchillingSongVondracek2012,MicchelliWilloughby1979}; the positive
half of Proposition~\ref{prop:bernsteinpassport} is not claimed as a new
matrix-function characterization.  The present increment is its closure
with the exact scalable reverse witness in
Proposition~\ref{prop:weightedpathobstruction}, the all-graph
spectral-power classification, the polynomial no-go, and Part I's nonlinear
statistical escape.

At the fixed-graph level, Proposition~\ref{prop:dirichletcone} replaces a
declared power family by \emph{every} penalty sharing the three-node
graph's eigenbasis and identifies the exact boundary as a cone in the modal
spectrum, with the threshold \eqref{eq:conethreshold} and its
unit-dependence made explicit.  Proposition~\ref{prop:threepathqualification}
adds the Hilbert-scale side, which the cone alone does not see.
Corollary~\ref{cor:generalgraphgate} removes the three-node restriction
from the local sign test, while
Proposition~\ref{prop:weightedpathobstruction} restores an exact threshold
for arbitrary spectral functions by allowing the three-node edge weights
to vary.

The remaining scope is still real.  These are not theorems about faster
statistical or minimax rates, and they do not prove that every
probability-valid estimator has order at most one.  The universal
classification concerns spectral powers, and the polynomial no-go concerns
finite-degree graph-universal spectral penalties.  More generally, the
obstruction applies to linear spectral mechanisms whenever the
weighted-path growth test fires.  Nonlinear projection or a different
admissible penalty cone is a different representation and must be audited
separately.  Part I gives one such nonlinear CDF repair, exhibits its raw
signed smoother as the unconstrained ambient solution of an induced penalty
that fails the same gate under the stated sign condition, and makes clear
that activating the repair destroys the raw linear modal fingerprint
without turning the repaired G field into the constrained V solution.
\end{remark}

\begin{remark}[Barrier--repair--rate closure across the series]
\label{rem:cdf-rate-closure}
The companion \citep{FrameworkPartI} makes the statistical distinction
explicit.  At \(T\) graph vertices, it applies the signed Richardson filter
\[
 S_{\tau,m}
 =
 \sum_{j=1}^mc_j(I+a_j\tau L)^{-1}
\]
to independent empirical CDFs and then projects coordinatewise onto the CDF
set.  Before projection, the CRPS risk has an exact bias--variance
decomposition whose variance is proportional to
\(\tr(S_{\tau,m}^2)\), not to the ambient number of vertices.  Projection
can only decrease that risk and subtracts an observable squared repair
distance.  On the rescaled growing path, over the real-order graph-source
ball \(T^{-1}\|L_T^sF\|^2\le R^2\), \(s>0\), with the same sample size
\(n\) at every vertex, \(N=nT\), any fixed integer Richardson
qualification \(m\ge s\), and
\(1\le K_N=(R^2N)^{1/(4s+1)}\le c_0\min\{T,\sqrt N\}\), the repaired
estimator attains
\[
 (nT)^{-1}
 +
 R^{2/(4s+1)}(nT)^{-4s/(4s+1)},
\]
and an endpoint-Bernoulli subexperiment gives the matching minimax lower
bound.  Thus Proposition~\ref{prop:threepathqualification} is not the
statistical lower bound: it is the representation barrier showing why a
raw higher-order linear claim requires an additional validity mechanism.
The companion proves that nonlinear CDF repair supplies that mechanism
without paying in minimax order.

The phase boundary is exact for this construction.  If \(0<s\le1\),
qualification \(m=1\) is enough; the ordinary Laplacian resolvent is
Markov and the CDF projection is inactive for every input.  If \(s>1\),
using the available source smoothness requires \(m\ge2\); the raw
Richardson core is not graph-universally probability valid and has a
negative entry at sufficiently small dimensionless coupling on every fixed
path with at least three vertices, while projection restores uniform
validity without changing the
minimax order.  This does not assert that repair is active at every
finite-sample risk-optimal tuning.  Corollary~\ref{cor:universalpowerclassification} shows that
fractional powers \(0<\nu\le1\) enlarge the graph-universal Markov class
but cannot cross this order-one spectral boundary.

The two halves meet exactly, not by analogy.  Part I shows that its signed
filter satisfies \(S_{\tau,m}=(I+\tau K_\tau)^{-1}\) for the induced
penalty \(K_\tau=\tau^{-1}(S_{\tau,m}^{-1}-I)\), which is positive
semidefinite with kernel the consensus line and equals
\(A_m\tau^{m-1}\mathbf L^m+O(\tau^m)\).  Hence the Richardson family is not
a second, unrelated operator family: it enters the contraction body of
Corollary~\ref{cor:contractionbody} and is judged by the same passport.
Corollary~\ref{cor:generalgraphgate}, combined with Part I's
small-\(\tau\) control of \(K_\tau\), then supplies the verdict for every
\(m\ge2\) on the growing path.  The barrier and the construction it
obstructs are one object seen from two sides.
\end{remark}

These semantic constraints delimit what a linear V mechanism may do on a
represented state space.  We next turn to a different question: what target
preservation forces after a coupled auxiliary field has been frozen.

\section{Exactification as a unique jet normal form}
\label{sec:exactification}

Work in a finite-dimensional local chart whose interior contains the anchor.
Fix an anchor \(\bar\Theta\) and freeze any feasible, possibly nonlocal,
auxiliary field \(\widetilde A\).  As a function of \(\Theta\), write
\begin{equation}
 d(\Theta)
 :=\cD_\Theta(\widetilde A)
 =\cH_0(\Theta,\widetilde A)-\cJ(\Theta).
 \label{eq:frozendefect}
\end{equation}
Throughout the jet statements, assume both
\(\Theta\mapsto\cH_0(\Theta,\widetilde A)\) and \(\cJ\) are \(C^k\)
near the anchor; hence \(d\) is \(C^k\) there.
The raw G model is \(\cH_0(\Theta,\widetilde A)=\cJ(\Theta)+d(\Theta)\).
Let \(j^k_{\bar\Theta}f\) denote the \(k\)-jet of a \(C^k\) function at
\(\bar\Theta\), and let \(T_kf(\Theta;\bar\Theta)\) denote its Taylor
polynomial through order \(k\).

\begin{proposition}[Jet splicing and the retained auxiliary tail]
\label{prop:jetsplicing}
For a fixed feasible field \(A\), write
\(H_A(\Theta)=\cH_0(\Theta,A)\).  Whenever \(H_A\) and \(\cJ\) are
\(C^q\) near \(\bar\Theta\), with \(q>k\), the canonical order-\(k\)
exactification is
\begin{align}
 Q^{\rm can}_{X_k,A}
 &:=
 H_A-T_k(H_A-\cJ)
 \label{eq:jetsplicingdefinition}\\
 &=
 T_k\cJ+(I-T_k)H_A.
 \label{eq:jetsplicingidentity}
\end{align}
It therefore splices the target's low-order jet to the frozen lift's
higher-order tail:
\begin{equation}
 \nabla^r Q^{\rm can}_{X_k,A}(\bar\Theta)
 =
 \begin{cases}
  \nabla^r\cJ(\bar\Theta),&0\le r\le k,\\
  \nabla^r H_A(\bar\Theta),&k<r\le q.
 \end{cases}
 \label{eq:jetsplicingderivatives}
\end{equation}
For two feasible frozen fields \(A,B\),
\begin{equation}
 Q^{\rm can}_{X_k,A}-Q^{\rm can}_{X_k,B}
 =
 (I-T_k)(H_A-H_B).
 \label{eq:betweenfieldtail}
\end{equation}
Thus raising the exactification order removes, rather than adds, the
corresponding low-order channels through which the auxiliary field can
affect the model.  An \(A\)-independent \(k\)-flat stabilizer leaves
\eqref{eq:betweenfieldtail} unchanged; an \(A\)-dependent stabilizer is a
separate design choice and contributes its own between-field difference.
\end{proposition}

\par\medskip

\begin{theorem}[Converse and uniqueness of jet exactification]
\label{thm:xnormalform}
Let \(d,c,S\) be \(C^k\) near \(\bar\Theta\), and consider every corrected
scalar model of the form
\begin{equation}
 Q(\Theta)
 =\cH_0(\Theta,\widetilde A)-c(\Theta)+S(\Theta).
 \label{eq:generalcorrectedmodel}
\end{equation}
Then
\begin{equation}
 j^k_{\bar\Theta}Q=j^k_{\bar\Theta}\cJ
 \quad\Longleftrightarrow\quad
 j^k_{\bar\Theta}(c-S)=j^k_{\bar\Theta}d.
 \label{eq:xconverse}
\end{equation}
In particular, if \(S\) is \(k\)-flat at the anchor, then target-jet
preservation forces
\begin{equation}
 j^k_{\bar\Theta}c=j^k_{\bar\Theta}d.
 \label{eq:xuniqueness}
\end{equation}
The canonical representative is
\begin{equation}
 Q_{X_k}(\Theta)
 =\cH_0(\Theta,\widetilde A)
 -T_kd(\Theta;\bar\Theta)+S_{k+1}(\Theta),
 \label{eq:xcanonical}
\end{equation}
where \(S_{k+1}\) is \(k\)-flat.  Every other target-\(k\)-jet-preserving
model of the form \eqref{eq:generalcorrectedmodel} differs from
\eqref{eq:xcanonical} only by a \(k\)-flat function.
\end{theorem}

\par\medskip

The theorem is a converse, not merely a recipe.  First-order surrogate
theory already studies models whose approximation error has zero value and
gradient at the anchor \citep{Mairal2013FirstOrder}.  The new point here is
that, once a nonlocal auxiliary lift is fixed, the exact defect identifies
the required correction, while Proposition~\ref{prop:jetsplicing} states
exactly which part of the auxiliary geometry can remain.

\begin{corollary}[First-order conservation of the target score]
\label{cor:scoreconservation}
For \(k=1\) and a first-order-flat stabilizer,
\begin{equation}
 \nabla Q_{X_1}(\bar\Theta)=\nabla\cJ(\bar\Theta).
 \label{eq:scoreconservation}
\end{equation}
Therefore the frozen nonlocal field cannot survive as a permanent
first-order force while the original target score is preserved.  A
correction-free plug-in model has the exact target score only in the special
case \(\nabla d(\bar\Theta)=0\).
\end{corollary}

\begin{proposition}[Conjugate-lift transmission and an affine no-go]
\label{prop:xtransmission}
Suppose, in a local finite-dimensional chart,
\begin{equation}
 \cH_0(\Theta,A)
 =
 c(\Theta)+\Phi(A)-\langle A,\eta(\Theta)\rangle,
 \label{eq:conjugatetransmissionlift}
\end{equation}
where \(\eta\) is \(C^q\), \(q>k\).  For the canonical models of
Proposition~\ref{prop:jetsplicing},
\begin{equation}
 Q^{\rm can}_{X_k,A}(\Theta)
 -
 Q^{\rm can}_{X_k,B}(\Theta)
 =
 -\left\langle
 A-B,\,
 \eta(\Theta)-T_k\eta(\Theta;\bar\Theta)
 \right\rangle.
 \label{eq:conjugatetaildifference}
\end{equation}
In particular, the between-field curvature carried by first-order X is
\begin{equation}
 \nabla^2 Q^{\rm can}_{X_1,A}(\bar\Theta)
 -
 \nabla^2 Q^{\rm can}_{X_1,B}(\bar\Theta)
 =
 -\langle A-B,D^2\eta(\bar\Theta)\rangle.
 \label{eq:xcurvaturetransmission}
\end{equation}
Thus the linear map
\begin{equation}
 \mathsf T_X:\delta A
 \longmapsto
 -\langle\delta A,D^2\eta(\bar\Theta)\rangle
 \in\operatorname{Sym}^2(T_{\bar\Theta}^\ast)
 \label{eq:xtransmissionmap}
\end{equation}
is the exact curvature-transmission channel.  Auxiliary perturbations in
\(\ker\mathsf T_X\) are invisible to the first-order X Hessian.  If
\(\eta\) is affine, all canonical first-order X models are identical
across \(A\): exactification is automatically field-safe but transports no
auxiliary design leverage.
\end{proposition}

\par\medskip

What survives X is algorithmic rather than estimand-changing information,
but survival is conditional rather than automatic.  The Hessian of the
frozen lift, a shared preconditioner, or a proposal geometry survives only
through the tail in \eqref{eq:betweenfieldtail}; in a conjugate lift its
first available channel is \(\mathsf T_X\).  Permanent first-order sharing
requires a changed scalar target, such as the V objective, or a non-scalar
mechanism with its own exactness device.

The defect remainder also supplies an exact acceptance certificate.  For a
positive-semidefinite trust matrix \(R\), define
\begin{align}
 Q_X(\Theta)
 ={}&\cH_0(\Theta,\widetilde A)
 -d(\bar\Theta)
 -\langle\nabla d(\bar\Theta),\Theta-\bar\Theta\rangle
 \nonumber\\
 &+\frac12\|\Theta-\bar\Theta\|_R^2.
 \label{eq:xonemodel}
\end{align}
For a candidate \(\Theta^+\), let \(s=\Theta^+-\bar\Theta\) and
\begin{equation}
 \chi
 =d(\Theta^+)-d(\bar\Theta)
 -\langle\nabla d(\bar\Theta),s\rangle
 +\frac12\|s\|_R^2.
 \label{eq:xremainder}
\end{equation}

The following finite-sample budget restates, in the present frozen-field
notation, the acceptance corollary deployed in Part~I
\citep{FrameworkPartI}; it is included so that the robustness theorem
below is self-contained.

\begin{proposition}[Exactification budget]
\label{prop:xbudget}
For every candidate,
\begin{equation}
 Q_X(\Theta^+)=\cJ(\Theta^+)+\chi.
 \label{eq:xexactledger}
\end{equation}
If
\(
 Q_X(\Theta^+)\le Q_X(\bar\Theta)-\Delta^Q
\)
and, for some \(0\le\eta<1\),
\begin{equation}
 \chi\ge-\eta\Delta^Q,
 \label{eq:xbudgetcondition}
\end{equation}
then
\begin{equation}
 \cJ(\Theta^+)
 \le\cJ(\bar\Theta)-(1-\eta)\Delta^Q.
 \label{eq:xbudgetdescent}
\end{equation}
\end{proposition}

\par\medskip

\begin{theorem}[Robustness to approximate defect jets]
\label{thm:approximatejet}
Fix an integer \(k\ge1\).  At iteration \(\ell\), freeze an auxiliary field \(\widetilde A_\ell\) at
\(\Theta_\ell\), let
\[
 d_\ell(\Theta)
 =\cH_0(\Theta,\widetilde A_\ell)-\cJ(\Theta),
\]
and suppose that only an approximate correction \(\widehat c_\ell\) is
available.  Assume, on a neighborhood containing the points used below,
that \(\cJ\), \(\Theta\mapsto\cH_0(\Theta,\widetilde A_\ell)\),
\(\widehat c_\ell\), and \(S_\ell\) are \(C^k\).  For any stabilizer
\(S_\ell\), define
\begin{align}
 \widehat Q_\ell(\Theta)
 &=\cH_0(\Theta,\widetilde A_\ell)
   -\widehat c_\ell(\Theta)+S_\ell(\Theta),\\
 r_\ell(\Theta)
 &=d_\ell(\Theta)-\widehat c_\ell(\Theta)+S_\ell(\Theta).
 \label{eq:approxjetresidual}
\end{align}
Then \(\widehat Q_\ell=\cJ+r_\ell\), and consequently, for
\(0\le j\le k\),
\begin{equation}
 \left\|
 \nabla^j\widehat Q_\ell(\Theta_\ell)
 -\nabla^j\cJ(\Theta_\ell)
 \right\|
 =\left\|\nabla^j r_\ell(\Theta_\ell)\right\|.
 \label{eq:approxjeterror}
\end{equation}
In particular, a certified bound
\(\|\nabla r_\ell(\Theta_\ell)\|\le\delta_{\ell,1}\) is an equally
sharp bound on target-score error.

Let a candidate \(\Theta_{\ell+1}\) satisfy
\begin{equation}
 \widehat Q_\ell(\Theta_{\ell+1})
 \le \widehat Q_\ell(\Theta_\ell)-\Delta_\ell,
 \qquad \Delta_\ell\ge0,
 \label{eq:approxmodeldecrease}
\end{equation}
and suppose its residual transport obeys
\begin{equation}
 r_\ell(\Theta_{\ell+1})-r_\ell(\Theta_\ell)
 \ge-\eta_\ell\Delta_\ell-\varepsilon_\ell,
 \qquad 0\le\eta_\ell<1,\quad \varepsilon_\ell\ge0.
 \label{eq:approxtransportbudget}
\end{equation}
Then the original target satisfies the quantitative accepted-descent bound
\begin{equation}
 \cJ(\Theta_{\ell+1})
 \le \cJ(\Theta_\ell)
 -(1-\eta_\ell)\Delta_\ell+\varepsilon_\ell.
 \label{eq:approxaccepteddescent}
\end{equation}

For a sequence of accepted steps, assume that \(\cJ\) is bounded below,
\(\sup_\ell\eta_\ell\le\bar\eta<1\), and
\(\sum_\ell\varepsilon_\ell<\infty\).  Then
\(\sum_\ell\Delta_\ell<\infty\) and \(\cJ(\Theta_\ell)\) converges.  If,
in addition, for constants \(c,K>0\) and sequences
\(\zeta_\ell,\delta_{\ell,1}\to0\),
\begin{equation}
 \Delta_\ell\ge c\|\Theta_{\ell+1}-\Theta_\ell\|^2,
 \qquad
 \|\nabla r_\ell(\Theta_\ell)\|\le\delta_{\ell,1},
 \label{eq:approxstationarityconditions}
\end{equation}
the gradient \(\nabla\widehat Q_\ell\) is \(K\)-Lipschitz on the segment
from \(\Theta_\ell\) to \(\Theta_{\ell+1}\), and
\begin{equation}
 \|\nabla\widehat Q_\ell(\Theta_{\ell+1})\|\le\zeta_\ell,
 \label{eq:approxsubproblemsolve}
\end{equation}
then every cluster point of \(\{\Theta_\ell\}\) is stationary for
\(\cJ\).
\end{theorem}

\par\medskip

In the quadratic laboratory of Section~\ref{subsec:quadraticlab}, freezing
\(\widetilde A_\ell\) makes \(d_\ell\) an explicit quadratic function of
\(\theta\).  Approximating its linear and quadratic coefficients therefore
turns \eqref{eq:approxjeterror} into computable matrix-norm score and
curvature bounds, while \eqref{eq:approxtransportbudget} is a scalar
acceptance test.  Thus the same laboratory passes from exact X to an
implementable inexact X without changing its estimand.

\section{Coupling accounting consequences and a scoped normal form}
\label{sec:conservation}

This section collects consequences that require an explicit accounting
statement rather than another representation invariant.  The first is
quantitative: in an aligned Gaussian experiment, one admissible quadratic
filter simultaneously determines statistical risk, target deformation, and
outer-curvature transmission.

\subsection{A statistical three-currency consequence}
\label{subsec:stochasticfilter}

\begin{theorem}[Stochastic defect risk and the three-currency filter]
\label{thm:stochasticleverage}
Work in the defect-whitened coordinate of the exact quadratic class and
suppose an oracle field is observed with error
\begin{equation}
 \widehat z=z^\star+\xi,\qquad
 \E\xi=0,\qquad
 \operatorname{Cov}(\xi)=\sigma^2\Xi.
 \label{eq:noisyoracle}
\end{equation}
For a fixed center \(m\) and any symmetric contraction
\(0\preceq F\prec I\), deploy
\begin{equation}
 \widetilde z_F=m+(I-F)(\widehat z-m).
 \label{eq:filteredoracle}
\end{equation}
Conditional on \(z^\star\), its exact expected insertion defect is
\begin{equation}
 \E_\xi\frac12\|\widetilde z_F-z^\star\|^2
 =
 \frac12\|F(z^\star-m)\|^2
 +\frac{\sigma^2}{2}
 \tr\{(I-F)\Xi(I-F)\}.
 \label{eq:stochasticrisklaw}
\end{equation}

Now additionally suppose
\begin{equation}
 z^\star\sim N(m,\Omega),\qquad
 \xi\sim N(0,\sigma^2I),\qquad
 \Omega\succ0,
 \label{eq:gaussianfieldexperiment}
\end{equation}
independently, and let \(Q=\Omega^{-1}\).  If the effective precision
\(\Lambda_{\rm B}:=\sigma^2Q\) belongs to the declared admissible centered
penalty cone, then
\begin{align}
 F_{\rm B}
 &:=
 \Lambda_{\rm B}(I+\Lambda_{\rm B})^{-1}
 =\sigma^2(\Omega+\sigma^2I)^{-1},
 \label{eq:bayesfilter}\\
 g_{\rm B}(\widehat z)
 &:=
 m+(I-F_{\rm B})(\widehat z-m)
 \label{eq:bayesrule}
\end{align}
is the unique Bayes rule, up to almost-sure equality, among all measurable
finite-risk rules using the same \(\widehat z\).  For
\(
 R(g)=\frac12\E\|g(\widehat z)-z^\star\|^2
\),
\begin{equation}
 R(g)-R(g_{\rm B})
 =
 \frac12\E\|g(\widehat z)-g_{\rm B}(\widehat z)\|^2.
 \label{eq:bayespythagoras}
\end{equation}
In the joint quadratic lift of
Proposition~\ref{prop:leveragefilter}, after translating \(z-m\), let
\(\mathfrak V_{\rm B}\) and \(\cJ_{\rm B}\) denote the excess value and
profiled target obtained by adding the centered effective penalty
\[
 \frac12(z-m)^\top\Lambda_{\rm B}(z-m).
\]
Then this same filter satisfies
\begin{align}
 \mathfrak V_{\rm B}(\Theta)
 &=
 \frac12\{z^\star(\Theta)-m\}^\top
 F_{\rm B}\{z^\star(\Theta)-m\},
 \label{eq:bayestargetshift}\\
 \nabla^2\cJ_{\rm B}-\nabla^2\cJ
 &=
 T_LF_{\rm B}T_L^\top.
 \label{eq:bayescurvaturedividend}
\end{align}
Hence one admissible filter simultaneously determines statistical
bias--variance, variational target deformation, and outer curvature
transmission.

More explicitly, on a prior eigendirection with variance \(\omega>0\),
\begin{equation}
 f_{\rm B}=\frac{\sigma^2}{\omega+\sigma^2},
 \qquad
 R_{\rm B}=\frac12\frac{\sigma^2\omega}{\omega+\sigma^2}.
 \label{eq:modalbayesrisk}
\end{equation}
The uncoupled rule \(F=0\) and hard centering \(F=I\) have excess risks
\begin{equation}
 R_{\rm P}-R_{\rm B}
 =\frac12\frac{\sigma^4}{\omega+\sigma^2},
 \qquad
 R_{\rm hard}-R_{\rm B}
 =\frac12\frac{\omega^2}{\omega+\sigma^2},
 \label{eq:modalrisktriangle}
\end{equation}
while the outer curvature dividend on that mode is
\(f_{\rm B}\cM_{L,\gamma}\), where
\(\gamma=\sigma^2/\omega\) is the corresponding effective-precision
eigenvalue.
\end{theorem}

\par\medskip

Corollary~\ref{cor:contractionbody} describes the unrestricted algebraic
design space of PSD quadratic penalties.  Graph locality,
sparsity, computational constraints, or physical semantics restrict the
attainable filters to a proper subset.  Scalarized design over the
unrestricted contraction body is a semidefinite program over a
spectrahedron, not a
license to ignore the admissibility gate in
Remark~\ref{rem:admissibleL}.

Theorem~\ref{thm:stochasticleverage} is deliberately narrower than a
``Bayes completeness of V'' claim.  Gaussian posterior means, Wiener
shrinkage, and their quadratic-penalty representation are classical
\citep{KimeldorfWahba1970,Wahba1978}.  Only the product
\(\Lambda=\lambda K\) is identified; \(\lambda\) and \(K\) separately are
not.  If \(\Lambda\succeq0\) is singular, flat consensus directions
correspond to a generalized-Bayes or infinite-variance limit, not a proper
Gaussian prior.  For nonwhite \(\Xi\), the Bayes filter need not be a
symmetric Euclidean contraction, and a noise-weighted fidelity changes the
certified metric.  Finally, \eqref{eq:bayespythagoras} compares rules with
the same information set under one declared prior and noise law; it does not
say that V pointwise dominates a G rule using additional data.

\subsection{A scoped architectural normal form}
\label{subsec:scopednormalform}

The letters P/G/X/V/C are useful only if their scope is explicit.  We now
define a restricted but broad grammar covering the auxiliary-field
algorithms that motivated the theory.

\begin{definition}[One-coupling architecture]
\label{def:onecoupling}
Starting from a certified elimination \((\cH_0,\cJ,\cD)\), a
\emph{one-coupling auxiliary architecture} introduces cross-index dependence
through one of the following primitives:
\begin{enumerate}
\item[(a)] an external deterministic field map applied to a selected oracle,
\(\widetilde A=\Gamma\{A^\star(\Theta)\}\), followed by evaluation of the
old lift;
\item[(b)] an additive auxiliary coherence cost \(\lambda\cR(A)\) inserted
before elimination;
\item[(c)] after (a), a smooth scalar correction declared to preserve a
specified target \(k\)-jet;
\item[(d)] after coordinate randomized laws have been fixed, a joint law in
their fixed-marginal coupling fiber.
\end{enumerate}
The architecture may compose one deterministic primitive with (d), but no
other source of cross-index dependence is admitted.  This is an algorithmic
class definition, not a claim about arbitrary optimization procedures.
\end{definition}

\begin{theorem}[PGXV-by-C normal form in the one-coupling class]
\label{thm:normalform}
Every architecture in Definition~\ref{def:onecoupling} has the following
normal form.
\begin{enumerate}
\item Identity external mapping and exact elimination give P.
\item A nonidentity external mapping followed by raw evaluation gives G and
pays \(\cD_\Theta(\widetilde A)\).
\item If a corrected G model preserves the target \(k\)-jet, it is an
\(X_k\) model modulo a \(k\)-flat term, by
Theorem~\ref{thm:xnormalform}.
\item An additive pre-elimination coherence cost gives V, whose target shift
is \(\mathfrak V_\Theta(\lambda)\) and whose supported solutions lie on the
frontier.
\item A fixed-marginal joint-law choice gives C.  It composes orthogonally
with any deterministic case, producing P/G/X/V\(\times\)C.
\end{enumerate}
Within this class the normal-form category is unique.  Its concrete
representation is identified only up to maps with the same output field,
\(k\)-flat changes of an X correction, equivalent representations of the
same auxiliary penalty, and different couplings in the same fixed-marginal
fiber.
\end{theorem}

\par\medskip

The theorem is syntactic; its substantive consequence is the following
no-free principle.  Its four records have different native units and are
not to be added without a problem-specific exchange theorem.

\begin{theorem}[No-free first-order coupling in smooth scalar auxiliary models]
\label{thm:conservation}
Fix an anchor and a nonlocal field \(\widetilde A\), and suppose cross-index
information enters a smooth scalar update only through
\(\cH_0(\Theta,\widetilde A)\) and an optional correction.  Then:
\begin{enumerate}
\item[(i)] without correction, the model pays the exact value defect
\(d(\bar\Theta)\) and the score discrepancy \(\nabla d(\bar\Theta)\);
\item[(ii)] preserving the original target value and score forces payment of
the value and score jets of \(d\), and cancels the nonlocal field's
first-order force;
\item[(iii)] any raw first-order force attributable to the frozen field is
exactly \(\nabla d(\bar\Theta)\): preserving the original score cancels it,
whereas retaining a nonzero such force changes the target score, as V does.
If \(\nabla d(\bar\Theta)=0\), there is no raw first-order force to retain;
\item[(iv)] moving only inside a fixed-marginal stochastic fiber cannot
change any additive coordinate target and can improve only a genuinely
joint functional.
\end{enumerate}
Consequently, modulo the zero-defect symmetries and same-output
equivalences already listed in Theorem~\ref{thm:normalform}, a nontrivial
coupling in this class leaves a nonzero record in at least one of four
accounts: insertion defect, target change,
exactification or acceptance cost, or joint dependence.  The theorem does
not assert that these accounts share a unit or admit a universal exchange
rate.
\end{theorem}

\par\medskip

\begin{table}[t]
\centering
\caption{Native coupling accounts inside the one-coupling class.  Entries
in the last column are not numerically comparable without an explicit
exchange condition.}
\label{tab:conservation}
\begin{tabular}{L{0.13\textwidth}L{0.23\textwidth}L{0.25\textwidth}L{0.20\textwidth}}
\toprule
Mechanism & What is coupled & What remains exact & Native payment \\
\midrule
G & externally inserted auxiliary field & old lift, not touching & realized defect \\
X & higher-order geometry of a G field & original target jet and accepted descent & defect-jet evaluation and acceptance \\
V & auxiliary block inside one objective & block descent for the changed objective & target shift on the frontier \\
C & joint realization of fixed coordinate laws & every protected coordinate law & total correlation or another dependence cost \\
\bottomrule
\end{tabular}
\end{table}

\section{Fixed-marginal dependence and V--C separation}
\label{sec:C}

Let \(P_1,\ldots,P_M\) be probability laws on coordinate spaces
\(\mathsf X_1,\ldots,\mathsf X_M\).  Their coupling fiber is
\begin{equation}
 \Gamma(P_1,\ldots,P_M)
 =\{Q:\ Q_m=P_m,\ m=1,\ldots,M\}.
 \label{eq:couplingfiber}
\end{equation}

\begin{proposition}[Coordinate safety]
\label{prop:coordinatesafety}
For every \(Q\in\Gamma(P_1,\ldots,P_M)\) and integrable coordinate
functions \(f_m\),
\begin{equation}
 \E_Q\sum_{m=1}^M f_m(X_m)
 =\sum_{m=1}^M\E_{P_m}f_m(X_m).
 \label{eq:coordinatesafety}
\end{equation}
Thus no C operation can strictly improve an additive coordinate objective.
Any strict gain must concern variance, coverage, diversity, meeting,
capacity, communication, or another joint functional.
\end{proposition}

\par\medskip

The quadratic laboratory makes the separation visible.  Suppose
\(B_L=\operatorname{diag}(B_1,B_2)\) and randomize
\(A=A^\star+\xi\) while fixing the two marginal laws of
\(\xi_1,\xi_2\).  Every C coupling preserves
\[
 \E\cD_\theta(A)
 =\frac12\sum_{m=1}^2\E(\xi_m^\top B_m\xi_m),
\]
but it can change a joint cross-term
\(\E(\xi_1^\top K_{12}\xi_2)\).  Thus the same example that generated the
representation gate, frontier, leverage spectrum, and X residual reaches C
only by changing dependence, never by silently improving its additive
defect.

For KL-generated auxiliary defects there is a stronger exact decomposition.

\begin{theorem}[V--C separation identity]
\label{thm:VCseparation}
Let \(Q\in\Gamma(P_1,\ldots,P_M)\), and let
\(A_1^\star,\ldots,A_M^\star\) be reference probability laws.  Assume all
displayed KL divergences are well defined.  Then
\begin{equation}
 \KL\!\left(
 Q\middle\|\bigotimes_{m=1}^M A_m^\star
 \right)
 =
 \underbrace{\KL\!\left(
 Q\middle\|\bigotimes_{m=1}^M P_m
 \right)}_{\TC(Q)}
 +\sum_{m=1}^M\KL(P_m\|A_m^\star).
 \label{eq:VCseparation}
\end{equation}
The first term depends only on the joint realization inside the C fiber; the
remaining terms are marginal auxiliary defects and are invariant over that
fiber.
\end{theorem}

\par\medskip

\begin{definition}[Dependence--benefit frontier]
\label{def:dependencefrontier}
Let
\(
 Q^\otimes=\bigotimes_{m=1}^M P_m
\)
and let \(c:\prod_m\mathsf X_m\to\R\) be a bounded measurable joint cost.
For \(Q\in\Gamma(P_1,\ldots,P_M)\), define its gain over independent
realization by
\begin{equation}
 g_c(Q):=\E_{Q^\otimes}c-\E_Qc.
 \label{eq:jointgain}
\end{equation}
For \(g\ge0\), define
\begin{equation}
 \mathfrak C_{P,c}(g)
 :=\inf\left\{
 \TC(Q):Q\in\Gamma(P_1,\ldots,P_M),\ g_c(Q)\ge g
 \right\},
 \label{eq:dependencefrontier}
\end{equation}
with the infimum of the empty set equal to \(+\infty\).
\end{definition}

\begin{proposition}[The least dependence required for joint gain]
\label{prop:dependencefrontier}
The function \(\mathfrak C_{P,c}\) is nondecreasing and
\(\mathfrak C_{P,c}(0)=0\).  If
\[
 \operatorname{osc}(c):=\sup c-\inf c\in(0,\infty),
\]
then
\begin{equation}
 \mathfrak C_{P,c}(g)
 \ge \frac{2g^2}{\operatorname{osc}(c)^2}.
 \label{eq:dependencepinsker}
\end{equation}
If \(c(x_{1:M})=\sum_m c_m(x_m)\), then \(g_c(Q)=0\) throughout the
coupling fiber and \(\mathfrak C_{P,c}(g)=+\infty\) for every \(g>0\).
Finally, if \(Q_\tau\) minimizes
\begin{equation}
 \E_Qc+\tau\TC(Q)
 \quad\text{over}\quad
 Q\in\Gamma(P_1,\ldots,P_M)
 \label{eq:dependencepenalized}
\end{equation}
for some \(\tau>0\), then \(g_c(Q_\tau)\ge0\) and it lies on the
dependence--benefit frontier:
\begin{equation}
 \TC(Q_\tau)
 =\mathfrak C_{P,c}\{g_c(Q_\tau)\}.
 \label{eq:dependencesupported}
\end{equation}
\end{proposition}

\par\medskip

\begin{remark}[Entropic C is recursively certified]
\label{rem:Cselfsimilar}
Let \(Q^\otimes=\bigotimes_mP_m\), \(\tau>0\), and suppose the
fixed-marginal problem
\[
 \inf_{Q\in\Gamma(P_1,\ldots,P_M)}
 \left\{\E_Qc+\tau\KL(Q\|Q^\otimes)\right\}
\]
has an entropic optimizer \(Q_\tau^\star\) with Schr\"odinger-potential
representation
\[
 \log\frac{dQ_\tau^\star}{dQ^\otimes}(x_{1:M})
 =
 -\frac{c(x_{1:M})}{\tau}
 +\sum_{m=1}^M\varphi_m(x_m)-\kappa.
\]
Whenever the displayed quantities are finite, every
\(Q\in\Gamma(P_1,\ldots,P_M)\) satisfies
\begin{align}
 &\E_Qc+\tau\KL(Q\|Q^\otimes)
 \nonumber\\
 &\quad -
 \left\{\E_{Q_\tau^\star}c
 +\tau\KL(Q_\tau^\star\|Q^\otimes)\right\}
 =
 \tau\KL(Q\|Q_\tau^\star).
 \label{eq:Cselfsimilar}
\end{align}
Indeed, the expectations of all coordinate potentials \(\varphi_m\) are
constant on the fixed-marginal fiber.  Thus entropic C design is itself a
certified elimination on measure space: inserting a suboptimal dependence
scheme pays the exact defect \(\tau\KL(Q\|Q_\tau^\star)\).  This is a
recursive use of the calculus, not a collapse of C into V: the protected
marginals remain the defining gauge of the C fiber.
\end{remark}

\begin{remark}[Elimination--tilt interchange on rate functions]
\label{rem:ldpinterchange}
The contraction principle of large deviations is itself a certified
elimination: \(J(y)=\inf_{\pi(x)=y}I(x)\), with exact excess
\(D_y(x)=I(x)-J(y)\ge0\) for every microscopic witness.  Tilting the
microscopic rate by \(\lambda R\) before the pushforward yields
\(J^{\rm pre}_\lambda(y)=\inf_{\pi(x)=y}\{I(x)+\lambda R(x)\}\), while
tilting after contraction yields \(J^{\rm post}_\lambda=J+\lambda r\).
The two agree for every base rate function precisely when
\(R=r\circ\pi\), that is, when the penalty is constant on each eliminated
fiber; otherwise the pre-contraction tilt changes which witness wins
inside the fiber.  This is the elimination--coupling interchange law at
the level of rate functions.
\end{remark}

\begin{remark}[A horizontal frontier]
\label{rem:horizontalfrontier}
The fidelity--coherence frontier is vertical: it asks how much marginal
target fidelity must be paid for coherence.  The dependence--benefit
frontier is horizontal inside a fixed-marginal fiber: it asks how much
dependence must be paid for a genuinely joint gain.  For recursive
algorithms, conditional-kernel protection remains the exactness condition,
and the corresponding dependence budget decomposes stagewise by the
conditional chain rule for relative entropy.
\end{remark}

For recursive algorithms, fixed one-layer marginals are insufficient.  A
simple sufficient condition is co-adapted kernel preservation.
The following path-law statement is the Vertical Fiber Ledger already
proved in Part~I \citep{FrameworkPartI}; it is restated here to keep the
fixed-marginal boundary self-contained.

\begin{proposition}[Conditional-kernel protection]
\label{prop:kernelprotection}
For each coordinate \(m\), let \(\mu_m\) be an initial law and \(K_{m,s}\)
a transition kernel.  Suppose a joint process satisfies
\begin{align*}
 X_0^m&\sim\mu_m,\\
 \mathcal L(X_{s+1}^m\mid\mathcal F_s)
 &=K_{m,s}(X_s^m,\cdot)
 \quad\text{a.s.}
\end{align*}
for the full joint history \(\mathcal F_s\).  Then every coordinate path has
exactly the Markov law generated by \((\mu_m,K_{m,0:S-1})\), regardless of
the dependence between coordinates.
\end{proposition}

\par\medskip

\begin{example}[Specification exactness does not identify the law]
\label{ex:dlrphases}
Proposition~\ref{prop:kernelprotection} is constructive: the joint law is
built forward, so coordinate exactness holds by the tower property.  The
converse specification problem behaves differently.  For the
ferromagnetic Ising model on \(\mathbb Z^d\), \(d\ge2\), at zero field
and sufficiently low temperature, the plus and minus phases
\(\mu^\pm\) satisfy the same one-site conditional specification
\(\gamma\), so every conditional insertion defect
\(\E_{\mu}\KL\{\mu(d\sigma_i\mid\sigma_{i^c})\,\|\,
\gamma_i(d\sigma_i\mid\sigma_{i^c})\}\) vanishes for both, while their
magnetizations differ by \(2m_\beta>0\) \citep{Georgii2011}.  Zero local
conditional defect therefore does not identify the global object.
Propagating local certificates to global observables requires a
uniqueness currency: inside a Dobrushin uniqueness regime the
influence-matrix bound \(\|(I-C)^{-1}\|\) is the exchange rate from local
specification error to global observable error.  That bound is unavailable
once the Dobrushin uniqueness condition fails, including coexistence
regimes.  Identifiability is an independent gate, not a consequence
of exactness.
\end{example}

\section{Three hard ledger tests}
\label{sec:examples}

\subsection{Hard instance I: a forward-KL distillation barrier in soft control}
\label{sec:rlbarrier}

The orientation of a certified defect matters when one auxiliary stage is
followed by an amortized projection.  Consider a discounted MDP with one
self-looping state, \(n\) actions, zero rewards, discount
\(\gamma\in(0,1)\), and entropy temperature \(\tau>0\).  A stationary policy
is a probability vector \(p\in\Delta_n\), and its entropy-regularized return
is
\begin{equation}
 \eta_{\tau,n}(p)
 =\frac{\tau}{1-\gamma}\mathsf H(p),
 \qquad
 \mathsf H(p)=-\sum_{a=1}^np(a)\log p(a).
 \label{eq:softbanditreturn}
\end{equation}
The exact soft teacher is the uniform law \(u_n\), and elimination gives the
identity
\begin{equation}
 \eta_{\tau,n}(u_n)-\eta_{\tau,n}(p)
 =\frac{\tau}{1-\gamma}\KL(p\|u_n).
 \label{eq:softbanditdefect}
\end{equation}
Thus performance is denominated in actor-to-teacher KL.  Standard
teacher-to-actor distillation instead reports \(\KL(u_n\|p)\).  The next
result shows that these currencies cannot be exchanged by a
dimension-independent modulus.

\begin{theorem}[No dimension-free forward-KL performance certificate]
\label{thm:forwardklbarrier}
Fix \(\tau>0\) and \(\gamma\in(0,1)\).  There is no function
\(\psi:[0,\infty)\to[0,\infty)\), independent of the action dimension and
satisfying \(\psi(x)\to0\) as \(x\downarrow0\), for which
\begin{equation}
 \eta_{\tau,n}(u_n)-\eta_{\tau,n}(p)
 \le \psi\!\left\{\KL(u_n\|p)\right\}
 \label{eq:impossibleforwardcertificate}
\end{equation}
holds for every \(n\) and every full-support \(p\in\Delta_n\).

Indeed, let
\begin{equation}
 \delta_n=(\log n)^{-1/2},\qquad
 p_n(1)=\delta_n,\qquad
 p_n(a)=\frac{1-\delta_n}{n-1},\quad a=2,\ldots,n.
 \label{eq:forwardklcounterpolicy}
\end{equation}
Then, as \(n\to\infty\),
\begin{align}
 \KL(u_n\|p_n)&\longrightarrow0,
 &
 \|u_n-p_n\|_{\rm TV}&\longrightarrow0,
 \label{eq:forwardklvanishes}\\
 \eta_{\tau,n}(u_n)-\eta_{\tau,n}(p_n)
 &=\frac{\tau}{1-\gamma}\KL(p_n\|u_n)
 \longrightarrow+\infty.
 \label{eq:softgapdiverges}
\end{align}
More precisely, for all sufficiently large \(n\),
\begin{align}
 \KL(u_n\|p_n)
 &\le \frac{2}{\sqrt{\log n}},
 \label{eq:forwardklrate}\\
 \KL(p_n\|u_n)
 &\ge
 \sqrt{\log n}
 -\frac{\log\log n}{2\sqrt{\log n}}
 -\frac1e,
 \label{eq:reverseklrate}
\end{align}
and the directional exchange ratio satisfies
\begin{equation}
 \frac{\KL(p_n\|u_n)}{\KL(u_n\|p_n)}
 \sim\log n.
 \label{eq:exchangegrowth}
\end{equation}
\end{theorem}

\par\medskip

\begin{corollary}[Support restores an explicit exchange rate]
\label{cor:supportexchange}
Let \(q,p\in\Delta_n\) have full support and suppose
\(\min_aq(a)\ge\alpha>0\).  Then
\begin{equation}
 \KL(p\|q)
 \le\frac{2}{\alpha}\KL(q\|p).
 \label{eq:supportexchange}
\end{equation}
Consequently, whenever an exact soft-performance defect is
\(c\KL(p\|q)\), teacher-to-actor distillation supplies the valid certificate
\[
 c\KL(p\|q)
 \le\frac{2c}{\alpha}\KL(q\|p).
\]
For the uniform teacher, \(\alpha=1/n\); Theorem~\ref{thm:forwardklbarrier}
shows that some dimension dependence is unavoidable, and
\eqref{eq:exchangegrowth} shows that the best uniform linear exchange
constant grows at least logarithmically.
\end{corollary}

\par\medskip

The theorem is a composition obstruction, not merely a warning that KL is
asymmetric.  Exact elimination fixes the orientation in
\eqref{eq:softbanditdefect}; the actor projection reports the other
orientation.  Even an exact teacher, zero approximation error, bounded
rewards, one state, and full support do not close that ledger.  A support
floor, a density-ratio bound, a symmetric divergence, or a separately proved
exchange inequality is therefore a necessary currency in any
dimension-uniform teacher--actor certificate.  The tempting proof
\[
 \textnormal{small distillation loss}
 \Longrightarrow
 \textnormal{small soft-return loss}
\]
silently substitutes \(\KL(u_n\|p_n)\) for the certified defect
\(\KL(p_n\|u_n)\).  Theorem~\ref{thm:forwardklbarrier} shows that no omitted
dimension-free constant can repair that substitution.

\subsection{Hard instance II: quantum Gibbs fibers and observability}

Let
\[
 \mathsf S_{++}^d
 =\{\rho\in\mathbb C^{d\times d}:\rho=\rho^\ast,\ \rho\succ0,\
   \tr(\rho)=1\}
\]
be the manifold of faithful density matrices.  The Bregman generator
\(\Phi(\rho)=\tr(\rho\log\rho)\) induces the Umegaki relative entropy
\[
 D_{\rm U}(\rho\|\sigma)
 =\tr\{\rho(\log\rho-\log\sigma)\}.
\]
It is the exact defect generated by the Gibbs variational principle
\citep{Watrous2018QuantumInformation}.  Consider faithful local Gibbs states
\(\rho_t=\exp(-H_t)/\tr\{\exp(-H_t)\}\), positive weights
\(w_t\), and \(\sum_tw_t=1\).  Hard coherence restricts all auxiliary
states to a common \(\rho\).

\begin{proposition}[Noncommutative hard-coherence endpoint]
\label{prop:quantumconsensus}
Define
\begin{equation}
 L_\infty=\sum_{t=1}^T w_t\log\rho_t,\qquad
 Z_\infty=\tr\{\exp(L_\infty)\},\qquad
 \rho_\infty=\frac{\exp(L_\infty)}{Z_\infty}.
 \label{eq:quantumcenter}
\end{equation}
Then, for every \(\rho\in\mathsf S_{++}^d\),
\begin{equation}
 \sum_{t=1}^T w_tD_{\rm U}(\rho\|\rho_t)
 =D_{\rm U}(\rho\|\rho_\infty)-\log Z_\infty.
 \label{eq:quantumdecomp}
\end{equation}
Consequently the hard-coherence endpoint is the normalized
log-Euclidean mean \(\rho_\infty\), and its frontier cost is
\begin{equation}
 \mathfrak F(0)=-\log Z_\infty.
 \label{eq:quantumfrontiercost}
\end{equation}
\end{proposition}

\par\medskip

This endpoint is the quantum-Gibbs construction of Part~I
\citep{FrameworkPartI}; related variational barycenter identities for
quantum relative entropy are developed by \citet{MosonyiBunthVrana2024}.

No commutativity among the \(H_t\) is required.  On the trace-zero tangent
space, the local Hessian of \(\Phi\) is the Bogoliubov--Kubo--Mori metric
\citep{Petz1994CanonicalCorrelation},
\begin{equation}
 g_\rho(U,V)
 =\int_0^\infty
   \tr\!\left\{
    U(\rho+sI)^{-1}V(\rho+sI)^{-1}
   \right\}\,ds.
 \label{eq:BKMmetric}
\end{equation}
Thus the frontier is genuinely non-Euclidean even infinitesimally.  The
faithfulness assumption is essential here: it excludes zero-temperature and
rank-deficient boundary limits, where the logarithmic chart and uniform
quadratic control fail.

\begin{proposition}[The local measurement-to-defect exchange tax]
\label{prop:quantumlocaltax}
Let \(\rho\in\mathsf S_{++}^d\), and let
\(U=U^\ast\) with \(\tr U=0\) be a tangent direction.  Define the symmetric
logarithmic derivative \(L_U\) by
\[
 U=\frac12(\rho L_U+L_U\rho)
\]
and its metric by
\(
 g_{{\rm SLD},\rho}(U,U)=\tr(UL_U)
\).
For a POVM \(\mathsf M=\{M_x\}\), write
\begin{equation}
 I_{\mathsf M,\rho}(U)
 =
 \sum_x
 \frac{\{\tr(M_xU)\}^2}{\tr(M_x\rho)}.
 \label{eq:measurementfisher}
\end{equation}
Then
\begin{equation}
 I_{\mathsf M,\rho}(U)
 \le
 g_{{\rm SLD},\rho}(U,U)
 \le
 g_{{\rm BKM},\rho}(U,U),
 \label{eq:BKMSLDorder}
\end{equation}
and the second inequality is an equality if and only if
\([\rho,U]=0\).  Consequently the defect exchange factor
\begin{equation}
 \kappa_Q(\rho,U)
 :=
 \frac{g_{{\rm BKM},\rho}(U,U)}
      {g_{{\rm SLD},\rho}(U,U)}
 \label{eq:quantumexchangefactor}
\end{equation}
is strictly larger than one in every noncommuting direction.

For the statistical consequence, let
\(\vartheta\mapsto\rho_\vartheta\) be \(C^3\) on
\([\theta-\delta,\theta+\delta]\), suppose
\(\rho_\vartheta\succeq mI\) there for some \(m>0\), and set
\(\dot\rho_\theta=U_\theta\ne0\).  For every \(N\), let
\(\widehat\theta_N\) be obtained by an arbitrary, possibly collective, POVM
on \(\rho_\vartheta^{\otimes N}\).  Assume explicitly that
\begin{enumerate}
\item[(a)] it is locally unbiased at \(\theta\):
\(\E_\theta\widehat\theta_N=\theta\) and
\(\left.\partial_\vartheta\E_\vartheta\widehat\theta_N
\right|_{\vartheta=\theta}=1\);
\item[(b)] it is localized:
\(\widehat\theta_N\in[\theta-\delta,\theta+\delta]\) almost surely; and
\item[(c)] its third-order remainder is uniformly negligible:
\(N\E_\theta|\widehat\theta_N-\theta|^3\to0\).
\end{enumerate}
Then
\begin{equation}
 \liminf_{N\to\infty}
 N\,\E_\theta
 D_{\rm U}(\rho_{\widehat\theta_N}\|\rho_\theta)
 \ge
 \frac12\kappa_Q(\rho_\theta,U_\theta).
 \label{eq:quantumdefectCR}
\end{equation}
\end{proposition}

\par\medskip

For the qubit
\[
 \rho_p=\diag(p,1-p),\qquad
 U=\begin{pmatrix}0&u\\\overline u&0\end{pmatrix},
\]
the tax is explicit:
\begin{equation}
 \kappa_Q(p)
 =
 \frac{\log\{p/(1-p)\}}{2(2p-1)}.
 \label{eq:qubitquantumtax}
\end{equation}
It tends to one as \(p\to1/2\) and diverges as \(p\to1\).  Thus the
measurement scale and the certified defect scale can separate without
bound near the low-temperature boundary even before an auxiliary coupling
is selected.

For a single pair, equality in measured-versus-quantum relative entropy is
classical sufficiency and recovery theory
\citep{HiaiPetz1991,BertaFawziTomamichel2017}.  The family-level increment
below is the requirement that one fixed measurement preserve the entire
ledger simultaneously.

\begin{theorem}[Lossless classicalization of a defect family]
\label{thm:quantumclassicalization}
Let
\(
 \mathcal S=\{\sigma_0,\rho_1,\ldots,\rho_T\}
\)
be faithful density matrices.  The following are equivalent.
\begin{enumerate}
\item Every pair of states in \(\mathcal S\) commutes.
\item There is one fixed POVM \(\mathsf M\) such that, for every \(t\),
\begin{equation}
 D_{\rm U}(\rho_t\|\sigma_0)
 =
 \KL\!\left(
 p_{\rho_t}^{\mathsf M}
 \middle\|
 p_{\sigma_0}^{\mathsf M}
 \right),
 \qquad
 p_\rho^{\mathsf M}(x)=\tr(M_x\rho).
 \label{eq:losslessclassicalization}
\end{equation}
\end{enumerate}
Therefore a noncommuting faithful family has no fixed measurement whose
classical KL values are the entire family's exact Umegaki defect ledger.
\end{theorem}

\par\medskip

For any fixed measurement, data processing gives the exact second-level
ledger
\begin{equation}
 \boxed{
 D_{\rm U}(\rho\|\sigma)
 =
 \KL(p_\rho^{\mathsf M}\|p_\sigma^{\mathsf M})
 +\Delta_{\mathsf M}(\rho,\sigma),
 \qquad
 \Delta_{\mathsf M}(\rho,\sigma)\ge0.
 }
 \label{eq:measurementledger}
\end{equation}
Theorem~\ref{thm:quantumclassicalization} identifies commutativity as the
family-level condition under which one fixed measurement can make every
remainder vanish.  An informationally complete measurement does not
contradict the theorem: tomography may reconstruct the state and then
evaluate \(D_{\rm U}\), but inversion, finite-sample confidence, and
relative-entropy evaluation are additional operations, not equality of the
raw outcome KL with the quantum defect.

\begin{corollary}[Measured KL cannot certify an upper defect budget]
\label{cor:quantumlowercertificate}
Suppose a globalization candidate has exact quantum insertion budget
\begin{equation}
 B_h
 =
 \sum_t w_t
 D_{\rm U}(\widetilde\rho_{t,h}\|\rho_t^\star)
 \label{eq:quantumtruebudget}
\end{equation}
and an acceptance proof requires
\begin{equation}
 B_h\le(1-\varrho)\Delta_h.
 \label{eq:quantumacceptancegate}
\end{equation}
For any fixed POVM, the observable surrogate
\begin{equation}
 B_h^{\mathsf M}
 =
 \sum_t w_t
 \KL\!\left(
 p_{\widetilde\rho_{t,h}}^{\mathsf M}
 \middle\|
 p_{\rho_t^\star}^{\mathsf M}
 \right)
 \label{eq:quantummeasuredbudget}
\end{equation}
satisfies \(B_h^{\mathsf M}\le B_h\).  Therefore
\(
 B_h^{\mathsf M}\le(1-\varrho)\Delta_h
\)
does not imply \eqref{eq:quantumacceptancegate}.  A measured-KL gate is
valid only with an additional simultaneous upper bound on
\(\Delta_{\mathsf M}\), an exact quantum/tomographic certificate, or an
open change of target to a measured-relative-entropy objective.
\end{corollary}

\par\medskip

\begin{remark}[Quantum boundary of the stochastic filter]
\label{rem:quantumstochasticfilter}
In BKM-whitened multiparameter coordinates, let the measurement-induced
noise covariance be \(\Sigma_Q\succ0\) and the oracle-field prior covariance
be \(\Omega\succ0\).  The Gaussian posterior-mean shrinkage is
\[
 S_{\rm B}=\Omega(\Omega+\Sigma_Q)^{-1}.
\]
It is symmetric if and only if
\([\Omega,\Sigma_Q]=0\).  Hence a misaligned quantum experiment generally
falls outside the symmetric Euclidean V contraction body
\(S_{\rm V}=(I+\Lambda)^{-1}\).  Noise-weighted fidelity can recover the
Bayes filter, but it changes the metric fixed by the certified Umegaki
defect.  This is a metric-alignment boundary, not a failure of Gaussian
conditioning.
\end{remark}

\subsection{Hard instance III: the classical Polyakov fiber}
\label{sec:polyakov}

The third test shows the calculus reaching a classical field theory, and
shows the quotient machinery of
Proposition~\ref{prop:quotientleverage} arising by necessity rather than
by choice.  Fix a common marked compact worldsheet \(\Sigma\), let
\(X_t:\Sigma\to M\) be immersions into a Riemannian target with induced
metrics \(g_t=X_t^\ast G\succ0\), and for an auxiliary Euclidean
worldsheet metric \(h\succ0\) define the Polyakov functional
\begin{equation}
 H_t(h)
 =\frac T2\int_\Sigma\sqrt{\det h}\;
 \tr(h^{-1}g_t).
 \label{eq:polyakovlift}
\end{equation}

\begin{proposition}[Certified elimination and conformal incompatibility]
\label{prop:polyakov}
Let \(\lambda_1,\lambda_2\) denote the eigenvalues of \(h^{-1}g_t\).
Then, pointwise on \(\Sigma\),
\begin{equation}
 \tr(h^{-1}g_t)-2\sqrt{\det(h^{-1}g_t)}
 =\bigl(\sqrt{\lambda_1}-\sqrt{\lambda_2}\bigr)^2,
 \label{eq:polyakovamgm}
\end{equation}
so that
\begin{equation}
 \cJ_t=\inf_{h\succ0}H_t(h)
 =T\int_\Sigma\sqrt{\det g_t},
 \qquad
 \cD_t(h)
 =\frac T2\int_\Sigma\sqrt{\det h}\,
 \bigl(\sqrt{\lambda_1}-\sqrt{\lambda_2}\bigr)^2\ \ge0:
 \label{eq:polyakovdefect}
\end{equation}
the eliminated value is the Nambu--Goto area and the defect vanishes
exactly on the Weyl orbit \(h\propto g_t\).  This orbit is formally
analogous to the orbit-degenerate setting of
Proposition~\ref{prop:quotientleverage}; no infinite-dimensional
Morse--Bott regularity claim is needed here.  For weights \(w_t>0\),
\(\sum_tw_t=1\), linearity of \(H_t\) in \(g_t\) gives the exact
hard-coherence price of one shared worldsheet metric,
\begin{equation}
 \mathfrak F(0)
 =T\int_\Sigma
 \left[
 \sqrt{\det\Bigl(\textstyle\sum_tw_tg_t\Bigr)}
 -\sum_tw_t\sqrt{\det g_t}
 \right]\ \ge0,
 \label{eq:polyakovfrontier}
\end{equation}
with equality if and only if the metrics \(g_t\) are pointwise
proportional almost everywhere.
\end{proposition}

\par\medskip

Three caveats are mandatory.  The statement is Euclidean and classical;
a Lorentzian action is not a nonnegative minimum problem.  All induced
metrics live on the same marked worldsheet; if each surface may
reparameterize independently, one must first quotient by
\(\mathrm{Diff}\times\mathrm{Weyl}\), otherwise coordinates can
manufacture coherence --- the no-go of
Proposition~\ref{prop:quotientleverage}(iv).  And the quantum theory
integrates over \(h\), gauge-fixes, and faces anomalies and moduli: it is
a pushforward, not a profiled infimum, and lies outside the present
certificate genus (Section~\ref{sec:boundaries}).  The classical
equivalence of \eqref{eq:polyakovlift} with the Nambu--Goto action and
the Weyl degeneracy of the auxiliary metric are standard
\citep{Polyakov1981}; the multi-embedding conformal-incompatibility
frontier \eqref{eq:polyakovfrontier} is the framework combination.

\section{Relation to neighboring theories}
\label{sec:related}

\subsection{Classical components and claimed increments}
\label{subsec:classicalnew}

Table~\ref{tab:classicalnew} states the novelty boundary directly.  The
right column is the contribution claimed by this manuscript; it is not a
claim that the ingredients in the middle column are new, nor an assertion
that every possible equivalent formulation is absent from the literature.
Part I of the series \citep{FrameworkPartI} owns the constructive
P/G/X/V/C algorithms, their descent certificates, statistical examples, and
numerical diagnostics.  This Part II owns the integrability, repair,
representation, frontier, leverage, semantic-admissibility, converse, and
impossibility results summarized below.

\small
\begin{longtable}{L{0.16\textwidth}L{0.34\textwidth}L{0.40\textwidth}}
\caption{Classical ingredients versus the manuscript's claimed increment.}
\label{tab:classicalnew}\\
\toprule
Layer & Classical component & Claimed increment here \\
\midrule
\endfirsthead
\multicolumn{3}{l}{\textit{Table~\thetable\ continued}}\\[0.3em]
\toprule
Layer & Classical component & Claimed increment here \\
\midrule
\endhead
\midrule
\multicolumn{3}{r}{\textit{continued on next page}}\\
\endfoot
\bottomrule
\endlastfoot
Defect geometry
& Bregman divergences, conjugate duality, directed centroids
& Treating the lift-forced orientation as an exact target-loss coordinate,
while keeping dependence and computation in separate native accounts \\
\addlinespace
Elimination order
& Exact edge potentials, cycle tests, Hodge projection
& Nonnegativity plus touching plus square flatness as an iff criterion for a
partial-minimum tower; nearest certified repair and its two-tax identity \\
\addlinespace
Representation refinement
& Coordinate changes, infimal projection, infimal postcomposition
& An iff quotient-faithfulness gate for surjective refinements, its
composition law, and a duplication no-go for raw-weight penalties \\
\addlinespace
Frontier geometry
& Constrained value functions, Lagrange duality, parametric sensitivity
& The exact V target-shift identity, matched-budget G inefficiency read
(for penalty-commuting smoothers) as
the signed oracle-weighted difference of squared shrinkages between an
induced and declared filter, and a work--susceptibility ledger for supported
paths \\
\addlinespace
Quadratic coupling geometry
& Schur complements, Tikhonov filters, Gram classification,
Wiener--Gaussian shrinkage
& A complete augmented modal fingerprint under origin-preserving linear
equivalence; its PSD block jointly prices target-shift value, gradient, and
outer-curvature transmission, with an aligned three-currency filter
consequence \\
\addlinespace
Probability-valued graph smoothing
& Bernstein subordination, Stieltjes-matrix preservation, CRPS propriety,
graph poly-Laplacian rates, qualification and saturation
& A graph-universal positive passport closed by an exact weighted-path
reverse witness, spectral-power iff classification, and polynomial no-go;
combined with Part I's real-smoothness theorem, a nonlinear CDF repair that
restores validity without losing the matching minimax rate \\
\addlinespace
Target-preserving surrogates
& First-order surrogate touching, additive multifidelity correction, and
trust-region model management
& The converse that a frozen coupled lift forces defect-jet subtraction,
unique modulo a flat term; the target-jet/lift-tail splicing identity,
conjugate-lift transmission map and affine no-go; a representation-level
perturbation theorem that subsumes the deployable first-order rule of
Part I \\
\addlinespace
Fixed marginals
& Probability couplings, total correlation, entropic transport and
Schr\"odinger projections
& Architectural V--C separation in the same defect ledger and a recursive
certificate for entropic C design \\
\addlinespace
Soft control
& KL asymmetry, support-dependent divergence comparisons, forward/reverse
policy-improvement distinctions
& An exact-teacher sequence with vanishing forward KL and total variation
but diverging certified soft-return defect; hence no dimension-free
composition modulus \\
\addlinespace
Quantum information
& BKM and SLD metrics, data processing, Petz recovery, no broadcasting
& Reading their combination as an auxiliary-certificate theorem:
commutativity iff one fixed measurement closes the whole Umegaki ledger,
and measured KL has the wrong bound direction for defect-budget acceptance \\
\addlinespace
Elimination topology
& Hodge decomposition, harmonic representatives
& Pricing both curvature and global periods in target-loss units; the
three-tax certified repair \\
\addlinespace
Operator semantics
& Minty resolvents, Moreau proximity operators, Beurling--Deny criterion
& One falsifiable V passport; the exact fixed-graph Dirichlet cone, scalable
spectral-function obstruction, and all-weighted-graph validity--locality--%
order boundary \\
\addlinespace
Degenerate oracles
& Morse--Bott critical manifolds, Schur's lemma
& Quotient leverage on gauge orbits and an isotypic selection rule for
curvature transmission \\
\addlinespace
Classical fields
& Polyakov--Nambu--Goto equivalence, corrector variational principles,
DLR specifications
& Closed-form conformal-incompatibility and common-corrector frontiers;
zero local defect without global identifiability \\
\end{longtable}
\normalsize

Exact variable elimination combines local factors associatively; its order
controls induced width and cost \citep{Dechter1999,PeyrardEtAl2015}.
Separately, a four-cycle condition is a classical integrability test for the
unilateral increments of an exact potential game
\citep{MondererShapley1996}.  Theorem~\ref{thm:flatness} claims neither
ingredient as new.  It adds nonnegativity and edge touching to square
exactness.  Together, these conditions are necessary and sufficient for
local gaps to form a certified partial-infimum tower.
Proposition~\ref{prop:holonomy} measures their failure directly in
objective-loss units.

Least-squares projection of an inconsistent edge flow onto a global
potential is classical in HodgeRank \citep{JiangLimYaoYe2011}; related
decompositions project finite games onto potential and harmonic components
\citep{CandoganEtAl2011}.  The flat projection in
Theorem~\ref{thm:curvaturerepair} uses this algebra and is not claimed as a
new Hodge decomposition.  The framework-specific step is the second
projection: nonnegativity and touching define a certified elimination cone,
the repaired potential recovers a complete partial-minimum tower, and
\eqref{eq:twotaxrepair} separates the unavoidable integrability cost from
the additional MM certification cost.  On a general complex the harmonic
component of the Hodge decomposition is classical; the increment of
Theorem~\ref{thm:threetax} is that periods are priced in the same
target-loss units as curvature and certification.  The operator-semantic
result combines Minty's resolvent correspondence, the
Moreau--Gribonval--Nikolova characterization of proximity operators, and
the Beurling--Deny criterion
\citep{BauschkeMoffatWang2012,GribonvalNikolova2020,
FukushimaOshimaTakeda2011} to form the V passport.  Morse--Bott
nondegeneracy \citep{Bott1954} and Schur's lemma support the separate
quotient-leverage and selection results.  The Hilbert-space and
cellular-sheaf statements \citep{HansenGhrist2019} are extensions of the
quadratic filter law under declared domain and transport conditions, not
additional representation gates.  Gibbs-phase theory
\citep{Georgii2011} supplies the identifiability no-go; the corrector
variational principle \citep{JikovKozlovOleinik1994} and the classical
Polyakov--Nambu--Goto equivalence \citep{Polyakov1981} supply the two
closed-form frontiers.  In every case the claimed increment is the
certified combination, not the component.

Infimal projection and infimal postcomposition are classical
\citep{RockafellarWets1998,CombettesEtAl2019Interpolation}.
Theorem~\ref{thm:refinementgate} uses that operation as the transport law
for a certified lift--coupling package: equality after fiberwise descent is
both necessary and sufficient for universal frontier invariance, and it
composes across refinement chains.  Corollary~\ref{cor:duplicationnogo}
does not say that ensembles or model averaging are new.  It says that a
nonconstant strictly convex separable penalty on raw simplex coordinates
cannot be intrinsic when output-equivalent components may be duplicated.
This complements prior analysis of duplicate-model effects in weighted
model averaging \citep{GarthwaiteMubwandarikwa2010} by stating the
obstruction as uniform quotient faithfulness.

Parametric optimization studies value functions, sensitivity, and solution
maps under perturbation \citep{BonnansShapiro2000}.  The present parameter
\(\lambda\) acts on a variable that has already been eliminated, and the
vertical displacement of the value function is therefore identified
exactly:
\(
 \cJ_\lambda-\cJ=\inf_A(\cD_\Theta+\lambda\cR)
\).
The resulting frontier is a property of the lift, not merely a generic
regularization path.

Schur complements, spectral Tikhonov factors, and Gram classification are
classical \citep{ChungEspanolNguyen2014,HornJohnson2013MatrixAnalysis}, as
is the Le Chatelier principle that restricting an optimized response
stiffens the retained problem \citep{MilgromRoberts1996}.
Theorem~\ref{thm:completequadraticfingerprint} uses these ingredients to
make the augmented modal collection complete under origin-preserving
linear equivalence and to obtain its value--gradient--curvature exchange
inequality.  Proposition~\ref{prop:leveragefilter} then shows that the same
modal account governs target deformation and Schur-complement curvature
transmission.  Corollary~\ref{cor:contractionbody} gives the corresponding
inverse design statement before semantic admissibility is imposed.

Qualification and saturation are classical parts of inverse-problem
regularization \citep{EnglHankeNeubauer1996}.  Graph poly-Laplacian
regression can exploit higher smoothness under model-specific geometric and
sampling assumptions \citep{GarciaTrillosMurrayThorpe2023}, while
higher-order graph generators need not preserve positivity immediately
\citep{GregorioMugnolo2020}.  CRPS propriety and minimax distributional
regression rates are likewise established in their native settings
\citep{GneitingRaftery2007,PicDombryNaveauTaillardat2023}.
Bernstein subordination and functions preserving Stieltjes matrices are
also classical \citep{SchillingSongVondracek2012,MicchelliWilloughby1979}.
Section~\ref{sec:vsemantics} turns these ingredients into a semantic stress test.
Propositions~\ref{prop:dirichletcone}--\ref{prop:threepathqualification}
do not reprove those theories.  Their increment is the exact collision with
the V passport: on one three-node graph the admissible penalties commuting
with the Laplacian form a cone whose boundary is exactly modal order one,
so probability validity and higher Hilbert-scale approximation order are
alternatives with an explicit coupling threshold, and
Corollary~\ref{cor:generalgraphgate} transfers the sign half of that
statement to an arbitrary graph; its witness catalogue is restated in
Part~I \citep{FrameworkPartI} for the induced Richardson penalty.
Proposition~\ref{prop:weightedpathobstruction} then varies the same
three-node graph's edge weights to give an exact scale- and ratio-tunable
obstruction for an arbitrary spectral function.  The positive Bernstein
passport plus that reverse formula yields the all-weighted-graph
spectral-power classification and polynomial no-go; the claim is this
closure, not a new characterization of Bernstein or Stieltjes-preserving
functions.

These are not statistical-rate claims.  Their role in the series becomes
sharper when paired with Part I: the latter proves an exact
effective-dimension CRPS ledger and a matching minimax rate for every real
source smoothness \(s>0\), using a projected signed Richardson rule of any
integer qualification \(m\ge s\).  The claimed combination is therefore
not a new generic Sobolev exponent; it is a graph-universal validity
barrier followed by a nonlinear repair that is provably free at minimax
order.  The boundary is uniform validity, not a claim that projection is
active at every finite-sample optimal tuning.

Gaussian posterior means, Wiener shrinkage, and quadratic regularization are
classically linked \citep{KimeldorfWahba1970,Wahba1978}.
Theorem~\ref{thm:stochasticleverage} does not relabel that correspondence as
a new completeness theorem.  Its added content is conditional on the
certified quadratic auxiliary fiber: the declared lift fixes the defect
metric and admissible penalty cone, while deployment on a represented state
space is separately subject to the V passport.  The same filter that is
Bayes-optimal for the noisy oracle also prices the V target shift and the
outer leverage dividend.  Removing any of those conditions returns only
the classical Gaussian sequence calculation.

Bregman geometry supplies many possible defect coordinates and centroids
\citep{BanerjeeEtAl2005}.  Quantum information geometry and matrix
divergences supply particularly rigid ones
\citep{Petz1994CanonicalCorrelation,KulisSustikDhillon2009}.  Our use of
these geometries is representation-conditional: the certified lift fixes
the direction of divergence, and hence fixes the hard-coherence endpoint.
Likewise, the BKM--SLD ordering, quantum Cram\'er--Rao bounds, Petz
recoverability, and no broadcasting are classical
\citep{Hayashi2002QuantumFisher,Petz1986Sufficiency,
BarnumEtAl1996Broadcasting}.  Propositions~\ref{prop:quantumlocaltax} and
\ref{cor:quantumlowercertificate} and
Theorem~\ref{thm:quantumclassicalization} combine them only after elimination
has fixed Umegaki relative entropy as the required insertion currency.  The
claimed consequence is the certificate statement: a measured lower bound
cannot be substituted into an upper-budget gate, and a single lossless
classical ledger for the family exists exactly in the commuting case.

In policy optimization, \citet{ChanEtAl2022Greedification} show that reverse
KL has stronger improvement guarantees than forward KL, construct a
one-state counterexample in which reducing forward KL worsens a policy, and
derive conversions involving an additional R\'enyi-divergence term.
Theorem~\ref{thm:forwardklbarrier} does not claim the qualitative asymmetry
as new.  Its distinct statement is a dimension-uniform composition
impossibility at an exact teacher: forward KL and total variation both
vanish, while the exact zero-reward soft-return defect diverges.  The
framework contribution is to derive the required orientation from
elimination and to identify the missing exchange rate as a separate
currency rather than an untracked projection error.

Target-preserving representation design is classical in parameter-expanded
EM, conditional and marginal augmentation, and parameter-expanded data
augmentation
\citep{LiuRubinWu1998,MengVanDyk1999,LiuWu1999}.  Those theories select an
augmentation to improve the convergence or mixing of a single inferential
problem.  Section~\ref{sec:liftgate} asks a different question: which
properties of an exact lift survive auxiliary reparameterization or
redundant refinement when its coordinates are coupled across problems?
Bijective transport, fiberwise coherence descent, and the
unrestricted-collapse result separate intrinsic coupling claims from
coordinate acceleration or duplication.  Section~\ref{sec:fingerprint}
asks the distinct local question of what remains after a representation has
passed that gate.  Standard Schur-complement curvature is completed there by
oracle energy and its cross block, yielding a linearly complete modal
fingerprint and identifying the curvature retained by first-order
exactification.

General first-order surrogate theory studies models whose approximation
error has zero value and gradient at an anchor \citep{Mairal2013FirstOrder};
multifidelity trust-region model management already uses additive Taylor
correction to match low- and high-fidelity models
\citep{RobinsonEtAl2008Multifidelity}; block-successive upper-bound
minimization supplies convergence conditions for broad block approximations
\citep{RazaviyaynHongLuo2013}.  Thus neither the naked Taylor correction nor
first-order agreement is claimed as new.  Theorem
\ref{thm:xnormalform} is instead a certified-elimination converse: once a
frozen coupled auxiliary fiber is declared, its exact nonnegative insertion
defect determines every target-jet-preserving correction modulo a flat
remainder.  Proposition~\ref{prop:jetsplicing} identifies the retained
lift tail, Proposition~\ref{prop:xtransmission} determines whether a
conjugate auxiliary field can reach that tail at all, and
Proposition~\ref{prop:xbudget} exposes the remaining higher-order debit.

Variational constraints on latent laws are also well established
\citep{GanchevEtAl2010}.  Here their role is sharpened by comparison with
the constrained frontier and with external coupling at a matched coherence
budget.

The static C layer is classical probability coupling, and total correlation
is classical information theory \citep{Watanabe1960,CoverThomas2006}.  The
role of Theorem~\ref{thm:VCseparation} is architectural: it prevents a
dependence cost from being confused with a marginal auxiliary defect and
supplies the algebraic meaning of PGXV\(\times\)C.  When the joint objective
is linear in \(Q\), the penalized problem
\eqref{eq:dependencepenalized} is an entropically regularized
fixed-marginal transport problem \citep{Cuturi2013}; the role of
\(\mathfrak C_{P,c}\) here is not to rename that literature, but to complete
the coupling ledger with the least dependence compatible with a declared
joint benefit.  The entropy-projection structure underlying
Remark~\ref{rem:Cselfsimilar} is classical in the Schr\"odinger problem
\citep{Leonard2014}; its role here is to show that a fixed-marginal C design
can recursively generate its own certified auxiliary defect.

The distinction between pointwise legality and global representability also
recurs in Bochner, Pick, moment, marginal-consistency, complete-positivity,
and Markov-embeddability tests.  These analogies motivate the word
\emph{gate}; they are not additional results of the present calculus.
Remark~\ref{rem:admissibleL} records which certificates a declared
lift--coupling package must carry.

\section{Boundaries and extensions}
\label{sec:boundaries}

Several boundaries are structural rather than technical.  We first delimit
the finite-dimensional core, then record two transport-stable quadratic
extensions before returning to the remaining statistical and structural
limits.

\subsection{Scope of the finite-dimensional core}
\label{subsec:corescope}

First, the normal-form theorem is relative to Definition~\ref{def:onecoupling}.
Algorithms may couple outer parameters directly, introduce several nested
auxiliary layers, use non-scalar game dynamics, or alter the data-generating
model.  Such methods need not reduce to one letter.  The point of the grammar
is auditability: within a declared auxiliary-only coupling architecture, the
location and currency of nonlocality are explicit.

The composition law \eqref{eq:qdescentcomposition} is functoriality of
representation refinement only.  It does not solve composition of several
algorithmic coupling sites, exchange different defect currencies, or imply
that a statistically richer augmentation is computationally free.  Compact
Hausdorff spaces and lower semicontinuity in
Theorem~\ref{thm:refinementgate} are clean sufficient conditions for
fiberwise attainment; noncompact refinements require inf-compactness,
closure qualifications, or an explicitly approximate descent.

Theorem~\ref{thm:flatness} gives the gate for several nested auxiliary
layers.  If their local gaps are nonnegative, touching, and flat, the
canonical defects already form the exact multi-stage ledger
\eqref{eq:multilayerchain}, and several variational penalties have the
path-independent target-shift identity \eqref{eq:multilayervshift}.  If a method
instead splices independently designed surrogates, uses opposite divergence
orientations, or inserts an amortized projection, its defect curvature may
be nonzero and no common scalar ledger exists.  The Bregman three-point
remainder in Remark~\ref{rem:canonicalinteraction} is one such warning, not a
universal interaction tax.  Small curvature controls order ambiguity through
\eqref{eq:holonomybound}; it does not by itself imply convergence.

Proposition~\ref{prop:leveragefilter} is exact for a joint quadratic lift
with constant blocks.  In a general smooth model, the penalized auxiliary
minimizer moves with \(\lambda\), so \(B_L\), \(C_L\), \(K_L\), and their
spectral projectors must be reevaluated along the path.  Freezing them at
one anchor gives a local second-order approximation, not a global filtering
identity.
Likewise, a quantum transfer uses the BKM metric only on the tangent space of
faithful states and does not cover rank-deficient or zero-temperature
limits.  Finally, \(\cM_{L,\gamma}\) is profile curvature.  Calling it
statistical information requires a correctly oriented likelihood and a
separate score-variance calculation; penalized curvature alone does not
determine sandwich uncertainty.

\subsection{Infinite-dimensional and transported extensions}
\label{subsec:transportextensions}

Two controlled extensions remain inside the self-adjoint quadratic
certificate genus.  They retain the closed-form and domain-compatibility
conditions needed by the finite-dimensional filter law.

\begin{proposition}[The filter law with continuous spectrum]
\label{prop:hilbertfilter}
Let the whitened auxiliary space be a real Hilbert space, let
\(K\succeq0\) be the self-adjoint operator of a closed nonnegative
quadratic form, and set \(S_\lambda=(I+\lambda K)^{-1}\),
\(F_\lambda=I-S_\lambda\).  With the spectral measure
\(\nu_z(B)=\|E_K(B)z\|^2\), the value, defect, and roughness of the
penalized minimizer are one half times the integrals of
\(\lambda\gamma/(1+\lambda\gamma)\),
\(\{\lambda\gamma/(1+\lambda\gamma)\}^2\), and
\(\gamma/(1+\lambda\gamma)^2\) against \(\nu_{z^\star}\), exactly as in
Proposition~\ref{prop:couplingspectrum}.  For every \(\lambda>0\), the
roughness integral is finite because
\(\gamma/(1+\lambda\gamma)^2\le1/(4\lambda)\).  If finite oracle
roughness at \(\lambda=0\) is invoked, additionally require
\(z^\star\in\operatorname{dom}K^{1/2}\).  The curvature-transmission
identity is asserted when the whitened cross operator \(T\) is a bounded
operator between the declared auxiliary and outer Hilbert tangent spaces;
under that condition,
\(\nabla^2\cJ_\lambda-\nabla^2\cJ=TF_\lambda T^\ast\); and
\(F_\lambda\to I-P_{\ker K}\) strongly as \(\lambda\to\infty\).  Graph
limits, Neumann Laplacians, and Sobolev-type penalties therefore obey the
same filter law.  The hypotheses are not cosmetic: the form must be
closed, the defect metric uniformly coercive, and the cross operator
bounded as just stated (or supplied with a separate domain theorem that
makes \(TF_\lambda T^\ast\) a bounded quadratic form on the declared outer
domain); a merely pointwise positivity assumption does not suffice in
infinite dimension.
\end{proposition}

\par\medskip

\begin{proposition}[Bundle-valued fields, transport, and holonomy]
\label{prop:sheafV}
When auxiliary values live in different fibers across the index set, the
difference \(A_s-A_t\) is meaningless and transport is part of the
coupling semantics: the package becomes
\(L=(E\to\mathsf T,\cH_0,\cR,\nabla,C_{\mathrm{comp}})\), with
restriction or transport maps supplied by \(\nabla\).  For a quadratic
defect \(\tfrac12(A-a)^\top Q(A-a)\), \(Q\succ0\), and the coboundary
penalty \(\cR(A)=\tfrac12\|\mathsf BA\|_W^2\) of a cellular sheaf, the
product fiber Hilbert metric is part of the declared package.  All norms,
adjoints, isometry statements, and the strong-convexity constant \(\mu\)
below are taken in that same metric (with the \(Q\)-whitened coordinates
used for the modal statement).  The
penalized minimizer is
\begin{equation}
 A_\lambda=(Q+\lambda\mathsf B^\ast W\mathsf B)^{-1}Qa,
 \label{eq:sheafresolvent}
\end{equation}
the \(Q\)-whitened generalized sheaf-Laplacian modes obey the filter law
of Proposition~\ref{prop:leveragefilter}, and the hard endpoint is the
\(Q\)-orthogonal projection onto \(\ker\mathsf B\), the space of global
sections; ordinary graph consensus is the constant-sheaf case
\citep{HansenGhrist2019}.  For isometric transports \(U_i\) around a cycle
with holonomy \(U_C=U_{m-1}\cdots U_0\),
\begin{equation}
 \|(I-U_C)A_0\|^2
 \le m\sum_{i=0}^{m-1}\|A_{i+1}-U_iA_i\|^2,
 \label{eq:holonomybill}
\end{equation}
where \(A_m=A_0\),
so nontrivial holonomy forces a positive coherence bill on any field not
fixed by \(U_C\); and if the defect is \(\mu\)-strongly convex, the
hard-coherence price obeys
\(\cJ_\infty-\cJ\ge\tfrac{\mu}2
\operatorname{dist}^2(A^\star,\ker\mathsf B)\).
\end{proposition}

\par\medskip

Both extensions retain closed quadratic-form and domain-compatibility
assumptions; neither licenses an arbitrary nonlinear or noncoercive
infinite-dimensional coupling.

\subsection{Remaining statistical and structural boundaries}
\label{subsec:remainingboundaries}

The stochastic filter theorem has a deliberately narrower information set.
Its Bayes statement assumes a Gaussian oracle-field experiment, known
covariances, squared defect loss, and an effective precision inside the
declared penalty cone.  Overlapping windows, misspecification, robust scores,
or shared samples can produce nonwhite \(\Xi\); a rule using raw data or an
extra pilot is not dominated by a posterior mean that sees only
\(\widehat z\).  The theorem supplies an exact conditional risk ledger, not
a universal statistical optimality claim.

The probability--approximation-order frontier is equally scoped.
Propositions~\ref{prop:dirichletcone}--\ref{prop:threepathqualification}
concern penalties that commute with one declared Laplacian on three nodes;
Corollary~\ref{cor:generalgraphgate} extends the local sign test.
Propositions~\ref{prop:bernsteinpassport} and
\ref{prop:weightedpathobstruction} instead quantify over all finite
connected weighted graphs, but only within the spectral-function class;
they do not characterize arbitrary noncommuting penalties or every
graph-universally valid scalar function.  The spectral-power iff and
polynomial no-go are complete only in their stated families.  Nonlinear CDF
projection can restore validity, as Part I shows, but an active projection
mixes graph modes and therefore falls outside the exact quadratic
fingerprint.

Part I now supplies a complete fixed-design growing-path experiment:
independent replicated sampling, the path eigenvalue growth, a real-order
graph-source class, empirical-CDF noise, an effective-dimension upper bound,
and a matching Bernoulli lower bound.  Random geometric graphs, shared
samples across vertices, data-dependent tuning without sample splitting,
and adaptation to unknown source smoothness remain outside that theorem.
Its fixed-filter unbiased risk identity does not close the selection
problem.  The weighted-path results supply the semantic obstruction for
linear spectral mechanisms; they are not being repackaged as the minimax
lower bound, and they do not imply that repair is active at every
finite-sample optimal tuning.

The quantum observability results also do not provide a free measurement
algorithm.  Theorem~\ref{thm:quantumclassicalization} concerns equality of
raw outcome KL and Umegaki relative entropy under one fixed measurement.
Tomography or shadow reconstruction can use a richer output to estimate the
state, but an upper confidence bound for all candidate defects requires a
separate shot-complexity and simultaneous-inference analysis.  The
multi-parameter measurement frontier may additionally be governed by
Holevo-type incompatibility rather than the scalar SLD bound used in
Proposition~\ref{prop:quantumlocaltax}.

The frontier need not be convex as a function of \(r\).  Penalized V
recovers supported frontier points; constrained formulations may reveal
unsupported points.  Strong dual recovery in
Theorem~\ref{thm:frontier}(iii) therefore requires its stated convexity and
qualification assumptions.

Exactification needs the appropriate target or defect jet.  When
\(\cJ\) is itself intractable, an unbiased stochastic score, a dual bound,
or another certificate may replace exact evaluation only after a separate
error ledger is proved.  The phrase ``preserves the target'' is not licensed
by an approximate correction alone.  The conjugate transmission law is
also conditional on the canonical X model and an \(A\)-independent
stabilizer; an auxiliary-dependent trust term creates an additional
channel that must be audited separately.

The current calculus is built on one-sided minimum elimination.  A
convex--concave saddle has two profiled certificates and a primal--dual gap.
Its coupled theory should be built from a two-sided frontier.  Alternating
descent--ascent does not inherit scalar block descent merely by changing
\(\inf\) to \(\inf\sup\).

Gauge-cohomological and pushforward eliminations lie one step further
outside the present certificate genus.  A BV master action constrains its
deformations through a graded bracket and cohomology; supersymmetric
localization cancels by exactness rather than by a nonnegative ledger;
T-duality eliminates through a gauged parent action whose functional
determinant shifts the dilaton; and quantum worldsheet integration is a
pushforward, not a profiled infimum.  A future calculus of
``pushforward elimination preserving a structural equation'' would need a
different notion of certificate; relabeling such structures as minimum
defects is a category error, not a generalization.

Finally, objective fidelity is not statistical, computational, or physical
optimality.  The quantum endpoint and classicalization criterion are exact
statements, not claims of efficient large-scale computation.  Risk bounds,
mixing-time results, shot complexity, and physical interpretations require
assumptions beyond the present calculus.

\section{Discussion}
\label{sec:discussion}

Part I takes a declared certified lift as an explicit interface and develops
the five-mechanism constructive calculus.  This Part II has audited that
interface rather than treated its existence or interpretation as automatic.
The deterministic audit proceeds from integrability, through representation
validity, to a scalar frontier and a non-scalar leverage fingerprint;
variational filters then face a separate semantic gate.  The fixed-marginal
layer adds an orthogonal horizontal dependence frontier.
The integrability gate decides whether local certificates belong to any
common scalar objective: nonnegativity and touching certify the edges, while
zero defect curvature glues them into a partial-elimination tower.  The
finite-state repair theorem makes this gate constructive.  It lower-bounds
the error of every scalar explanation by curvature energy and reconstructs
the nearest certified tower, with a separate debit for restoring
majorization and touching.  The
representation gate then removes coordinate aliases, descends coherence
through redundant fibers, and excludes target-calling or dummy lifts from
claiming free coherence.  Its duplication corollary gives a concrete audit
for model averaging and ensembles: copying an output-equivalent component
cannot legitimately change the coherence price.
Conditional on an admissible package, the primitive is the defect
\(\cD_\Theta\), the exact height of an auxiliary state above its eliminated
target.  Constraining coherence produces \(\mathfrak F_\Theta\); penalizing
coherence produces its supported dual.  Every externally selected field can
then be located vertically above the frontier rather than judged by an
operator-specific criterion.

The marked collection
\[
 \{(\gamma,m_\gamma,\mathcal I_\gamma)\}_\gamma
\]
is a quadratic fingerprint complete under origin-preserving linear
auxiliary equivalence.  It adds the non-scalar diagnostic.
Its blocks separate oracle energy \(w_\gamma\), outer leverage
\(\cM_{L,\gamma}\), and their missing orientation channel \(b_\gamma\).
The last block matters: energy and curvature alone do not determine the
first-order target response.  The PSD identity
\eqref{eq:fingerprintexchangeblock} is a genuine exchange condition, not
just parallel bookkeeping; it bounds the target-shift gradient by the value
paid and curvature transmitted.  Within the same quadratic germ, the
leverage filter law makes the transmission exact, while the zero mode is
conserved.  Under an aligned Gaussian experiment, the same filter also gives
optimal Bayes bias--variance allocation.  A lift with a small coupling
defect but negligible leverage may simply have hidden the structure that
coupling was meant to exploit.  Lift choice must therefore trade fidelity
and coherence against leverage and computation; no scalar ``best
representation'' exists without an application-specific admissible class.

The three-node path supplies the local semantic stress test for V, and the
test is a cone rather than a list of examples.  Varying its two edge weights
then turns that local picture into a graph-universal spectral witness.  The
positive side is the classical Bernstein passport: subordination converts a
local graph into a possibly dense long-range Laplacian.  The exact reverse
formula fires whenever \(\psi(rx)>r\psi(x)\).  Together they show that
spectral powers are valid on every weighted graph exactly up to order one
and that no finite-degree graph-universal polynomial can exceed the ordinary
Laplacian.  This is a validity--locality--order theorem for stated linear
spectral families, not a universal impossibility theorem and not, by itself,
a statistical lower bound.

Part I's signed Richardson smoother followed by CDF projection is a
constructive escape from this linear family: it restores probability
validity while retaining an exact finite-sample effective-dimension bound.
On growing paths it attains the matching minimax CRPS rate for every real
source smoothness \(s>0\).  At \(s\le1\) the Markov resolvent is already
enough; at \(s>1\) nonlinear repair permits higher qualification without
losing rate.  The price is representational rather than statistical: once
projection is active, the estimator leaves the exact quadratic fingerprint
and no raw modal interpretation may be carried through it.  Uniform
necessity of the interface does not mean that projection activates at every
finite-sample tuning.

The normal-form theorem identifies where coupling enters.  External
insertion selects a generally inefficient point; variational coupling traces
supported frontier points; exactification splices the target's low-order jet
to the frozen lift's higher-order tail; and fixed-marginal coupling controls
joint realization after coordinate laws have been fixed.  Raising the X
order deletes additional auxiliary channels.  In a conjugate lift, the map
\(\delta A\mapsto-\langle\delta A,D^2\eta\rangle\) decides whether a
first-order X field changes curvature at all; affine \(\eta\) is the
automatic-safety/no-leverage endpoint.  The horizontal frontier
\(\mathfrak C_{P,c}\) additionally asks for the least dependence needed to
produce a declared joint gain.  Entropic C design is recursively certified
by \eqref{eq:Cselfsimilar}: a suboptimal joint realization pays an exact
relative-entropy defect while its coordinate marginals remain protected.
The soft-control barrier shows that this ledger remains necessary even when
each stage is individually exact.  It catches a direction error: the
projection reports \(\KL(\textnormal{teacher}\|\textnormal{actor})\), while
elimination certifies
\(\KL(\textnormal{actor}\|\textnormal{teacher})\).  The quantum test catches
a different proof error: measurement supplies a lower bound on the
Umegaki defect, while monotone acceptance needs an upper bound.
Theorem~\ref{thm:quantumclassicalization} shows that one fixed measurement
can erase this second ledger only for a commuting family.  These are not
decorative applications; they are two places where an otherwise plausible
performance or descent proof has the inequality in the wrong currency or
the wrong direction.

The three pillars now separate cleanly.  Flatness decides whether reported
local gaps are integrable and repairable.  Within the second pillar,
quotient faithfulness decides whether a coupling comparison is intrinsic;
the frontier prices its scalar deformation; the modal fingerprint records
how auxiliary modes reach the outer target; and the V passport tests whether
the induced operator respects the represented state space.  The duplication
and weighted-path no-go results expose failure at the first and last of
these stages.  The three hard tests then show what cannot be exchanged
silently, while approximate X closes the algorithmic loop: a computed
defect jet comes with an exact score-error identity, a residual-transport
acceptance test, and summability conditions leading to stationarity.

There is deliberately no universal scalar conservation law.  Insertion
defect and target shift are measured in objective units; dependence is
measured by a joint-law divergence; evaluation, sampling, and computation
have still other units.  Theorem~\ref{thm:conservation} says only that a
nontrivial coupling leaves a record in at least one native account.  A
support condition in soft control or a measurement remainder bound in the
quantum test is an example of the additional exchange theorem required
before two accounts can be compared.  With several coupling sites,
flatness is prior even to that question: without it, the reported objective
payments depend on the order in which the local gaps are read.

Beyond the declared class, payments need not add.  Nested auxiliary and
projection stages may have nonzero defect curvature.  A composition theory
for several coupling sites therefore requires not another flat mechanism
label, but a distinction between canonical conditional defects, which are
path independent, and externally reported gaps, whose holonomy must be
audited.  The refinement composition law
\eqref{eq:qdescentcomposition} does not remove this open boundary: it
composes representations, not algorithmic coupling payments.

The series therefore proceeds in the natural reading direction.  Part I
moves forward from a declared certified lift to algorithms, certificates, and
statistical constructions.  Part II then moves backward from that interface
to its structural preconditions and boundaries: integrability and repair,
representation validity, frontier and leverage geometry, semantic
admissibility, forced converses, and impossibility tests.  The second
direction limits the first without being a prerequisite for reading it.

\bibliographystyle{plainnat}
\bibliography{references}

\end{document}